\documentclass[a4paper,12pt]{article}

\usepackage{amsmath,amsfonts,graphicx}

\usepackage{float}

\numberwithin{equation}{section}
\numberwithin{equation}{subsection}
\numberwithin{table}{section}

\title{Elliptic Curves in Recreational Number Theory }
\author{Allan J. MacLeod,\\Statistics, O.R. and Mathematics Group, (Retired)\\
University of the West of Scotland,\\High St., Paisley,\\Scotland.  PA1 2BE\\
(e-mail: peediejenn@hotmail.com)}

\begin{document}

\maketitle

\begin{abstract}
Several problems which could be thought of as belonging to recreational mathematics are described. They are all such that solutions to the
problem depend on finding rational points on elliptic curves. Many of the problems considered lead to the search for points of very large height on the
curves, which (as yet) have not been found.
\end{abstract}

\newpage

\section{Introduction}
Elliptic curves underpin some of the most advanced mathematics currently being pursued, most notably (probably) the proof
of Fermat's Last Theorem by Sir Andrew Wiles. Recreational number theory is enjoyed by many people with {\bf NO} professional
mathematics qualifications. The purpose of the present survey is to show how, even a small knowledge of elliptic curves and
a computer, can help solve many problems.

These problems are all reasonably simple to understand but can lead to some very interesting mathematics. We start with a simple

{\bf Problem:}  If possible, find rational numbers $a,b,c$ with
\begin{equation}
a\,b\,c = N = a+b+c
\end{equation}
where $N$ is a specified non-zero integer.

The first stage is to analyze the problem. We will find that solutions are related to rational points on the elliptic curve
\begin{equation}
y^2=x^3+N^2x^2+8N^2x+16N^2=x^3+N^2(x+4)^2
\end{equation}

The rational points form a group, often denoted $\Gamma$, which is finitely generated. This means that every element $P \in \Gamma$
can be written
\begin{equation}
P=n_1P_1+\ldots+n_rP_r+T \hspace{2cm} n_i \in \mathbb{Z}
\end{equation}
where $P_1,\ldots,P_r$ are points of infinite order, called generators, and $T$ is a point of finite order, called a torsion point.

The subset of points of finite order is small and
easily determined. For example, the above curve has $(0,\pm 4N)$ as elements of $T$, together with the point at infinity which is
the group identity. The torsion points often correspond to trivial solutions of the problem, and only rarely give a solution we seek.

The quantity $r$ is called the {\bf rank} of the elliptic curve. It can be zero, so that there are only torsion points, which usually means
no non-trivial solutions.

We can estimate the rank using the famous Birch and Swinnerton-Dyer (BSD) Conjecture and a moderate amount of computing. If the rank is
computed to be greater than one, it is often reasonably easy to find at least one generator, though there are exceptions to this rule.
The real computational problems are with
curves where the rank is estimated to be exactly $1$.

A by-product of the BSD calculations is an estimate of the height of the point - the higher
the height, the more digits in the rational coordinates and so the point is harder to find (in general). There are
two possible height normalizations available and I use the one described in Silverman \cite{silv}, which gives values
half that of Pari's current {\bf ellheight} command.

Let me state at the outset that my primary interest is in computing {\bf actual} numerical values for solutions of the many
problems discussed rather than just proving a solution exists.

To get anywhere in this subject, you need a computer and some software. To help with the algebra in the analysis of the problems, some
form of Symbolic Algebra package would be very useful. I tend to use an old MS-Dos version of Derive but this falls down sometimes, so I have
been known to use Maxima (free) or even Mathematica (not-free).

More importantly, you need a package which does the number theory calculations on elliptic curves. The most powerful is Magma, but this
costs money. Pari is a free alternative but does not have so many built-in capabilities. The vast bulk of my code is written using
Pari, since I only recently purchased Magma as a retirement present to myself!
If you have an old 32-bit machine, I would recommend UBASIC for development work as it is very fast.

I have tried to give the correct credit where appropriate. If you think I have failed in this task, please contact me. I have
also tried to make each problem description self-contained so there is quite a lot of repetition between sub-sections. Any mistakes are,
of course, purely my fault. I would be extremely pleased to hear from readers with problems I have not covered.

\section{Elliptic Curves $101$}

What is an {\bf elliptic curve}?

Since I am aiming to make most of this report understandable to non-professionals, I will take a very simplistic approach,
which will certainly appal some professionals. But you can get quite far in using elliptic curves, without worrying about
topics such as cohomology or Galois representations.

There are several books that can be recommended on elliptic curves. At an introductory level Silverman and Tate \cite{siltate} is excellent.
More advanced are the books by Husem\"{o}ller \cite{husem}, Knapp \cite{kn} and Cassels \cite{cassbk}.
The two volumes by Silverman \cite{silvbk1} \cite{silvbk2} have become the standard mathematical introduction.

Consider the simple cubic
\begin{equation}
f(x)=x(x-2)(x-5)=x^3-7x^2+10x
\end{equation}
which has zeros at $x=0,2,5$. These are the only zeros, so the curve is non-negative for $x \ge 5$ and $0 \le x \le 2$, but
strictly negative for $2 < x < 5$ and $x<0$.

Make the very simple change of $y$ to $y^2$ giving
\begin{equation}
y^2=x(x-2)(x-5)=x^3-7x^2+10x=f(x)
\end{equation}
which gives the curve in the following graph.

\includegraphics{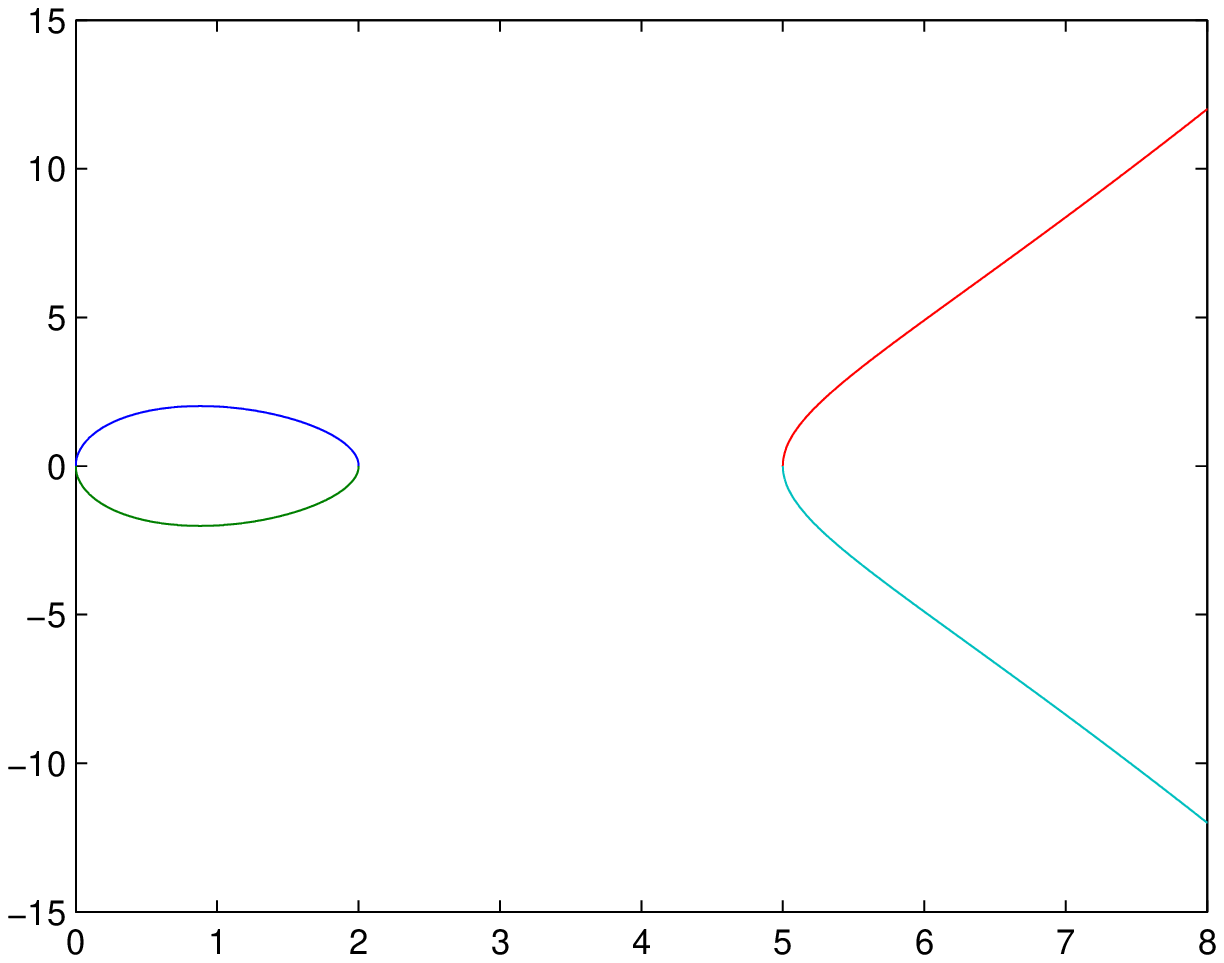}

Firstly, this curve is symmetric about the x-axis, because if $(u,v)$ lies on the curve so does $(u,-v)$. Secondly, there are
no points for $2 < x < 5$ or $x<0$, since $f(x)<0$ in these intervals. This means that the curve has two disconnected components.
The first for $0 \le x \le 2$ is a closed shape, usually called the "egg" - for obvious pictorial reasons. The second for $x \ge 5$
is usually called the infinite component.

This is all very nice and pretty, but it does not help us solve our problems. The problems we consider look for integer or rational
solutions, and we can link these to finding rational points on the elliptic curve. A rational point is one where both coordinates
are rational. There are $3$ obvious rational points on this curve, namely $(0,0)$, $(2,0)$, and $(5,0)$. A little simple search
shows that $x=8$ gives $y^2=8 \times 6 \times 3=144$, so $(8,12)$ and $(8,-12)$ are also rational points on the curve.

Let $A=(8,12)$ and $B=(2,0)$. The line joining $A$ to $B$ has equation $y=2x-4$, and meets the curve where
\begin{equation}
y^2=(2x-4)^2=x^3-7x^2+10x
\end{equation}
which gives the cubic equation
\begin{equation*}
x^3-11x^2+26x-16=0
\end{equation*}

This has at most $3$ real solutions, and we already know two of them, $x=2$ and $x=8$. We also know from basic algebra that the product of the
three roots is $-(-16)$. Thus the third point of intersection must be $x=1$ which gives $y=-2$. So from two rational points we get a third
rational point!

Thinking about this procudure, if $A$ and $B$ are any rational points, the line joining them will have rational coefficients. So the
intersection with the curve will be a cubic with rational coefficients. The product of the three roots will be the negative of the constant
in the cubic and so rational. Since we know two of the roots are rational, the third must be rational. But the line has rational coefficients,
so the y-coordinate must also be rational. So from two rational points we get a third rational point. This operation is thus closed on the
set of rational points $\Gamma$.

A closed binary operation is the first requirement for a GROUP. With some further fiddling, we can define an "addition" operation on $\Gamma$,
which gives a commutative group. This is discussed further in the appendices.

The crucial point is that this group is finitely generated. This means that there are points $G_1,\ldots,G_r \in \Gamma$ such that
\begin{equation}
P=n_1G_1+\ldots+n_rG_r+T
\end{equation}
where $P$ is any rational point, $n_i$ are integers and $T$ is an element of a small subset of $\Gamma$ which is called the torsion subgroup
and which can be easily computed.

Thus, all we need to do to find rational points is to find $G_1,\ldots,G_r$ and hence use these to find a solution to our original problem.
The value of $r$ is called the {\bf rank} of the elliptic curve, and it can be $0$. For example, $y^2=x(x-2)(x-4)$ has rank zero.

The challenge is to compute these generator points. It can be very difficult, but that is what provides the challenge.

\newpage

\section{Introductory Problem}
Let $x=1$, $y=2$ and $z=3$. Then $x+y+z=6$ and $x\,y\,z=6$. It is a natural question to ask if this happens for other integers apart from $6$.
\begin{equation}
x+y+z=N=x\,y\,z
\end{equation}

The following lemma shows, however, that we are doomed to failure.

{\bf Lemma.} The only integer solution is the one above.

{\bf Proof.} We can assume, without loss of generality, that $N>0$ and so, either all of $x,y,z$ are positive or two of them are negative and one positive.

Suppose, firstly, that all the variables are positive and we can assume $1 \le x \le y \le z$. We have
\begin{equation}\label{espz}
z=\frac{x+y}{x\,y-1}
\end{equation}
and
\begin{equation*}
x+y-(x\,y-1)=(1-x)(y-1)+2
\end{equation*}

If the right-hand-side is negative then $0 < z < 1$, contradicting $z$ being an integer. If $x\ge2$, the only value of $y$ giving a positive RHS
is $y=2$, but this gives $z=4/3$. $x=2,y=3$ gives $z=1$ and $N=6$. If $x=1$, $z=(y+1)/(y-1)$, which only gives integer values when $y=2, z=3$,
again for $N=6$.

If we now assume two of $x,y,z$ are negative, let them be $x$ and $y$, with $z>0$. But \eqref{espz} implies $z < 0$, a contradiction.

Thus, the only purely integer solution is the one at the start of the section.

\vspace{1cm}

We now consider $x,y,z \in \mathbb{Q}$.

We have $z=N/(x\,y)$, and, substituting, we have that $x,y$ must satisfy the quadratic equation
\begin{equation}\label{eqspquad}
x\,y^2+x(x-N)y+N=0
\end{equation}

For fixed $N$, we can generate rational values of $x$, and test if this quadratic factors into linear factors.
For $1 \le N \le 19$, the following results were obtained in a simple search

\begin{table}[H]
\begin{center}
\caption{Equal sum and product results}
\begin{tabular}{lrrr}
$\,$&$\,$&$\,$&$\,$\\
N&x&y&z\\
$\,$&$\,$&$\,$&$\,$\\
6&1&2&3\\
7&$4/3$&$9/2$&$7/6$\\
9&$1/2$&$4$&$9/2$\\
13&$36/77$&$121/42$&$637/66$\\
14&$1/3$&$9$&$14/3$\\
15&$1/2$&$12$&$5/2$\\
16&$-2/3$&$18$&$-4/3$\\
19&$121/234$&$324/143$&$3211/198$
\end{tabular}
\end{center}
\end{table}

For those values of $N$ without a solution there are $2$ possibilities:
\begin{itemize}
\item A solution exists but we haven't searched far enough to find it.
\item A solution doesn't exist at all.
\end{itemize}
Elliptic curves allow us to find which alternative is true for a value of $N$, and also to attempt to find a solution if one does exist.

For the quadratic in $y$ \eqref{eqspquad} to have rational solutions, the discriminant must be a rational square. Thus
there must exist $D \in \mathbb{Q}$ such that
\begin{equation}
D^2=x^2(x-N)^2-4Nx=x^4-2Nx^3+N^2x^2-4Nx
\end{equation}

Since $x \ne 0$, we have
\begin{equation}
\frac{D^2}{x^4}=-4N\frac{1}{x^3}+N^2\frac{1}{x^2}-2N\frac{1}{x}+1
\end{equation}

Multiply both sides by $(-4N)^2$ and define $G=4ND/x^2$, $H=-4N/x$ to give
\begin{equation}\label{eceqsp}
G^2=H^3+N^2H^2+8N^2H+16N^2=H^3+(NH+4N)^2
\end{equation}
which is the elliptic curve for this problem. Since $x,D$ are rational, $(H,G)$ must be a rational point on the curve and vice-versa.

The curve has the rational points $(0,\pm 4N)$, but these do not give finite solutions since $x=-4N/H$.
Thus, we need to know whether there are other rational points on the curve for specific values of $N$.

The points $(0, \pm 4N)$ are points of inflexion on the curve and are thus of order $3$. This reduces the possibilities for the torsion
subgroup $T$ to be isomorphic to $\mathbb{Z}/3\mathbb{Z}$, $\mathbb{Z}/6\mathbb{Z}$, $\mathbb{Z}/9\mathbb{Z}$, $\mathbb{Z}/12\mathbb{Z}$, or
$\mathbb{Z}/2\mathbb{Z} \oplus \mathbb{Z}/6\mathbb{Z}$. Numerical tests suggest that only the first possibility occurs, but this is,
of course, nothing like a proof.

The curves have discriminant
\begin{equation*}
\Delta=4096N^4(N^2-27)
\end{equation*}
so are non-singular for all $N$ in $[1,999]$. We have $\Delta >0$ for $N \ge 6$, so, in this case, the curves have two components - an "egg" component
and an infinite component to the right of the egg.

In general, therefore, we need the elliptic curve to have rank greater than $0$ for a solution to the problem, but, before discussing the search for points,
we consider the possibility of parametric solutions. The search method is extremely simplistic. We generate simple polynomial forms for $N$ and $H$ in
\eqref{eceqsp} and test whether they give a square, ignoring those $H$ which give a torsion point.

A very limited search using Pari-GP gives $N=2k^2-2$ and $H=8k(k+1)$ which gives
\begin{equation*}
(x,y,z) = \left(\, \, \frac{1-k}{k} \, \, , \, \, 2k^2 \, \, , \, \, \frac{-(k+1)}{k} \, \, \right)
\end{equation*}

There is no known
guaranteed method for determining the rank, so we used the Birch and Swinnerton-Dyer conjecture to estimate the rank of the curves.
If we throw enough computer time at the calculations, we can be fairly certain that the results are correct.

Applying the BSD calculations for this problem, with $N \in [1,999]$, we find $374$ curves with rank zero, $514$ curves with
rank one, and $111$ curves with rank greater than one.

A by-product of the rank-one computations is an estimate of the height
of the rational point. Very simply, the greater the height the more digits are in the numerator and denominator of $H$. There is a
drawback, however, in that the height calculation involves a term for the size of the Tate-Safarevic group, which we don't
know how to calculate, so we give it the value $1$. This can sometimes lead to huge over-estimates of the true height of a point.

It is usually assumed that at least one point on a curve of higher rank can be easily calculated, though there are counter-examples,
such as $N=733$ in the current problem.

For this problem, the smallest estimated rank-1 height is $0.194$ with $N=6$, whilst the largest is $549.3$ for $N=977$.
Readers should note that I use the height normalisation which is the smaller of the two available, so, for example, this is half that
given by more recent versions of Pari.

These rank calculations allow us to concentrate on the curves where a point can be found, and not waste time on curves where
$(0, \pm 4N)$ are the only points.
To compute rational points, we use a whole arsenal of approaches.

For small heights, we search for solutions where $H$ might be
an integer or a rational $H=u/v^2$ with $u,v$ small integers. This usually gives most heights up to about $5$, and several of the higher rank
curves. Searching eventually becomes time-consuming, so we need a more sophisticated approach.

The elliptic curves \eqref{eceqsp} do not have points of order $2$ so cannot be expressed in the form $y^2=x^3+ax^2+bx$. Such curves
are much easier to deal with. For curves without 2-torsion, John Cremona's {\bf mwrank} package still works but can take some time.
I am a great fan of the method described by Joseph Silverman in \cite{silv}.
I have both a Ubasic and Pari implementation which works efficiently up to heights of $16$ and $17$ and can be stretched to about $20$
if the curve has a moderately-sized conductor and a smallish number of local heights.

This method finds a large number of points before the heights become too large. We are helped, however, by one of the most useful
theoretical ideas - an isogenous curve. Very very roughly, the torsion points allow us to consider a related elliptic curve,
which has the same rank and a computable relationship between rational points. It is often the case that a rank-$1$ isogenous
curve has a point of smaller height, which can be computed fairly easily. There is a discussion of isogeny in the Appendix.

For this problem, the $3$-isogeny curve from \eqref{eceqsp} is
\begin{equation}
V^2=U^3-27N^2(U-4(N^2-27))^2
\end{equation}
which only has the point at infinity as a torsion point.

We used this to find a further large group of rational points on the isogenous curve, which are easily transformed
to points on the original curve. There were still a considerable number of unsolved values left and I moved on to
the method of 4-descent originally described by Merriman, Siksek and Smart \cite{mss} and explained very nicely in
Tom Womack's Ph.D thesis \cite{wom}. Magma contains software for this procedure, but Magma is not free. I have developed
my own Pari code for doing (most of) the computations and the code works very well most of the time.

With this, I reduced the number of unsolved values to less than $50$. For most of the other problems discussed in this report,
I would be able to use the method of Heegner points, see Watkins \cite{watk}, for some points with very large height but smallish conductor.
For these curves, however, the conductors are large and Heegner point methods are essentially ruled out.

As of the time of writing, the following Table gives the unsolved values of $N$. I bought myself Magma as a retirement present
and have used it to resolve some problem values.

\begin{table}[H]
\begin{center}
\caption{Equal sum and product unsolved}
\begin{tabular}{lrr}
$\,$&$\,$&$\,$\\
N&Height&3-isog ht.\\
$\,$&$\,$&$\,$\\
773&404&135\\
802&193&579\\
895&131&393\\
947&548&183\\
956&150&451\\
977&549&183
\end{tabular}
\end{center}
\end{table}

A natural further question to ask is: for those values of $N$ with a solution, can we find one with $x,y,z > 0$?

From $x=-4N/H$, it is
clear that $H<0$ is a necessary requirement. We also require $0<x<N$. From the equation of the curve, the gradient is
\begin{equation*}
\frac{dG}{dH}=\frac{3H^2+2N^2H+8N^2}{2G}
\end{equation*}
which is greater than zero if $H,G \ge 0$. Thus $(0, \pm 4N)$ must lie on the infinite component, and all $3$ zeroes
of the elliptic curve are negative. Since $G^2=-64$ for $H=-4$, at least one of the zeroes must lie in $(-4,0)$.
But, if $-4 <u <0$ we have $x>N$, so an all positive solution can only come from a point on an egg component. Thus, there
are no positive solutions if $N<6$. From now on, we assume $N \ge 6$.

If $H=-2N$, then $G^2=4N^2(N^2-6N+4)$ which is positive if $N \ge 6$. Thus there is at least one point (not necessarily rational)
on the egg with $H<-4$. Thus the whole egg must exist in the region $H<-4$, where $0<x<N$, so we have a positive rational solution from
any rational point on the egg.

This means that at least one generator must lie on the egg to give a positive solution. It is impossible for the chord-and-tangent
process to get a point on the egg from two points on the infinite component, since the egg is closed and a line going into the egg
must come out of the egg, giving $2$ intersections.

\newpage

\section{Problems from Triangles}

\subsection{Congruent Numbers}
Consider a right-angled triangle

\begin{center}
\begin{picture}(10,8)
\put(0,0){\line(1,0){12}}
\put(12,0){\line(0,1){8}}
\put(0,0){\line(3,2){12}}
\put(12,0.5){\line(-1,0){0.5}}
\put(11.5,0){\line(0,1){0.5}}
\put(6,0.25){b}
\put(12.25,4){a}
\put(6,4.5){h}
\end{picture}
\end{center}

By Pythagoras, $a^2+b^2=h^2$ and the area $\Delta=ab/2$. The most well-known such triangle is probably the $(3,4,5)$ triangle
introduced at school. It has area $6$. A natural question is {\bf what other integers can occur as the area of a right-angled
triangle with integer sides?}

If we scale the $(3,4,5)$ sides by $2$, we get $(6,8,10)$ with $6^2+8^2=10^2$ and area $24$. Scale by $3$ and the area goes up by
factor of $9$ to $54$. The next smallest integer right-angled triangle is $(5,12,13)$ with area $30$. So restricting to integer
sides is, perhaps, too restrictive.

What if we allow $a,b,c$ to be rational? Let $\Delta=5$ and $a=A/D$, $b=B/D$ and $h=H/D$ with $A,B,H,D \in \mathbb{Z}$. Then
\begin{equation*}
A\,B=10D^2 \hspace{2cm} A^2+B^2=H^2
\end{equation*}
so that we require
\begin{equation*}
A^4+100D^4=A^2H^2
\end{equation*}

Searching for $(A,D)$ pairs which give $A^4+100D^4$ equal to a square, which is divisible by $A$, quickly gives $A=9$ and $D=6$, giving
$a=3/2, b=20/3, h=41/6$. If we miss out the requirement that the square is divisible by $A$, we get $A=3, D=2, B=40/3$ quicker, and
the same triangle eventually. It is simple to check $a^2+b^2=h^2$ and $ab/2=5$.

If, however, we try to do the same procedure for $\Delta=3$ we will be waiting forever - there is no solution. There is thus a dual problem.
For $\Delta=N$, does a solution exist? If it does, what is it?

We now show the link between the triangle problem and elliptic curves. There are many approaches, and I particularly like this one, from
Keith Conrad \cite{con}.

Let $a,b,h > 0$ and $N>0$. Then $t=h-a>0$, and
\begin{equation*}
2\,a\,t=b^2-t^2
\end{equation*}
so that
\begin{equation*}
\frac{4Nt}{b}=b^2-t^2
\end{equation*}

Multiply by $N^3b/t^3$ to give
\begin{equation*}
\frac{4N^4}{t^2}=\frac{N^3b^3}{t^3}-N^2\frac{Nb}{t}
\end{equation*}
and define $Y=2N^2/t$ and $X=Nb/t$ so that
\begin{equation}\label{econg}
Y^2=X^3-N^2X=X(X-N)(X+N)
\end{equation}

Conversely, let $(u,\pm v)$ be rational points on the elliptic curve with $|v|>0$. Define $a=v/u$ choosing the sign
of $v$ to make $a>0$. Define $b=2N/a$, so $a,b$ are strictly positive with $N=a\,b/2$ and
\begin{equation*}
a^2+b^2=\frac{v^2}{u^2}+\frac{4N^2u^2}{v^2}=\frac{v^4+4N^2u^4}{u^2v^2}
\end{equation*}

Using $v^2=u^3-N^2u$, we have
\begin{equation*}
a^2+b^2=\frac{(u^3-Nu)^2+4N^2u^4}{u^2v^2}=\frac{(u^3+N^2u)^2}{u^2v^2}=\left( \frac{u^2+N^2}{v} \right)^2
\end{equation*}
so we can set $h=|(u^2+N^2)/v|$.

The elliptic curve \eqref{econg} has $3$ finite torsion points when $X=0,N,-N$.
By the area scaling mentioned earlier, we can assume that $N$ is squarefree.

Since doubling a point on a curve of this form gives an X-coordinate which is a rational square, a point of order $4$ would have to solve
\begin{equation*}
\frac{(X^2+N^2)^2}{4(X^3-N^2X)}=0
\end{equation*}
which is clearly impossible. A point of order $3$ would be a point of inflexion so would satisfy
\begin{equation*}
3x^4-6N^2x^2-N^4=0
\end{equation*}
which factors to
\begin{equation*}
3(x^2-(\frac{2\sqrt{3}}{3}+1)n^2)(x^2+(\frac{2\sqrt{3}}{3}-1)n^2)
\end{equation*}
showing there are no rational roots. The lack of points of orders $3$ and $4$ shows that the torsion subgroup must be isomorphic to
$\mathbb{Z}/2\mathbb{Z} \oplus \mathbb{Z}/2\mathbb{Z}$. None of the torsion points provide a triangle with area $N$, so we must find
points of infinite order.

The rank of the congruent number elliptic curve has been the subject of a great deal of work over the years. There is no way I can
cover all of this, especially as some of it is very advanced. If you are interested in further reading, just put "congruent number"
into a search-engine.

In \cite{steph}, Nelson Stephens
showed that a form of the Birch and Swinnerton-Dyer Conjecture implies that $N$ is a congruent number if $N$ is squarefree and congruent
to $5,6,7$ modulo $8$. If $N$ is squarefree and congruent to $1,2,3$ modulo 8, then most $N$ are not congruent numbers but some are.
Tunnell in \cite{tunn} provided a simple criterion for determining whether $N$ is a congruent number, involving computations with modular forms.
A very readable summary of all this is in the book by Neal Koblitz \cite{kob1}.

Knowing a number $N$ is congruent can still be a long way from finding a specimen triangle. Over the years, many computational approaches
have been suggested for determining solutions. Examples can be found by Alter and Curtz \cite{altcur} and Godwin \cite{godw} - I particularly
like Godwin's method.
For larger solutions, Noam Elkies' \cite{elk1} showed how to apply essentially Heegner-point methods
to the elliptic curve. The calculations of all solutions for the first one million congruent numbers has been completed,
according to a recent e-mail from Randall Rathbun, the principal searcher for several years.

For those who like pencil-and-paper algebra, this problem provides almost unlimited opportunities.
We use the fact that, for elliptic curves of the form
$y^2=x^3+Ax^2+Bx$, with $A,B \in \mathbb{Z}$, a rational point $(x,y)$ has the structure $x=du^2/v^2$ and $y=duw/v^3$, where
$d,u,v,w \in \mathbb{Z}$, $d$ is squarefree and $\gcd(u,v)=\gcd(d,v)=1$.

Thus
\begin{equation*}
d\,w^2=d^2u^4-N^2v^4
\end{equation*}
implying that $d|N$.

Suppose $d=-1$, so
\begin{equation*}
w^2=N^2v^4-u^4=(Nv^2+u^2)(Nv^2-u^2)
\end{equation*}
and suppose $Nv^2+u^2=r^2$ and $Nv^2-u^2=s^2$.

Thus
\begin{equation*}
Nv^2=r^2-u^2=(r+u)(r-u)
\end{equation*}
and suppose $r+u=Np^2$ and $r-u=q^2$. Thus, we have
\begin{equation*}
r=\frac{Np^2+q^2}{2} \hspace{1.5cm} u=\frac{Np^2-q^2}{2} \hspace{1.5cm} v=pq
\end{equation*}
which gives the quartic equation
\begin{equation}
-N^2p^4+6Np^2q^2-q^4=4s^2
\end{equation}

It is straightforward to program a search on this quartic which, very quickly, finds a solution
with $p=5, \, q=89, \, s=1361$, leading eventually to
\begin{equation*}
a=\frac{6428003}{1423110} \hspace{2cm} b=\frac{173619420}{6428003}
\end{equation*}
showing that this algebraic "descent" procedure can be very effective. Note that, because of all the suppositions, there
is no guarantee that this method will always work.

\subsection{Variants of the Congruent Number Problem}
We will describe two extensions of the Congruent Number Problem which have had reasonable prominence in recent years. A very nice
summary of such problems is given by Jaap Top and Noriko Yui in \cite{topyui}.

Firstly, consider a triangle with rational sides and an internal angle $\theta$, then $\cos \theta$ is rational. Suppose $\cos \theta=s/r$,
with $s,r,\in \mathbb{Z}$ and $\gcd(s,r)=1$. An integer $N$ is called $\theta$-congruent if there exists a rational-sided triangle with
one angle $\theta$ and area $N\sqrt{r^2-s^2}$.

This idea was discussed by Fujiwara \cite{fuji1}, Kan \cite{kan1} and Yoshida \cite{yosh1,yosh2}. The classical congruent number problem
corresponds to $\theta=\pi/2$.

We have, for sides $(a,b,c)$,
\begin{equation*}
\cos\theta = \frac{b^2+c^2-a^2}{2bc}=\frac{s}{r} \hspace{2cm} \sin\theta=\frac{r^2-s^2}{r}
\end{equation*}
and $bc=2rN$ so that
\begin{equation*}
a^2=b^2+c^2-4sN
\end{equation*}
and
\begin{equation}
a^2b^2=b^4-4sNb^2+4r^2N^2
\end{equation}

This quartic can be transformed into the elliptic curve
\begin{equation}
G^2=H^3+2NsH^2+N^2(s^2-r^2)H=H(H+N(r+s))(H+N(s-r))
\end{equation}
with $b=G/H$.

The most-studied problems are $\pi/3$ and $2\pi/3$ congruent number problems where $\cos\theta=\pm 1/2$,
giving respectively
\begin{equation*}
\pi/3: G^2=H^3+2NH^2-3N^2H \hspace{1cm} 2\pi/3: G^2=H^3-2NH^2-3N^2H
\end{equation*}

Both curves have discriminant
\begin{equation*}
\Delta=2304N^6
\end{equation*}
and each has $3$ points of order $2$ at $H=0,N,-3N$ and $H=0,-N,3N$ respectively.

The rank distributions for both curves, for $N \in [1,999]$ are given in
\begin{table}[H]\label{tabtheta}
\begin{center}
\caption{Rank summary for $\theta$-congruent curves}
\begin{tabular}{lcccrr}
$\,$&$\,$&$\,$&$\,$ &$\,$&$\,$\\
$\theta$&Rank$=0$&Rank$=1$&Rank$\ge 2$&Ave. ht&Max. ht.\\
$\,$&$\,$&$\,$&$\,$ &$\,$&$\,$ \\
$\pi/3$&458&481&60&9.2&103.0\\
$2\pi/3$&456&480&63&11.3&96.9
\end{tabular}
\end{center}
\end{table}

We can apply a similar approach to these equations as Elkies did to the congruent number equations. Define $H=-Nu$ and
$G=N\sqrt{-N}v$ which gives
\begin{equation*}
V^2=U^3-2U^2-3U \hspace{2cm} V^2=U^3+2U^2-3U
\end{equation*}
and we can apply the same arithmetic ideas in the quadratic field $\mathbb{Q}(\sqrt{-N})$.

With this idea, I have solutions for all $N \in [1,999]$ for both problems.

\vspace{0.5cm}

The second variant is the following:
an integer $N$ is called {\it t-congruent} we we can find $a,b,c \in \mathbb{Q}$ such that
\begin{equation}
a^2=b^2+c^2-2bc\frac{1-t^2}{t^2+1} \hspace{1cm} 2N=bc\frac{2t}{1+t^2}
\end{equation}

I first encountered this idea in a different form in a June 2001 posting to the NMBRTHRY listserver by Jim Buddenhagen.
His mailing included the statement

{\it I now know the arcane fact
that probably the smallest integer sided triangle with square integer area
and with the tangent of half of one of its angles equal to 1/23 has sides:

$a=14254641987126588457485660564361$,

$b=26615740418486740330088123851895$,

$c=40581522360016211774464382403072$,

and area = $6846355192444296386702668725552^2$}.

If we define $\alpha=\tan^{-1} \,t$, then this defines a triangle with rational sides $a,b,c$, area $N$, and one angle $2\alpha$.
Thus, the standard congruent number problem corresponds to $t=1$. In the following, we fix $N=1$ and try to find which rational
$t$ give a triangle, which we can then scale up to having integer sides, rational area, and angle $2\alpha$.

With $N=1$, we have $c=(t^2+1)/(b\,t)$, so if we substitute into the first equation and clear denominators, we must have
\begin{equation*}
t^2b^4+2t(t^2-1)b^2+(t^2+1)^2=\Box
\end{equation*}
which gives the quartic
\begin{equation}
y^2=x^4+2t(t^2-1)x^2+t^2(t^2+1)^2
\end{equation}

Mordell's method gives the elliptic curve with $t=m/n, \gcd(m,n)=1$
\begin{equation}\label{ecbudd}
G^2=H(H+m\,n^3)(H-m^3\,n)
\end{equation}
with
\begin{equation}
b=\frac{G}{mnH}
\end{equation}

The elliptic curve \eqref{ecbudd} has discriminant
\begin{equation*}
\Delta=16m^{10}n^{10}(m^2+n^2)^2
\end{equation*}
so is non-singular if $t$ is non-zero. There are $3$ clear torsion points of order $2$ at $H=0$, $H=m^3n$ and $H=-mn^3$, and
these seem to be the only ones, so we assume the torsion subgroup is isomorphic to $\mathbb{Z}/2\mathbb{Z} \oplus \mathbb{Z}/2\mathbb{Z}$.
None of the torsion points lead to a non-trivial solution.

Since there are essentially two independent parameters $m$ and $n$, we cannot consider the same sort of ranges as in other examples. There are just
over $6000$ possible pairs with $1 \le m,n \le 99$. We selected $t=m/1, m=1,\ldots,99$, $t=1/n, n=2,\ldots,99$ and $t=m/n, 2 \le m,n \le 19,
\gcd(m,n)=1$ for investigation.

\begin{table}\label{budduns}
\begin{center}
\caption{Unsolved values and heights for \eqref{ecbudd}}
\begin{tabular}{llrr}
$\,$&$\,$&$\,$&$\,$\\
m&n&Ht&2-is. Ht\\
$\,$&$\,$&$\,$&$\,$\\
61&1&470&235\\
67&1&360&180\\
74&1&435&218\\
85&1&184&92\\
88&1&107&213\\
92&1&102&204\\
94&1&441&881
\end{tabular}
\end{center}
\end{table}

\subsection{Base/Altitude}
Consider the triangle

\begin{center}
\begin{picture}(10,6)
\put(1,1){\line(1,0){12}}
\put(13,1){\line(-2,1){8}}
\put(1,1){\line(1,1){4}}
\put(5,1){\line(0,1){4}}
\put(5.25,1){\line(0,1){0.25}}
\put(5,1.25){\line(1,0){0.25}}
\put(6.5,0.25){b}
\put(3,1.25){x}
\put(8,1.25){b-x}
\put(4.5,3){a}
\put(2.5,3){y}
\put(9.5,3){z}
\end{picture}
\end{center}

We wish to find integer sides $(b,y,z)$ such that the altitude $a$ satisfies $b/a=N$ for $N$ a strictly positive integer.
Since the ratio $b/a$ is preserved by scaling, we can consider rational sides and scale.

Then, we must have
\begin{equation*}
a^2+x^2=y^2 \hspace{2cm} a^2+(b-x)^2=z^2=a^2+(Na-x)^2
\end{equation*}
and we can assume that $a=k^2-1$ and $x=2k$ with $k$ rational. Thus,
\begin{equation*}
z^2=(N^2+1)k^4-4Nk^3+2(1-N^2)k^2+4Nk+N^2+1
\end{equation*}

This quartic has $z^2=4$ when $k=1$ so we define $j=1/(k-1)$ giving
\begin{equation*}
s^2=4j^4+8(1-N)j^3+4(N^2-3N+2)j^2+4(N^2-N+1)j+N^2+1
\end{equation*}
and, defining $s=Y/2$, $j=X/2$, we have
\begin{equation*}
Y^2=X^4+4(1-N)X^3+4(N^2-3N+2)X^2+8(N^2-N+1)X+4N^2+4
\end{equation*}

This quartic can be transformed to the elliptic curve
\begin{equation}\label{ecba}
g^2=h^3+(N^2+2)h^2+h
\end{equation}
with the reverse transformation
\begin{equation}
k=\frac{g+(N+1)h+1}{g+(N-1)h-1}
\end{equation}

The curve \eqref{ecba} has the point $(0,0)$ of order $2$. There will be $3$ points of order $2$ if $(N^2+2)^2-4=N^2(N^2+4)$ is an
integer square, which clearly cannot happen when $N \in \mathbb{Z}$.

Points of order $4$ would require $h^2=1$, and we find $(-1, \pm N^2)$ are of order $4$. The fact that the points of order $4$
occur at $h=-1$ means that we cannot have points of order $8$.

Points of order $3$ would happen at a point of inflexion where
\begin{equation*}
3h^4+4(N^2+2)h^3+6h^2-1=0
\end{equation*}
and the only possible integer roots are $\pm 1$, neither of which give the quartic a zero value.

Thus, the torsion subgroup is isomorphic to $\mathbb{Z} /4\mathbb{Z}$ and all the finite torsion points lead to $a=0$.
So a possible solution requires a curve of rank at least $1$.

The computations find the first rank $1$ curve when $N=5$. The point $(-25,35)$ lies on $y^2=x^3+27x^2+x$ and gives
$k=19/11$. We thus have $a=240/121$, $b=1200/121$, $y=482/121$ and $z=818/121$, giving the basic integer triangle
with sides $(241,409,600)$. We can check that this triangle does have $b=5a$ by calculating the area using Heron's
formula. We find $\Delta=36000$, implying $a=2\Delta/b=72000/600=120$.

The discriminant of the elliptic curve is $16N^2(N^2+4)$, so the curves are non-singular for $N>0$. The BSD calculations
suggest that for $N \in [1,999]$, there are $391$ rank zero curves, $498$ rank one curves and $110$ curves with higher rank.
The average height of the rank-one curves is $103.5$ and the maximum height is $1466$. These values are somewhat misleading,
since the largest height of an unsolved value is $630$. We have been able to use both the 2-isogenous curve
\begin{equation*}
y^2=x^3-2(N^2+2)x^2+N^2(N^2+4)x
\end{equation*}
and a 4-isogenous curve
\begin{equation*}
g^2=h^3+2(4-N^2)h^2+(N^2+4)^2h
\end{equation*}
which also has its own 2-isogeny,  to find solutions with large height.

With all these curves, I have reduced the number of values of $N$ with no known solution to the $7$
given in Table $4.3$.

\begin{table}\label{bauns}
\begin{center}
\caption{Unsolved values and heights for $b/a=N$}
\begin{tabular}{lrrrr}
$\,$&$\,$&$\,$&$\,$ &$\,$\\
N&Ht&2-is. Ht&4-is. Ht&2-is of 4-is Ht\\
$\,$&$\,$&$\,$&$\,$ &$\,$\\
599&630&315&158&315\\
683&471&236&118&236\\
907&368&184&92&184\\
919&556&278&139&278\\
921&569&284&142&284\\
947&383&182&96&192\\
956&433&216&108&216
\end{tabular}
\end{center}
\end{table}
As can be seen, the heights of the original elliptic curves seem to be higher than for the isogenous curves.

It is clearly an alternative problem to look for integer triangles where $a/b=N$. Rather than go
through the whole process of deriving an elliptic curve again, it is easier to realise that this
is just $b/a=1/N$. Replacing $N$ by $1/N$ in \eqref{ecba} and simplifying gives
the elliptic curve
\begin{equation}\label{ecab}
g^2=h^3+(2N^2+1)h^2+N^4h
\end{equation}
with $a=k^2-1$, $b=a/N$, $x=2k$, $a^2+x^2=y^2$, and $a^2+(b-x)^2=z^2$ as before. The triangle $(b,y,z)$
can be scaled to integer sides and
\begin{equation*}
k=\frac{g+(N+1)h+N^3}{g+(1-N)h-N^3}
\end{equation*}

These elliptic curves have discriminant
\begin{equation*}
\Delta=16N^8(4N^2+1)
\end{equation*}
and torsion points $(0,0)$ of order $2$ and $(-N^2, \pm N^2)$ of order $4$. The BSD calculations predict
$454$ rank zero curves, $468$ rank one and $77$ with rank at least two. The average height of the rank
one curves is $117$ with the maximum height $3578$.

Although these heights are similar to the first problem, currently there are $78$ unsolved values. These values
seem to have the majorities of the heights of the original curves smaller than those of the 2-isogeny and 4-isogeny curves.
The unsolved values of $N$ together with the estimated heights for the original curves are given in the following Table $4.4$. The asterisk
indicates those values of $N$ where the height on the 4-isogeny curve is a quarter of the original curve.

\begin{table}\label{abuns}
\begin{center}
\caption{Unsolved values and heights for $a/b=N$}
\begin{tabular}{llllllrrrrr}
$\,$&$\,$&$\,$&$\,$&$\,$&$\,$&$\,$&$\,$&$\,$&$\,$&$\,$\\
N&Ht&$\,$&N&Ht&$\,$&N&Ht&$\,$&N&Ht\\
$\,$&$\,$&$\,$&$\,$&$\,$&$\,$&$\,$&$\,$&$\,$&$\,$&$\,$\\
233&279&$\,$&317&463&$\,$&337&129&$\,$&389&496*\\
401&914*&$\,$&439&569*&$\,$&471&187&$\,$&481&216\\
499&830*&$\,$&502&160&$\,$&514&544*&$\,$&521&825*\\
537&103&$\,$&556&115&$\,$&562&469&$\,$&565&171\\
583&151&$\,$&586&576*&$\,$&596&97&$\,$&599&899*\\
601&723*&$\,$&604&225&$\,$&607&1060*&$\,$&622&150\\
628&253&$\,$&633&177&$\,$&639&129&$\,$&647&399\\
653&256&$\,$&674&593*&$\,$&677&463&$\,$&683&862\\
691&667*&$\,$&692&106&$\,$&721&129&$\,$&723&284\\
739&554*&$\,$&749&179&$\,$&757&611*&$\,$&772&180\\
773&1140&$\,$&778&108&$\,$&786&233&$\,$&791&101\\
796&218&$\,$&803&171&$\,$&815&108&$\,$&817&110\\
827&237&$\,$&829&900*&$\,$&838&452&$\,$&839&419*\\
842&745&$\,$&845&110&$\,$&871&126&$\,$&873&113\\
878&340&$\,$&883&167&$\,$&887&1430&$\,$&889&267\\
913&119&$\,$&917&146&$\,$&919&3580*&$\,$&921&155\\
922&962*&$\,$&926&372*&$\,$&931&102&$\,$&933&180\\
934&394*&$\,$&943&120&$\,$&947&1380&$\,$&953&478*\\
961&1090*&$\,$&963&162&$\,$&978&142&$\,$&982&1400\\
983&563&$\,$&989&267&$\,$&$\,$&$\,$&$\,$&$\,$&$\,$\\
\end{tabular}
\end{center}
\end{table}

\subsection{Leech's Problem}
Consider the following problem.

{\it Find two integer right-angled triangles on the same base
whose heights are in the ratio N:1 for N an integer.}

\begin{center}
\begin{picture}(10,7)
\put(1,1){\line(1,0){12}}
\put(1,1){\line(6,1){12}}
\put(1,1){\line(2,1){12}}
\put(13,1){\line(0,1){6}}
\put(12.75,1){\line(0,1){0.25}}
\put(12.75,1.25){\line(1,0){0.25}}
\put(7,0.5){b}
\put(13.25,2){a}
\put(7,2.25){c}
\put(6,4.25){d}
\put(13.75,4){Na}
\end{picture}
\end{center}

I first read of this in a paper by Chris Smyth \cite{smy}, where he attributes the problem to the late John Leech.

Thus,
\begin{equation*}
b^2+a^2=c^2 \hspace{3cm} b^2+N^2a^2=d^2
\end{equation*}

So, $b=\alpha (p^2-q^2)$ and $a=\alpha 2pq$, leading to
\begin{equation}
p^4+(4N^2-2)p^2q^2+q^4=\Box
\end{equation}

This quartic can be transformed to the equivalent elliptic curve
\begin{equation}\label{ecleech}
g^2=h^3+(N^2+1)h^2+N^2h=h(h+1)(h+N^2)
\end{equation}
with
\begin{equation}
\frac{p}{q}=\frac{g}{h+N^2}
\end{equation}

The elliptic curve has points of order $2$ at $(0,0)$, $(-1,0)$ and $(-N^2,0)$ which give $p/q$ either $0$ or undefined. There
are also points of order $4$ at $(N, \pm N(N+1))$ and $(-N,\pm N(N-1))$ all of which give $|p/q|=1$, so none of these torsion points lead
to a non-trivial solution.

There are only two possible structures for the torsion subgroup - isomorphic to $\mathbb{Z}/2\mathbb{Z} \oplus \mathbb{Z}/4\mathbb{Z}$
or to $\mathbb{Z}/2\mathbb{Z} \oplus \mathbb{Z}/8\mathbb{Z}$. For the latter $(N,\pm N(N+1))$ must be double a point of order $8$.
Thus $N$ must be an integer square, say $M^2$.

If $(r,s)$ has order $8$, then
\begin{equation*}
\frac{r^4-2M^4r^2+M^8}{4r(r+1)(r+M^4)}=M^2
\end{equation*}
giving
\begin{equation*}
(r^2-2M(M^2+M+1)r+M^4)(r^2+2M(M^2-M+1)+M^4)=0
\end{equation*}
and we have integer roots only if $M^2+1=\Box$, which does not happen.

Thus the torsion subgroup is isomorphic to $\mathbb{Z}/2\mathbb{Z} \oplus \mathbb{Z}/4\mathbb{Z}$, and none of the
finite torsion points give a non-trivial answer. So we need the rank to be at least one.

As an example, for $N=28$, we have the elliptic curve $g^2=h^3+785h^2+784h$, which has a point of infinite order at
$(112/9,9856/27)$. This gives $p/q=11/14$, giving $b=455, a=528$. The hypoteneuse of the smaller triangle
is $697$ and $14791$ for the larger.

The curves \eqref{ecleech} have discriminant
\begin{equation*}
\Delta=16N^4(N^2-1)^2
\end{equation*}
and so are singular if $N=1$.

Searching for a parametric solution, we find that $N=4k^2+3k$ and $h=k$ give $g=\pm k(4k^2+5k+1)$. These lead to $b=4k(2k+1)$ and
$a=4k+1$.

The BSD calculations for $N \in [2,999]$ give $398$ curves of rank $0$, $514$ curves of rank $1$
and $86$ curves with rank at least $2$. I have solutions for all the positive rank curves.

\subsection{$R/r=N$}
Given a triangle with integer sides $(a,b,c)$, the circumcircle is the circle through the $3$ vertices of the triangle, with radius denoted $R$.
The incircle is the circle contained inside the triangle with each side a tangent to the circle, with radius denoted $r$.
Trigonometry gives several formulae for these radii, amongst which are
\begin{equation*}
R=\frac{a\,b\,c}{4\Delta} \hspace{2cm} r=\frac{\Delta}{s}
\end{equation*}
where $\Delta$ is the area of the triangle and $s=(a+b+c)/2$ is the semi-perimeter.

Thus, the dimension-less quantity $R/r$ satisfies
\begin{equation*}
\frac{R}{r}=\frac{a\,b\,c\,s}{4\Delta^2}=\frac{a\,b\,c}{4(s-a)(s-b)(s-c)}=\frac{2\,a\,b\,c}{(a+b-c)(b+c-a)(c+a-b)}
\end{equation*}

An equilateral triangle has $R/r=2$, and it can be shown that this is the smallest value for the ratio. The question we ask is, what other
integer ratios can $R/r$ take with integer triangles? Since the ratio is independent of any scaling of the sides, we can find rational
sides and scale to integers. Thus we wish to solve
\begin{equation}
N=\frac{2\,a\,b\,c}{(a+b-c)(b+c-a)(c+a-b)}
\end{equation}
if possible.

We have,
\begin{equation*}
Na^3-Na^2(b+c)-a(Nb^2-2(N+1)bc+Nc^2)+N(b+c)(b-c)^2=0
\end{equation*}

Cubics are not very amenable objects, but we can reduce the problem to a quadratic by setting $c=P-a-b$, where $P$ will be
the perimeter of the triangle. The quadratic in $a$ is
\begin{equation*}
2(2NP-(4N+1)b)a^2-2(b^2(4N+1)-bP(6N+1)+2NP^2)a+
\end{equation*}
\begin{equation*}
NP(4b^2-4bP+P^2)=0
\end{equation*}

For this to have a rational solution $a$, the discriminant must be a rational square, so there exists $d \in \mathbb{Q}$ with
\begin{equation*}
d^2=-2NbP^3+(4N^2+8N+1)b^2P^2-2(2N+1)(4N+1)b^3P+(4N+1)^2b^4
\end{equation*}

Setting $d=b^2v$ and $P=bu$ and assuming $b \ne 0$, we have
\begin{equation}
v^2=-2Nu^3+(4N^2+8N+1)u^2-2(8N^2+6N+1)u+16N^2+8N+1
\end{equation}

Define $v=y/(2N)$ and $u=-z/(2N)$ to give
\begin{equation*}
y^2=z^3+(4N^2+8N+1)z^2+4N(2N+1)(4N+1)z+4N^2(4N+1)^2
\end{equation*}

The right-hand side factors and defining $z=x-4N-1$ we have
\begin{equation}\label{ecbw}
y^2=x^3+2(2N^2-2N-1)x^2+(4N+1)x
\end{equation}
with the reverse transformation
\begin{equation}
\frac{P}{b}=\frac{4N+1-x}{2N}
\end{equation}

The elliptic curve has discriminant
\begin{equation*}
\Delta = 256N^3(N-2)(4N+1)^2
\end{equation*}
so is singular when $N=2$. For $N >2$, we have $\Delta > 0$, so the curves \eqref{ecbw} all have two components, with the "egg" lying
in the negative $x$ half-plane.

There is an obvious point of order $2$ at $(0,0)$. There are other points of order $2$ when $N(N-2)=\Box$, which does not happen when $N$ is an integer.
There are points of inflexion and hence points of order $3$ at $(1, \pm 2N)$. There are also points of order $6$ at $(4N+1, \pm 2N(4N+1))$,
so the torsion subgroup is either $\mathbb{Z} / 6 \mathbb{Z}$ or $\mathbb{Z} / 12 \mathbb{Z}$. The latter would need $4N+1=(2M+1)^2=J^2$
with $J,M \in \mathbb{Z}$,
and a rational point of order $12$ would satisfy
\begin{equation*}
(x^2+2J(2M^2-1)x+J^2)(x^2-2J(2M^2+4M+1)x+J^2)=0
\end{equation*}
The first quadratic has rational solutions if $M^2-1=\Box$ whilst the second requires $M(M+2)=\Box$, neither of which occurs if $N > 2$.
Thus the torsion subgroup is isomorphic to $\mathbb{Z} / 6 \mathbb{Z}$.

None of the torsion points lead to non-trivial solutions, so we need curves of rank greater than zero. using the BSD conjecture to estimate the rank, we find
$415$ curves of rank zero, $502$ curves of rank one and $80$ curves with rank greater than one for $1 \le N \le 999$. The average height for the rank $1$ curves
is $43.8$ and the maximum estimated height is $447.2$.

There are several isogenous curves we can use to find generators. The curve \eqref{ecbw} has a 2-isogeny with
\begin{equation}\label{ecbw2}
v^2=u^3-4(2N^2-2N-1)u^2+16N^3(N-2)u
\end{equation}

The fact that there is a point of inflexion at $x=1$ allows us to set $x=z+1$, giving the equivalent curve
\begin{equation*}
y^2=z^3+(\,(2N-1)z+2N \,)^2
\end{equation*}
and the formulae for a 3-isogeny given in the appendices eventually lead (after some transformations) to
\begin{equation}\label{ecbw3}
g^2=f^3+18(2N^2+10N-1)f^2+81(4N+1)^3f
\end{equation}
which has its own 2-isogenous curve,
\begin{equation}\label{ecbw6}
i^2=j^3-36(2N^2+10N-1)j^2+1296N(N-2)^3j
\end{equation}

Using all of these curves, I have reduced the number of unsolved values of $N$ to four. The values of $N$ and the particular version of elliptic curve
giving the smallest height (for that $N$) are given in Table \ref{bwuns}.

\begin{table}[H]\label{bwuns}
\begin{center}
\caption{Unsolved values of $N$}
\begin{tabular}{lrr}
$\;$&$\;$&$\;$\\
$N$&Height&Curve\\
683&150.9&\eqref{ecbw2}\\
619&168.9&\eqref{ecbw2}\\
883&197.4&\eqref{ecbw2}\\
773&223.6&\eqref{ecbw2}\\
\end{tabular}
\end{center}
\end{table}

There is a problem in just finding points on the elliptic curve - most points give one of $(a,b,c)$ negative! Thus we do not get a
triangle. We have to analyze deeper to find out which points (if any) give a real-life triangle.

We have
\begin{equation*}
u=\frac{P}{b}=\frac{4N+1-x}{2N}=\frac{a+b+c}{b}=1+\frac{a+c}{b} > 2
\end{equation*}
for a triangle. Thus $x<1$ is a requirement.

The quadratic in $a$ can be written
\begin{equation*}
a^2+b(1-u)a+\frac{b^2Nu(u-2)^2}{2(2N(u-2)-1)}
\end{equation*}
and the two roots will actually give the values for $a$ and $c$.

Thus we need
\begin{equation*}
1-u<0 \hspace{2cm} \frac{u}{2N(u-2)-1}>0
\end{equation*}

The fact that we need $u>2$ from the triangle inequality gives the first and the second implies that $x<0$. This means that the point must lie on the egg component.
Since all the torsion points lie on the infinite component, adding a point on the infinite component to a torsion point gives another point on the infinite component.
Also, doubling a point always gives a point on the infinite component.

Thus, if all the generators are on the infinite component, we {\bf never} have a point on the egg. So, we have a real-life triangle only if there is a generator point
which lies on the egg component. These are thin on the ground. So far, the only solutions found are

\begin{table}[H]
\begin{center}
\caption{Integer sided triangles with $R/r=N$}
\begin{tabular}{lrrr}
$\;$&$\;$&$\;$&$\;$ \\
N&$f$&$g$&$h$\\
2 & 1 & 1& 1 \\
 26 &     11 &     39  &    49\\
  74  &    259&     475 &    729\\
 218  &   115&     5239 &   5341\\
 250  &   97 &     10051&   10125\\
 314  &  177487799 &  55017780825 &    55036428301\\
 386  &   1449346321141 &  2477091825117 &  3921344505997\\
 394  &   12017&   2356695 &    2365193\\
 458  &   395  &   100989  &    101251\\
 586  &   3809 &   18411 &  22201\\
 602  &   833  &   14703 &  15523\\
 634  &   10553413  &  1234267713   &   1243789375\\
 674  &   535  &   170471 &     170859\\
 746  &   47867463&    6738962807 &     6782043733\\
 778  &   1224233861981  & 91266858701995   &   92430153628659\\
 866  &   3025  &  5629  &  8649
\end{tabular}
\end{center}
\end{table}

\subsection{Relations between the altitude, angle-bisector and median}
Consider the triangle $ABC$

\begin{center}
\begin{picture}(10,6)
\put(1,1){\line(1,0){12}}
\put(13,1){\line(-5,2){10}}
\put(1,1){\line(1,2){2}}
\put(3,1){\line(0,1){4}}
\put(3,5){\line(1,-2){2}}
\put(3,5){\line(1,-1){4}}
\put(0.5,1){C}
\put(3,5.2){A}
\put(13.2,0.9){B}
\put(6,0.5){a}
\put(1.6,3){b}
\put(8,3.1){c}
\put(2.5,2.5){h}
\put(3.8,2.6){t}
\put(5.3,2.8){m}
\end{picture}
\end{center}
where we assume $a,b,c \in \mathbb{Q}$. $h$ is the altitude, $t$ is the angle-bisector and $m$ is the median.

If $b=c$, then $h=t=m$, and it is a natural question to ask whether the ratio of two of these lengths can be an integer
greater than $1$.

Standard trigonometry gives
\begin{equation}
m^2=\frac{2b^2+2c^2-a^2}{4} \hspace{2cm} t^2=b\,c\,\left( 1-\frac{a^2}{(b+c)^2} \right)
\end{equation}
\begin{equation}
 h^2=\frac{(a+b+c)(a+b-c)(b+c-a)(c+a-b)}{4a^2}
\end{equation}
and it is straightforward, using these formulae, to show
\begin{equation*}
h^2 \le t^2 \le m^2
\end{equation*}

Consider, first, the relation $m/t=N$ where $N$ is an integer greater than $1$. Then we have the equation
\begin{equation}\label{mteq}
\frac{(a^2-2(b^2+c^2))(b+c)^2}{4bc(a^2-(b+c)^2)}=N^2
\end{equation}
which gives the quadratic equation for $a$
\begin{equation}
(\, b^2+2bc(1-2N^2)+c^2 \,)\, a^2=2(b+c)^2(b^2-2bcN^2+c^2)
\end{equation}

For $a$ to be rational, we must have
\begin{equation}
d^2=2(b^2+2bc(1-2N^2)+c^2)(b^2-2bcN^2+c^2)
\end{equation}
with $d \in \mathbb{Q}$. Defining $d=c^2j$ and $b=cz$ gives
\begin{equation}
j^2=2z^4+4(1-3N^2)z^3+4(4N^4-2N^2+1)z^2+4(1-3N^2)z+2
\end{equation}

For $z=1$ we have $j=\pm 4(N^2-1)$, so the quartic is birationally equivalent to an elliptic curve. Using
standard transformations described in Appendix B, we can find this elliptic curve to be
\begin{equation}\label{eccmt}
v^2=u^3+2(2N^4-N^2-1)u^2+(N^2-1)^2u
\end{equation}
with the reverse transformation given by
\begin{equation}
z=1+\frac{1}{w} \hspace{2cm} w=\frac{v}{4(N^2-1)u}-\frac{1}{2}.
\end{equation}

The elliptic curve \eqref{eccmt} has a point of order $2$ at $(0,0)$ and numerical tests suggest
this is the only torsion point, but I have not tried to prove this. This point
gives an undefined value for $z$, suggesting we need to look for points of infinite order.

Numerical tests suggested that, for $N$ a positive integer, the rank is always strictly positive,
and we quickly found the points $(\, (N-1)^2, \pm 2N^2(N-1)^2 \,)$ and $(\, (N+1)^2, \pm 2N^2(N+1)^2 \,)$.
These give $z=b/c=2N^2-1$, and using $b=2N^2-1, \, c=1$ gives $a=2N$ which is a numerical solution
of \eqref{mteq} but does not give an acceptable triangle.

The tangent to the elliptic curve at $(\, (N-1)^2, 2N^2(N-1)^2 \,)$ has equation
\begin{equation*}
v=2N^2(N-1)^2+\frac{2N^3-1}{N}(u-(N-1)^2)
\end{equation*}
and meets the curve again at $u=1/N^2$ with $v = \pm (N^4+N^2-1)/N^3$. This gives
\begin{equation}
\frac{b}{c}=\frac{N^4+2N^3+N^2-2N-1}{N^4-2N^3+N^2+2N-1}
\end{equation}

Using the numerator for $b$ and denominator for $c$ in equation \eqref{mteq} gives
\begin{equation}
a=\frac{2(N^8-3N^4+1)}{N^4-N^2+1}
\end{equation}
and scaling by $N^4-N^2+1$ gives
\begin{equation}
a=2(N^8-3N^4+1)
\end{equation}
\begin{equation}
b=(N^4+2N^3+N^2-2N-1)(N^4-N^2+1)
\end{equation}
\begin{equation}
c=(N^4-2N^3+N^2+2N-1)(N^4-N^2+1)
\end{equation}

It is easy to see that, if $N \ge 2$, then all three expressions are positive and that they satisfy
the triangle inequalities. These, therefore, give a parametric solution to the problem. We have
\begin{equation*}
m=2N^2\sqrt{(N^2-1)(N^4+2N^3+N^2-2N-1)(N^4-2N^3+N^2+2N-1)}
\end{equation*}
and numerical investigations have been unable to find a value of $N$ giving $m \in \mathbb{Z}$.

\vspace{1cm}

For $m/h=N$, we have
\begin{equation*}
(N^2-1)a^4+2(1-N^2)(b^2+c^2)a^2+N^2(b^2-c^2)^2=0
\end{equation*}
which initially might seem intractable as each of $a,b,c$ are fourth powers. But, if we set
$c=p-a-b$, we get a quadratic in $b$, though still a quartic in $a$! The quadratic is
\begin{equation*}
4(a^2-2aN^2p+N^2p^2)b^2+4(a^3-a^2p(2N^2+1)+3aN^2p^2-N^2p^3)b+
\end{equation*}
\begin{equation*}
a^4-4a^3p+2a^2p^2(2N^2+1)-4aN^2p^3+N^2p^4 = 0
\end{equation*}

To have $b \in \mathbb{Q}$, the discriminant must be a rational square. Clearing out any square terms,
we must then have $d \in \mathbb{Q}$ such that
\begin{equation*}
d^2=2a^3p(1-N^2)+a^2p^2(N^2-1)(4N^2+1)+4N^2ap^3(1-N^2)+N^2p^4(N^2-1)
\end{equation*}

Define $z=a/p$, $k=d/p^2$, $k=y/(2-2N^2)$, $z=w/(2-2N^2)$, and $w=x-(N^2-1)$ leading to
\begin{equation}\label{ecmh}
y^2=x^3+(4N^4-6N^2+2)x^2+(N^2-1)^2x
\end{equation}
with
\begin{equation}
\frac{a}{p}=z=\frac{x}{2(1-N^2)}+\frac{1}{2}
\end{equation}

The equation $y^2=0$ has roots at $x=0$ and
\begin{equation*}
-(2N^4-3N^2+1) \pm 2N(N^2-1)\sqrt{N^2-1}
\end{equation*}
and so the elliptic curve has $2$ components, the infinite one for $x \ge 0$ and the "egg"
which is a closed convex curve in $x<0$.

The elliptic curve has the obvious point at $(0,0)$ giving $a=1,p=2$ which give $b=c=1/2$,
again an arithmetic solution but not a triangular one. There are also $2$ points of order $4$
when $x=1-N^2$ with $y=\pm 2(N^2-1)^2$, giving $a/p=1$ which certainly does not give a real-life
triangle. Numerical tests suggest these are the only torsion points, implying that we would need
points of infinite order.

The numerical experiments also suggested that the ranks of \eqref{ecmh} can often be zero.
The BSD computations suggested that the heights of rank-1 generators can be very large even for small $N$.
For example, $N=58$ has an estimated height of $280.3/560.6$. For $N>99$, the heights can be enormous, so
we restricted the range to $N \in [2,99]$. Most of the large heights come when $N$ is even with at least one
of $N-1$ and $N+1$ prime.

\begin{table}[H]
\begin{center}
\caption{Small example triangles }
\begin{tabular}{lrrr}
$\,$&$\,$&$\,$&$\,$\\
N&a&b&c\\
$\,$&$\,$&$\,$&$\,$\\
4 & 238 & 241 & 31\\
5 & 50 & 59 & 11\\
8 & 23838 & 21191 & 2921\\
13 & 18067634 & 15502013 & 2621453\\
15 & 2938 & 3151 & 239\\
25 & 1750 & 1699 & 61\\
26 & 43654 & 42673 & 1273\\
29 & 1682 & 2029 & 349\\
31 & 315854 & 313007 & 5777\\
49 & 5899006 & 5814197 & 103253\\
55 & 4502900 & 4399621 & 110581
\end{tabular}
\end{center}
\end{table}

\vspace{1cm}
For completeness, we finish by discussing
\begin{equation*}
\frac{t}{h}=N
\end{equation*}
even though the analysis does not use elliptic curves.

We have
\begin{equation}
\frac{t^2}{h^2}=\frac{4a^2bc}{(a+b-c)(a+c-b)(b+c)^2}=N^2
\end{equation}
giving the equation for $a$
\begin{equation}
(b^2N^2+2bc(N^2-2)+c^2N^2)a^2=N^2(b^2-c^2)^2
\end{equation}

Thus for $a$ to be rational we require
\begin{equation}
d^2=(b^2N^2+2bc(N^2-2)+c^2N^2)
\end{equation}
to have rational solutions. Define $d=cy$ and $b=cx$, giving the quadric
\begin{equation}
y^2=N^2x^2+2(N^2-2)x+N^2
\end{equation}
which has rational solutions $x=0, y=\pm N$.

The line $y=N+kx$ through $(0,N)$ meets the curve again at
\begin{equation}
x=\frac{b}{c}=\frac{2(Nk+2-N^2)}{N^2-k^2}
\end{equation}

Substituting $b=Nk+2-N^2$ and $c=N^2-k^2$ into the equation for $a$ we have
\begin{equation}
a=\pm \frac{N(k-N-2)(k-N+2)(k^2+2kN-3N^2+4)}{k^2N+2k(2-N^2)+N^3}
\end{equation}
so we are guaranteed a positive value of $a$.

For positive $b$ we must have $k > N-2/N$ whilst positive $c$ requires $|k| < N$, so we are
restricted to $N-2/N < k < N$. Choosing $k=N-1/N$ and clearing denominators gives
\begin{equation}
a=4N^2-1 \hspace{1cm} b=2N^2(4N^2-3) \hspace{1cm} c=(2N^2-1)(4N^2-3)
\end{equation}
and it is easy to check that these values satisfy the triangle inequalities.

Other Diophantine problems are possible for these three lengths, see the recent paper of Bakker, Chahal and Top \cite{bct}.

\newpage

\section{Two Quadrics Simultaneously Square}

\subsection{Concordant Forms and Related Problems}
Let $M, N \in \mathbb{Z}$ with $M \ne N$ and $MN \ne 0$, and consider the problem of finding integers $x,y$ with $|y|>0$ such that
\begin{equation}
x^2+My^2=\Box \hspace{2.5cm} x^2+Ny^2=\Box
\end{equation}

If we can find $(x,y)$, these expressions are called concordant forms, a phrase given to them by Euler in the paper {\it De binis formulis
speciei xx+myy et xx+nyy inter se concordibus et disconcordibus}.

A very interesting modern
paper is by Ken Ono \cite{ono}. Concordant forms take up $6$ pages in Chapter XVI of Dickson \cite{dick1}.

Define $X=y/x$ and consider the equation $Y^2=M\,X^2+1$. This has solution $X=0, \, Y=1$, and the line $Y=1+kX$ meets the quadric
in one further point
\begin{equation*}
X=\frac{2k}{M-k^2}
\end{equation*}

The second identity requires $N\,X^2+1=\Box$, so the parameter $k$ must lie on the quartic
\begin{equation*}
d^2=k^4+2(2N-M)k^2+M^2
\end{equation*}
and there is an obvious rational solution at $k=0, \, d=M$, so the quartic is birationally equivalent to the elliptic curve
\begin{equation}
G^2=H^3+(M+N)H^2+M\,N\,H=H(H+M)(H+N)
\end{equation}
with
\begin{equation*}
k=\frac{G}{H+N}
\end{equation*}
and we get from this
\begin{equation}
\frac{y}{x}=\frac{2G}{MN-H^2}
\end{equation}

There are $3$ clear points of order $2$, at $(0,0)$, $(-M,0)$ and $(-N,0)$, which give $k=0$ or $k$ undefined. Some choices
of $M$ and $N$ give other torsion points and this is thoroughly discussed by Ono.

It should be noted that $M=-N$ gives the Congruent Number Problem, and $M=n^2, \, N=1$ gives Leech's problem, with both problems discussed
in the previous section.

Because of the two parameters $M,N$, there are a huge number of possible combinations that can arise. We can assume, without loss of
generality, that $|M|<|N|$ and that $\gcd(M,N)$ is squarefree. For $1 \le M , N \le 99$, there are $4517$ elliptic curves with $2000$ having rank zero,
$2208$ having rank one, and $309$ having rank greater than one.

The largest rank-1 height, according to the BSD calculations, comes from $M=74, N=97$ with a value of $13.01$. We easily find a generator and the
fact that $x=23697472157355594548677$ and $y=3456643292842216826580$ give
\begin{equation*}
x^2+74y^2=38023026153122318249173^2
\end{equation*}
\begin{equation*}
x^2+97y^2=41479673618314533708877^2
\end{equation*}

To reduce the variables to one parameter, we selected the following choices of $M,\,N$ to be considered
\begin{table}[H]
\begin{center}
\caption{Choices of $M, \, N$}
\begin{tabular}{ccccc}
$\;$&$\;$&$\;$&$\;$&$\,$\\
M&$\;$&N&$\,$&$\Delta$\\
$\;$&$\;$&$\;$&$\;$&$\,$\\
$1$&$\;$&$N$&$\,$&$16N^2(N-1)^2$\\
$N+1$&$\;$&$N$&$\,$&$16N^2(N+1)^2$\\
$N^2$&$\;$&$N$&$\,$&$16N^8(N-1)^2$\\
$1/N$&$\;$&$N$&$\,$&$16N^{10}(N^2-1)^2$
\end{tabular}
\end{center}
\end{table}

The sharp-eyed reader will have noticed that the final example does not have $M$ an integer. It is easy to see, however,
that the problem with $(M,N)=(1/N,N)$ is equivalent to $(M,N)=(N^3,N)$.

The first problem is clearly the easiest to deal with. Values of $N$ satisfying $x^2+y^2=\Box$ and $x^2+Ny^2=\Box$ with $x\,y \ne 0$
are known as Concordant Numbers. Dickson states that Brooks and Watson, in $1857$, gave the following as the
values of $N$ up to $100$ admitting a solution: 1, 7, 10, 11, 17, 20, 22, 23, 24, 27, 30, 31, 34, 41, 42, 45, 49, 50, 52, 57, 58, 59,
60, 61, 68, 71, 72, 74, 76, 77, 79, 82, 85, 86, 90, 92, 93, 94, 97, 99, 100.

As pointed out in the Wolfram Mathworld entry on Concordant Forms, this list misses out $N=47, 53, 83$. This should, in no way,
take anything away from Brooks and Watson who, presumably, did all their calculations by hand. It is, in many respects, amazing
that they only missed out three values. The Mathworld entry also contains a list of $16$ primes in $[1,999]$ which are not
known to be concordant or discordant. All of these can be rejected on the basis of Birch and Swinnerton-Dyer conjecture computations.
Presumably, Mathworld refers to the lack of a complete proof of their discordance.

Based on the BSD conjecture, I have computed solutions for all the values in $[-999,9999]$ which are predicted to have a non-trivial solution.
For example, the first example not found by Brooks and Watson is for $N=47$. This corresponds to the curve $G^2=H(H+1)(H+47)$, and it is
reasonably easy to find $H=-1296/169, G=\pm 98532/2197$, giving $x=14663, y=111384$ and
\begin{equation*}
x^2+y^2=112345^2 \hspace{2cm} x^2+47\, y^2= 763751^2
\end{equation*}

Similarly, I have solutions for all $N$ in $[1,999]$ for the $(N,N+1)$ problem. Here, the curve is $G^2=H^3+(2N+1)H^2+N(N+1)H$. As an example,
the BSD computations predict rank $1$ for $N=52$ with height of generator $5.08$. We find $H=13*77^2/18^2$ and $G= \pm 13*77*26095/18^3$ which leads
to $x=434734621$ and $y=72335340$, with
\begin{equation*}
x^2+52y^2=679028029^2 \hspace{2cm} x^2+53y^2=682870021^2
\end{equation*}

The other two problems lead to much higher heights, as can be seen in the summary Table.

\begin{table}[H]
\begin{center}
\caption{Rank summary for Concordant and related curves}
\begin{tabular}{lcccrr}
$\,$&$\,$&$\,$&$\,$ &$\,$&$\,$\\
(M,N)&Rank$=0$&Rank$=1$&Rank$\ge 2$&Ave. ht&Max. ht.\\
$\,$&$\,$&$\,$&$\,$ &$\,$&$\,$ \\
$(1,N)$&413 &500 &85&3.2&24.3\\
$(N,N+1)$&472 &487 &40&8.9&85.9\\
$(N,N^2)$&399 &504 &95&28.6&458.6\\
$(N,1/N)$&414&510&74&182.9&5101.2
\end{tabular}
\end{center}
\end{table}

The $(N,N^2)$ problem has curve $G^2=H^3+(N+N^2)H^2+N^3H$. As an example, $N=74$ has rank $1$, with $H=12823561/11664$
and $G=\pm 116021624725/1259712$. This gives $x=1497444368329$ and $34329686220$, with
\begin{equation*}
x^2+74y^2=3311108888329^2 \hspace{2cm} x^2+74^2y^2=25448063182921^2
\end{equation*}

The much larger heights for this problem mean that there are still unsolved values of $N$. Currently, there are $16$ values of $N$, which are
predicted to give a solution, but no point has been found. These values are given in Table $5.3$. The 2-isogenous curve is
$Z^2=W^3-2(N+N^2)W^2+(N^2-N)^2W$. In all cases, the 2-isogenous curve estimate is half that of the original curve.

\begin{table}[H]\label{nnsq}
\begin{center}
\caption{Unsolved values for $(N,N^2)$ problem}
\begin{tabular}{lrr}
$\,$&$\,$&$\,$\\
$N$&Isog. Ht.&Conductor\\
$\,$&$\,$&$\,$\\
727	&93.94&558126624\\
479	&96.80&1754764768\\
998	&98.38&15888255808\\
503	&101.7&2032168288\\
829	&103.1&758714064\\
982	&104.6&5045343168\\
863	&107.7&10271854048\\
843	&107.9&9573863328\\
797	&123.3&2022505456\\
839	&125.4&9438172768\\
557	&132.3&689993776\\
653	&140.2&1112074672\\
857	&150.7&1257376688\\
997	&198.2&3960131856\\
983	&221.4&15182332768\\
563&229.3&2850185248
\end{tabular}
\end{center}
\end{table}

\vspace{1cm}

The $(N,1/N)$ problem has elliptic curve
\begin{equation*}
V^2=U(U+N)(U+N^3)
\end{equation*}
with
\begin{equation}
\frac{y}{x}=\frac{2NV}{N^4-U^2}
\end{equation}

For example, $N=23$ has a point $(532900/169,860597730/2197)$, which eventually lead to the solutions
\begin{equation*}
x=919\,\,9662\,\,3733 \hspace{2cm} y=1715\,\,4581\,\,4180
\end{equation*}

The elliptic curve has $3$ finite torsion points of order $2$, at $(0,0)$, $(-N,0)$ and $(-N^3,0)$. The discriminant
is
\begin{equation*}
\Delta=16N^{10}(N^2-1)^2
\end{equation*}
and so the curves are singular for $N=\pm 1$. Apart from these values and $N=0$, $\Delta>0$ and so the curves
have two components. The distribution of rank estimates is given in the last row of Table 5.2.

As can be seen, the average heights are much larger than the other curves in this section. There are, thus, far more
unsolved values of $N$ than these other curves. This set of curves also posed more problems in computation. The $64$
values of $N$, which are predicted to give a solution, but for which no solution has been found are given in Table $5.4$.
In all cases, the given height is that of the 2-isogenous curve which is smaller than that of the original curve.
\begin{equation*}
g^2=h^3-2(N+N^3)h^2+N^2(1-N^2)^2h
\end{equation*}

\begin{table}\label{dd8uns}
\begin{center}
\caption{Unsolved values and heights for the $(N,1/N)$ problem}
\begin{tabular}{llllllrrrrr}
$\,$&$\,$&$\,$&$\,$&$\,$&$\,$&$\,$&$\,$&$\,$&$\,$&$\,$\\
N&2-is Ht&$\,$&N&2-is Ht&$\,$&N&2-is Ht&$\,$&N&2-is Ht\\
$\,$&$\,$&$\,$&$\,$&$\,$&$\,$&$\,$&$\,$&$\,$&$\,$&$\,$\\
282&130&$\,$&283&271&$\,$&302&176&$\,$&317&142\\
347&122&$\,$&393&127&$\,$&394&97.5&$\,$&427&111\\
457&184&$\,$&478&262&$\,$&488&239&$\,$&514&222\\
518&97.8&$\,$&537&225&$\,$&548&172&$\,$&557&352\\
562&114&$\,$&563&415&$\,$&569&102&$\,$&613&453\\
618&150&$\,$&622&218&$\,$&632&305&$\,$&642&140\\
653&710&$\,$&658&110&$\,$&669&93.2&$\,$&719&95\\
723&256&$\,$&730&108&$\,$&733&434&$\,$&734&152\\
742&107&$\,$&743&131&$\,$&745&112&$\,$&752&109\\
787&346&$\,$&788&555&$\,$&808&495&$\,$&811&222\\
821&221&$\,$&822&151&$\,$&829&230&$\,$&852&352\\
858&154&$\,$&859&102&$\,$&862&318&$\,$&865&456\\
877&1290&$\,$&886&225&$\,$&893&160&$\,$&907&276\\
908&265&$\,$&913&217&$\,$&920&205&$\,$&922&318\\
932&269&$\,$&933&485&$\,$&937&130&$\,$&962&99.2\\
977&638&$\,$&982&448&$\,$&992&219&$\,$&994&115
\end{tabular}
\end{center}
\end{table}

\subsection{A Related Problem}
Immediately after the discussion of concordant forms in Dickson, there is a short mention of the work of Lucas, Gerono and Pepin on
\begin{equation*}
t^2+u^2=2v^2 \hspace{2cm} t^2+2u^2=3w^2
\end{equation*}
which we can generalize to
\begin{equation}
x^2+Ny^2=(N+1)z^2 \hspace{2cm} x^2+(N+1)y^2=(N+2)w^2
\end{equation}

This has the obvious solutions $x=\pm 1, y=\pm 1$, but are there others, and for which values of $N$?

Proceeding as in the previous section, define $Y=z/x$ and $X=y/x$, so the first equation is $Y^2=NX^2+1$, which has the solution
$X=0, Y=1$. The line $Y=1+kX$ meets the curve in one further point where
\begin{equation*}
X=\frac{k^2(N+1)-2k(N+1)+N}{k^2(N+1)-N}
\end{equation*}

The second quadric is $1+(N+1)X^2=(N+2)\Box$ which gives the rational quartic
\begin{equation}
D^2=(N+1)^2(N+2)^2k^4-4(N+2)(N+1)^3k^3+
\end{equation}
\begin{equation*}
2(N+1)(N+2)(3N^2+4N+2)k^2-4N(N+2)(N+1)^2k+N^2(N+2)^2
\end{equation*}

Define $V=D/(N+1)(N+2)$ and $U=k/(N+1)(N+2)$ giving
\begin{equation}
V^2=U^4-4(N+1)^2U^3+2(N+1)(N+2)(3N^2+4N+2)U^2-
\end{equation}
\begin{equation*}
4N(N+2)^2(N+1)^3U+N^2(N+1)^2(N+2)^4
\end{equation*}

This quartic can be transformed to the elliptic curve
\begin{equation}
G^2=H^3+(1-N)(N^2+3N+2)H^2-N(N+1)^2(N+2)^2H
\end{equation}
which can be written
\begin{equation}
G^2=H(H+N^2+3N+2)(H-N(N^2+3N+2))
\end{equation}

So far, this is similar to previous analyses. The interesting change comes when we derive the reverse transformation
\begin{equation}
k=\frac{G+(N+1)^2H}{(H+N+1)(N+1)(N+2)}
\end{equation}

Years of experience of computing elliptic curves and transformations have shown that poles of the transformation function
very often give interesting information. We find $H=-(N+1)$ gives $G=\pm (N+1)^3$, which is not a torsion point.

Adding $(-(N+1),(N+1)^3)$ to $(0,0)$, give $H=N(N+1)(N+2)^2$ and $G=\pm N(N+1)^3(N+2)^2$. The positive value gives
\begin{equation*}
k=\frac{2N(N+2)}{N^2+3N+1}
\end{equation*}
which gives
\begin{equation}
X=\frac{y}{x}=\frac{N^4+2N^3-5N^2-14N-7}{3N^4+14N^3+21N^2+10N-1}
\end{equation}
which gives a parametric solution to the original problem. Adding other torsion points just changes the signs on $x$ and $y$.
Doubling $(-(N+1),(N+1)^3)$ gives a parametric solution of degree $12$ for both $x,y$. We can, thus, generate an infinite
number of parametric solutions of increasing degree.

\subsection{$x^2+exy+fy^2=\Box$ and $x^2+gxy+hy^2=\Box$}
As we saw in the previous sections, we can sometimes have non-trivial solutions for all values of a parameter and sometimes
for only certain values. In the following section, we consider
\begin{equation}
x^2+exy+fy^2=\Box \hspace{2cm} x^2+gxy+hy^2=\Box
\end{equation}
where $e,f,g,h, \in \mathbb{Z}$. We look for conditions on these parameters which give parametric solutions. To prevent
each quadric being a square $\forall \, x,y$, we assume $e^2 \ne 4f$ and $g^2 \ne 4h$.

In $x^2+exy+fy^2=z^2$, define $t=z/x, s=y/x$, so the we have $t^2=fs^2+es+1$. This has the
solution $(0,1)$, and consider the line $t=1+ms$, this meets the curve at $(0,1)$ and
\begin{equation}
s=\frac{y}{x} = \frac{e-2m}{m^2-f} = \frac{q(eq-2p)}{p^2-fq^2}
\end{equation}
if we set $m=p/q$.

Taking $x=p^2-f q^2$ and $y=q(eq-2p)$, and substituting into $x^2+gxy+hy^2=w^2$, gives
\begin{equation}
w^2=p^4-2g p^3q+(eg-2f+4h)p^2q^2+2(fg-2eh)p q^3+(e^2h-efg+f^2)q^4
\end{equation}

This quartic is birationally equivalent to an elliptic curve and using the formulae in Mordell \cite{mord},
we get the elliptic curve $E$ to be
\begin{equation}\label{efgheq}
v^2=u^3+2(2(f+h)-eg)u^2+(e^2-4f)(g^2-4h)u
\end{equation}
with the relationship
\begin{equation}
\frac{p}{q}=\frac{e(g^2-4h)-gu-v}{2(g^2-4h-u)}
\end{equation}
The discriminant of this curve is
\begin{equation}
\Delta=256(e^2-4f)^2(g^2-4h)^2((f-h)^2+(e-g)(eh-fg))
\end{equation}
so the curve is non-singular if $(f-h)^2+(e-g)(eh-fg) \ne 0$.

The elliptic curve $E$ clearly has one point of order $2$, at $(0,0)$, which gives $p/q=e/2$ and so $y=0$.
Taking out the common factor $u$ in the right-hand-side leaves a quadratic with discriminant
\begin{equation*}
4((f-h)^2+(e-g)(eh-fg))
\end{equation*}
which we have just assumed to be non-zero.

Thus, if $e=g$ or $eh=fg$, $E$ has certainly $3$ points of order $2$. If $e=g$ the curve $E$ is
\begin{equation*}
v^2=u(u-(e^2-4f))(u-(e^2-4h))
\end{equation*}
whilst if $eh=fg$ the curve $E$ reduces to
\begin{equation*}
v^2=u(u-(eg-4h))(u- \frac{(g^2-4h)e}{g})
\end{equation*}
It is difficult in such a general form to say anything about other torsion points.

If we look at the transformation relationship we see the denominator is zero when
$u=g^2-4h$. Substituting into $E$ we get $v=\pm(e-g)(g^2-4h)$. Adding one of these points to $(0,0)$, gives
$u=e^2-4f, v=\pm(e-g)(e^2-4f)$.

These points are points of order $2$ when $e=g$. When $e \ne g$, let $P=(g^2-4h,(e-g)(g^2-4h))$, then the u-coordinate of $2P$ is
\begin{equation}
u= \frac{(e^2-4f-g^2+4h)^2}{4(e-g)^2}
\end{equation}

For almost all $(e,f,g,h)$ this will be a non-integer fraction and hence $P$ will be a point of infinite order. The only problems are
when $e^2-4f=g^2-4h$, when $u=0$, so that $P$ is a point of order $4$.

Now, assuming $e \ne g$ and $e^2-4f \ne g^2-4h$, we put $u=e^2-4f$ and $v=-(e-g)(e^2-4f)$ into the transformation formulae and we find
\begin{equation}
x = e^4-4ge^3-2(4f-3g^2+4h)e^2+
\end{equation}
\begin{equation*}
4g(4f-g^2+4h)e+16f^2-8f(g^2+4h)+(g^2-4h)^2
\end{equation*}
and
\begin{equation}
y=8(e-g)(e^2-4f-g^2+4h)
\end{equation}
with $y$ clearly non-zero.
We can extend the elliptic curve computations and essentially generate an infinite set of solutions to the problem.

Summarizing, we seem to have an infinite number of solutions unless $e=g$ or $e^2-4f=g^2-4h$.

\subsection{$|e|=|g|=1$}
The simplest pair of quadrics, in this group, to consider are probably
\begin{equation}\label{dd100}
x^2+xy+y^2=\Box \hspace{2cm} x^2+xy+Ny^2=\Box
\end{equation}
where $e=g=1$ and $f=1, \, h=N$. Clearly $x=\pm k, y=0$ is a solution, but are there others?

We can parameterize the first quadratic by
\begin{equation*}
\frac{y}{x}=\frac{1-2m}{m^2-1}
\end{equation*}
and substituting into the second gives that $m$ must satisfy the quartic relation
\begin{equation}
D^2=m^4-2m^3+(4N-1)m^2+2(1-2N)m+N
\end{equation}
with $D,m \in \mathbb{Q}$.

Mordell's method \cite{mord} gives the equivalent elliptic curve as
\begin{equation}\label{ecdd100}
v^2=u^3+2(2N+1)u^2+(12N-3)u=u(u+3)(u+4N-1)
\end{equation}
with
\begin{equation*}
m=\frac{v+u+4n-1}{2(u+4n-1)}
\end{equation*}

The curve \eqref{ecdd100} has
\begin{equation*}
\Delta=2^8\,3^2\,(N-1)^2\,(4N-1)^2
\end{equation*}
and is singular at $N=1$, and, in general, has a torsion subgroup isomorphic to $\mathbb{Z}/ 2\mathbb{Z} \oplus \mathbb{Z}/2\mathbb{Z}$
with torsion points $(0,0)$, $(-3,0)$ and $(1-4N,0)$. None of these torsion points give a non-trivial solution.

All the curves have conductor less than $100$ million. The distribution of ranks is given
in the Table at the end of the section.
All values of $N$, with predicted rank greater than zero, have a solution which was found fairly easily.

As an example, $N=38$ gives the elliptic curve $v^2=u^3+154u^2+453u$ which has $(961/1764,1267435/74088)$ as a point of infinite order. This gives
$m=3973/7140$ and $y/x=5754840/35194871$ with
\begin{equation*}
x^2+x\,y+y^2=38397109^2 \hspace{2cm} x^2+x\,y+38y^2=51958741^2
\end{equation*}
with $x=35194871$ and $y=5754840$.

Several parametric solutions can be found. One is $N=3k^2-2, u=4k^2-3, v=4k(4k^2-3)$, which leads to $x=k^2+2k-3$ and $y=-4k$ and
\begin{equation*}
x^2+xy+y^2=(k^2+3)^2 \hspace{2cm} x^2+xy+Ny^2=(7k^2-3)^2
\end{equation*}

\vspace{1cm}

One stage more complex than \eqref{dd100} is the pair of quadrics
\begin{equation}
x^2+xy+Ny^2=\Box \hspace{2cm} x^2+xy-Ny^2=\Box
\end{equation}
so that $e=g=1$, $f=N$ and $h=-N$.

This gives the elliptic curve
\begin{equation}\label{ecdd110}
v^2=u^3-2u^2+(1-16N^2)u=u(u-(1-4N))(u-(1+4N))
\end{equation}
with
\begin{equation}
\frac{y}{x}=\frac{1-2m}{m^2-N} \hspace{2cm} m=\frac{v+u-(1+4N)}{2(u-(1+4N))}
\end{equation}

This curve has discriminant
\begin{equation*}
\Delta=2^{10}\,N^2\,(4N+1)^2\,(4N-1)^2
\end{equation*}
so is non-singular for $N \in [1,999]$.

The torsion subgroup seems to be isomorphic to $\mathbb{Z}/ 2\mathbb{Z} \oplus \mathbb{Z}/2\mathbb{Z}$
with torsion points $(0,0)$, $(1+4N,0)$ and $(1-4N,0)$, which do not give solutions with $y \ne 0$.

As an example, $N=21$ gives the rank $1$ curve $v^2=u^3-2u^2-7055u$, which has a point of infinite order at
$(-4067/81\, , \,343952/729)$. This gives $y/x=-18648/103945$ and
\begin{equation*}
x^2+x\,y+21y^2=127157^2 \hspace{1cm} x^2+x\,y-21y^2=39541^2
\end{equation*}
with $x=103945$ and $y=-18648$.

There are a few parametric solutions to be found. One such comes from $N=6k^2+6k+2, \, u=-(8k^2+8k+3), \, v=4(2k+1)(8k^2+8k+3)$.
This gives $x=5k^2+6k+2, \, y=-(2k+1)$ and
\begin{equation*}
x^2+xy+Ny^2=(7k^2+7k+2)^2 \hspace{2cm} x^2+xy-Ny^2=(k^2+k)^2
\end{equation*}

Very recently, the set of solutions for $N \in [1,999]$ was completed.

\vspace{1cm}

The final pair of quadrics considered were
\begin{equation}
x^2+xy+Ny^2=\Box \hspace{2cm} x^2-xy+Ny^2=\Box
\end{equation}
so that $e=1$, $g=-1$, and $f=h=N$ so that $e^2-4f=g^2-4h$.

We have the equivalent elliptic curve
\begin{equation}\label{ecdd10}
v^2=u^3+2(4N+1)u^2+(4N-1)^2u
\end{equation}
with
\begin{equation}
\frac{y}{x}=\frac{1-2m}{m^2-N} \hspace{2cm} m=\frac{v-u+4N-1}{2(u+4N-1)}
\end{equation}

$N=12$ gives the curve $v^2=u^3+98u^2+2209u$ which has rank $1$ and  a point of infinite order
with $u=-15228/289$ and $v=476298/4913$. This gives $m=-4111/238$ and $y/x=42840/345119$. Setting
$x=345119$ and $y=42840$ gives
\begin{equation*}
x^2+x\,y+12y^2=394861^2 \hspace{2cm} x^2-x\,y+12y^2=355451^2
\end{equation*}

A complete set of solutions for $N \in [1,999]$ is available from me, if wanted.

\begin{table}[H]
\begin{center}
\caption{Rank summary for $|e|=|g|=1$ curves}
\begin{tabular}{lcccrr}
$\,$&$\,$&$\,$&$\,$ &$\,$&$\,$\\
Curve&Rank$=0$&Rank$=1$&Rank$\ge 2$&Ave. ht&Max. ht.\\
$\,$&$\,$&$\,$&$\,$ &$\,$&$\,$ \\
\eqref{ecdd100}&409 &491 &98&4.1&36.3\\
\eqref{ecdd110}&427 &481 &91&9.5&132.5\\
\eqref{ecdd10}&443 &511 &45&15.7&196.6
\end{tabular}
\end{center}
\end{table}

\subsection{$|e|=|g|=N$}

The simplest pair of quadrics of this form is
\begin{equation}
x^2+Nxy=z^2 \hspace{2cm} x^2+Nxy+y^2=w^2
\end{equation}
and look for solutions with $x\,y\,z\,w \ne 0$.

Considering the first quadric, define $t=y/x$ and $s=z/x$ so that
\begin{equation*}
s^2=nt+1
\end{equation*}
which has the obvious solution $t=0, s=1$. The line $s=1+mt$ meets the curve at a second point
\begin{equation*}
t=\frac{y}{x}=\frac{N-2m}{m^2}
\end{equation*}

The second quadric is $t^2+Nt+1=\Box$, so that $m$ must satisfy the rational quartic
\begin{equation}
D^2=m^4-2Nm^3+(N^2+4)m^2-4Nm+N^2
\end{equation}

This quartic can be transformed, by Mordell's method, to the equivalent elliptic curve
\begin{equation}\label{ecdd50}
v^2=u(u-N^2)(u-N^2+4)
\end{equation}
with
\begin{equation}
m=\frac{v+Nu-N(N^2-4)}{2(u+4-N^2)}
\end{equation}

The elliptic curve has discriminant
\begin{equation*}
\Delta=256N^4(N+2)^2(N-2)^2
\end{equation*}
so is singular for $N=2$ in the range $[1,999]$. There are $3$ finite points of order $2$ at $u=0$, $u=N^2$ and $u=N^2-4$. There
are also $4$ points of order $4$ at $(\,N(N+2), \pm 2N(N+2)\,)$ and $(\,N(N-2), \pm 2N(N-2)\,)$. There will be points of order $8$
if $N(N \pm 2)$ is an integer square, which cannot happen. Thus, the torsion subgroup is isomorphic to
$\mathbb{Z}/2\mathbb{Z} \oplus \mathbb{Z}/4\mathbb{Z}$. All the torsion points give $x\,y\,z\,w =0$. Thus we would need a curve with rank
greater than zero.

For example, $N=24$ gives a curve of rank $1$ with a point of infinite order $(176,5280)$. This gives $m=16/3$ and $y/x=15/32$ and
\begin{equation*}
32^2+24\times32\times15=112^2 \hspace{1.5cm} 32^2+24\times32\times15+15^2=113^2
\end{equation*}

\begin{table}[H]
\begin{center}
\caption{Rank summary for $|e|=|g|=N$ curves}
\begin{tabular}{lcccrr}
$\,$&$\,$&$\,$&$\,$ &$\,$&$\,$\\
Curve&Rank$=0$&Rank$=1$&Rank$\ge 2$&Ave. ht&Max. ht.\\
$\,$&$\,$&$\,$&$\,$ &$\,$&$\,$ \\
\eqref{ecdd50}&411 &506 &81&9.5&71.6\\
\eqref{ecdd3}&450 &478 &70&15.0&308.9\\
\eqref{ecdd120}&461 &493 &44&26.1&254.8\\
\eqref{ecdd20}&467&458&73&36.4&503.4\\
\eqref{ecdd60}&423&473&102&18.8&329.8
\end{tabular}
\end{center}
\end{table}

From Table $5.6$ we have small heights, in general, and we have computed solutions for all $N \in [1,999]$ predicted to
give answers.

\vspace{1cm}

The next pair considered was
\begin{equation}
x^2+Nxy+y^2=\Box \hspace{2cm} x^2-Nxy+y^2=\Box
\end{equation}
which gives the elliptic curve
\begin{equation}\label{ecdd3}
v^2=u^3+2(N^2+4)u^2+(N^2-4)^2u=u(u+(N+2)^2)(u+(N-2)^2)
\end{equation}
with
\begin{equation}
\frac{y}{x}=\frac{N-2m}{m^2-1} \hspace{1.5cm} m=\frac{v-Nu+N(4-N^2)}{2(u+4-N^2)}
\end{equation}

For example, $N=35$ gives $v^2=u^3+2458u^2+1490841u$. If $u=40071$ we have $v = \pm 8267280$. Taking the positive value gives $m=439/5$
and hence $y/x=-95/5208$. Defining $x=5208$ and $y=95$ gives
\begin{equation*}
x^2+35xy+y^2=6667^2 \hspace{2cm} x^2-35xy+y^2=3133^2
\end{equation*}

The elliptic curve \eqref{ecdd3} has $3$ finite points of order $2$, at $(0,0)$, $(-(N-2)^2,0)$ and $(-(N+2)^2,0)$. There are also
points of order $4$ at $(N^2-4,\pm 2N(N^2-4))$ and $(4-N^2,\pm 4(N^2-4))$. There can only be points of order $8$ if $N^2-4$ is an
integer square, which does not happen. Thus the torsion subgroup is isomorphic to $\mathbb{Z}/2\mathbb{Z} \oplus \mathbb{Z}/4\mathbb{Z}$.
None of the torsion points lead to suitable values of $m$, so we need curves with rank at least $1$.

The curves have discriminant
\begin{equation*}
\Delta=1024N^2(N+2)^4(N-2)^4
\end{equation*}
so, in $[1,999]$, are only singular at $N=2$.

With the help of Randall Rathbun, I managed to reduce the unsolved set of values to $N=809$ with estimated height $308.9$. The 2-isogenous
curve had height half this value, but this was still way beyond my software's range. The conductor was too large to make a Heegner point
computation feasible. Thankfully, Tom Fisher at Cambridge used this problem as an example in \cite{fish}, and I now have a complete
set of solutions in $[1,999]$.

\vspace{1cm}

A minor variation on the last problem is
\begin{equation*}
x^2+Nxy+y^2=\Box \hspace{2cm} x^2+Nxy-y^2=\Box
\end{equation*}
which gives the elliptic curve
\begin{equation}\label{ecdd120}
v^2=u^3+(N^2+12)u^2+8(N^2+4)u=u(u+8)(u+N^2+4)
\end{equation}
with
\begin{equation}
\frac{y}{x}=\frac{N-2m}{m^2-1} \hspace{1.5cm} m=\frac{v+Nu}{2u}
\end{equation}

For example, $N=17$ gives the curve $v^2=u^3+301u^2+2344u$ which has rank $1$ and rational points $(20449/16, \pm 3252249/64)$. Taking
the negative value gives $m=-13019/1144$ and $y/x=912912/2950625$, with
\begin{equation*}
x^2+17xy+y^2=7438537^2 \hspace{1cm} x^2+17xy-y^2=7325641^2
\end{equation*}
with $x=2950625$ and $y=912912$.

The curve \eqref{ecdd120} has $3$ finite points of order $2$ at $(0,0)$, $(-8,0)$ and $(-N^2-4,0)$. There will only be
points of order $4$ when $8(N^2+4)$ is an integer square, which cannot happen. There is always the possibility of points of order $3$
to be considered. These will be points of inflexion of the curve so will be integer roots of
\begin{equation*}
3u^4+4(N^2+12)u^3+48(N^2+4)u^2-64(N^2+4)^2=0
\end{equation*}
There are at least $2$ real roots of this quartic, and I think there are never integer roots, but cannot show it.
Thus, I assume the torsion subgroup is isomorphic to $\mathbb{Z}/2\mathbb{Z} \oplus \mathbb{Z}/2\mathbb{Z}$. None of these torsion
points lead to non-trivial solutions, so we need curves of rank at least $1$.

The curve has discriminant
\begin{equation*}
\Delta=1024(N+2)^2(N-2)^2(N^2+4)^2
\end{equation*}
so is singular at $N=\pm 2$, as we might expect.

The unsolved values of $N$ are given in Table $5.7$.
For $N=589, 631, 685, 813$, the heights for the 2-isogenous curve
\begin{equation*}
g^2=h^3-2(N^2+12)h^2+(N^2-4)^2h
\end{equation*}
are half the given value. For all the other values, the 2-isogenous curve has double the given height.

\begin{table}\label{dd120uns}
\begin{center}
\caption{Unsolved values and heights for \eqref{ecdd120}}
\begin{tabular}{llllllrrrrr}
$\,$&$\,$&$\,$&$\,$&$\,$&$\,$&$\,$&$\,$&$\,$&$\,$&$\,$\\
N&Ht&$\,$&N&Ht&$\,$&N&Ht&$\,$&N&Ht\\
$\,$&$\,$&$\,$&$\,$&$\,$&$\,$&$\,$&$\,$&$\,$&$\,$&$\,$\\
357&128&$\,$&507&142&$\,$&540&99.8&$\,$&589&223\\
624&92.8&$\,$&631&197&$\,$&645&145&$\,$&685&254\\
699&96.5&$\,$&748&93.5&$\,$&792&149&$\,$&809&100\\
813&262&$\,$&825&247&$\,$&837&166&$\,$&917&107\\
931&122&$\,$&951&116&$\,$&981&96.1&$\,$&993&174
\end{tabular}
\end{center}
\end{table}

\vspace{1cm}

Another problem which attracted my attention is
\begin{equation*}
x^2+Nxy+(N+1)y^2=\Box \hspace{2cm} x^2+Nxy+(N-1)y^2=\Box
\end{equation*}
with $xy \ne 0$.

This gives the elliptic curve
\begin{equation}\label{ecdd20}
v^2=u^3+(N^2-4N+12)u^2+8(N-2)^2u=u(u+8)(u+(N-2)^2)
\end{equation}
with
\begin{equation}
\frac{y}{x}=\frac{N-2m}{m^2-N-1} \hspace{1.5cm} m=\frac{v+Nu}{2u}
\end{equation}

As an example, $N=13$ gives the curve $v^2=u^3+129u^2+968$. The points $(313^2/56^2, \pm 75866505/56^3)$ lie on the curve.
The positive v-value gives $m=470249/35056$ giving $x=1804683169$ and $y--150390240$, with
\begin{equation*}
x^2+13xy+14y^2=212685041^2 \hspace{1.5cm} x^2+13xy+12y^2=691441^2
\end{equation*}

As with the previous example, I assume the torsion subgroup is isomorphic to $\mathbb{Z}/2\mathbb{Z} \oplus \mathbb{Z}/2\mathbb{Z}$,
with finite torsion points $(0,0)$, $(-8,0)$ and $(-(N-2)^2,0)$, none of which give a non-trivial answer.

The curve has discriminant
\begin{equation*}
\Delta=1024(N-2)^4(N^2-4N-4)^2
\end{equation*}
so is only singular at $N=2$.

There are still $19$ unsolved values of $N$. Unlike other problems, the unknown values are almost
evenly distributed in cases where the original curve or the 2-isogenous curve
\begin{equation*}
g^2=h^3-2(N^2-4N+12)h^2+(N^2-4N-4)^2h
\end{equation*}
gives the smaller height estimate. Table $5.x$ gives the data, with the asterisk indicating values of $N$
where the 2-isogenous curve gives the smaller height.

\begin{table}\label{dd20uns}
\begin{center}
\caption{Unsolved values and heights for \eqref{ecdd20}}
\begin{tabular}{llllllrrrrr}
$\,$&$\,$&$\,$&$\,$&$\,$&$\,$&$\,$&$\,$&$\,$&$\,$&$\,$\\
N&Ht&$\,$&N&Ht&$\,$&N&Ht&$\,$&N&Ht\\
$\,$&$\,$&$\,$&$\,$&$\,$&$\,$&$\,$&$\,$&$\,$&$\,$&$\,$\\
403&98.4&$\,$&516&157&$\,$&543&96.8*&$\,$&615&135*\\
628&104&$\,$&663&179*&$\,$&676&222&$\,$&697&108*\\
813&103*&$\,$&823&97.3*&$\,$&831&252*&$\,$&841&181\\
861&91.2*&$\,$&883&98.8&$\,$&901&148&$\,$&907&90.2*\\
931&109*&$\,$&939&112&$\,$&940&106&$\,$&$\,$&$\,$
\end{tabular}
\end{center}
\end{table}

\vspace{1cm}

The final problem I include in this section is
\begin{equation*}
x^2+2Nxy+Ny^2=\Box \hspace{2cm} x^2-2Nxy+Ny^2=\Box
\end{equation*}

This gives the elliptic curve
\begin{equation}\label{ecdd60}
v^2=u^3+2N(N+1)u^2+N^2(N-1)^2u
\end{equation}
with
\begin{equation}
\frac{y}{x}=\frac{2(N-m)}{m^2-N} \hspace{1.5cm} m=\frac{v-Nu-N^2(N-1)}{u-N(N-1)}
\end{equation}

For example, $N=42$ gives the curve $v^2=u^3+3612u^2+2965284u$. The points $(-35574/25, \pm 1823976/125)$ lie on the curve.
The positive v=value gives $m=-42/65$ and hence $x=4183$ and $y=8580$, ignoring signs. We have
\begin{equation*}
x^2+84xy+42y^2=78157^2 \hspace{2cm} x^2-84xy+42y^2=9727^2
\end{equation*}

The curves have discriminant
\begin{equation*}
\Delta=256N^7(N-1)^4
\end{equation*}
so are only singular at $N=0,1$.

For general $N$, there is only one point of order $2$, at $(0,0)$. There are, however, points of order $4$ at $(N(N-1), \pm 2N^2(N-1))$.
Investigating the conditions for $3$ points of order $2$, we find this occurs iff $N=M^2$ with $M$ an integer. For this case we have
$(0,0)$, $(-M^2(M+1)^2,0)$ and $(-M^2(M-1)^2,0)$ of order $2$, with $(M^2(M^2-1), \pm 2M^4(M^2-1))$ and $(-M^2(M^2-1), \pm 2M^3(M^2-1))$
of order $4$. Points of order $8$ would require $M^2(M^2-1)=\Box$ which cannot happen if $M>1$, so if $N=M^2$ the torsion
subgroup is isomorphic to $\mathbb{Z}/2\mathbb{Z} \oplus \mathbb{Z}/4\mathbb{Z}$. For $N \neq M^2$, there is the possibility of
points of order $3$ but I conjecture these cannot happen, and, in this case, the torsion subgroup is isomorphic to $\mathbb{Z}/4\mathbb{Z}$.

The points of order $4$ allow us to find a 4-isogenous curve
\begin{equation*}
g^2=h^3+2N(1-2N)h^2+N^2h
\end{equation*}

There are several possible parametric solutions. From $N=4-k^2, \, u=k^2-4, \, v=k(k^2-4)^2$, we find $x=k(k^2-2), \, y=2$ with
\begin{equation*}
x^2+2Nxy+Ny^2=A^2 \hspace{1.5cm} x^2-2Nxy+Ny^2=B^2
\end{equation*}
where
\begin{equation*}
A=k^3-2k^2-4k+4 \hspace{3cm} B=k^3+2k^2-4k-4
\end{equation*}

In fact, the vast majority of parametric solutions found come from the fact that $u=-N$ gives a point when $4-N=\Box$, which has an infinite
number of solutions.

The, as yet, unsolved values are given in Table $5.9$. In all cases the height estimate from the original curve is smaller than those
from either the 2-isogenous curve or the 4-isogenous curve.

\begin{table}\label{dd60uns}
\begin{center}
\caption{Unsolved values and heights for \eqref{ecdd60}}
\begin{tabular}{llllllrrrrr}
$\,$&$\,$&$\,$&$\,$&$\,$&$\,$&$\,$&$\,$&$\,$&$\,$&$\,$\\
N&Ht&$\,$&N&Ht&$\,$&N&Ht&$\,$&N&Ht\\
$\,$&$\,$&$\,$&$\,$&$\,$&$\,$&$\,$&$\,$&$\,$&$\,$&$\,$\\
398&99.6&$\,$&422&149&$\,$&614&142&$\,$&642&141\\
662&251&$\,$&674&153&$\,$&733&147&$\,$&734&134\\
787&115&$\,$&842&119&$\,$&857&127&$\,$&917&125\\
947&119&$\,$&978&101&$\,$&$\,$&$\,$&$\,$&$\,$&$\,$
\end{tabular}
\end{center}
\end{table}

\subsection{Other Problems with Pairs of Quadrics}
It is, of course, perfectly possible to have problems which are reasonable to ask but do not conform to the
format of the previous sections. Each problem must be considered on its' own merits.

The first such problem is
\begin{equation}
x^2+Nxy=z^2 \hspace{2cm} y^2-Nxy=w^2
\end{equation}
where $z\,w \ne 0$.

For the first equation define $s=z/x$ and $t=y/x$, so $s^2=Nt+1$, which has a solution $t=0, s=1$. The line $s=1+mt$ meets the curve
again where
\begin{equation*}
t=\frac{N-2m}{m^2}
\end{equation*}

Using this in the second quadric, we must have $m$ satisfy the rational form
\begin{equation}
D^2=2Nm^3+(4-N^2)m^2-4Nm+N^2
\end{equation}
which is an elliptic curve immediately. Define $D=v/(2N)$ and $m=u/(2N)$ giving
\begin{equation}
v^2=u^3+(4-N^2)u^2-8N^2u+4N^4
\end{equation}

The right-hand side cubic is zero when $u=N^2$, so defining $u=h+N^2$ gives
\begin{equation}\label{ecdd40}
v^2=h^3+2(N^2+2)h^2+N^4h
\end{equation}
and we have
\begin{equation}
\frac{y}{x}=\frac{-4Nh}{(h+N^2)^2}
\end{equation}

The curve has discriminant
\begin{equation*}
\Delta=256N^8(N^2+1)
\end{equation*}
so is non-singular in $[1,999]$. The point $(0,0)$ is the only point of order $2$ as more would require $N^2+1$
to be an integer square. $(-N^2, \pm 2N^2)$ are points of order $4$, and since $-N^2<0$ there cannot be points of order $8$.
This means that the torsion subgroup will be isomorphic to $\mathbb{Z}/4\mathbb{Z}$ or possibly $\mathbb{Z}/12\mathbb{Z}$.
The latter would require a point of order $3$, which would be a point of inflexion satisfying
\begin{equation*}
3x^4+8(N^2+2)x^3+6N^4x^2-N^8=0
\end{equation*}

Numerical evidence is that such a point does not exist, so we assume that there are only $3$ finite torsion points, none of which
give a non-trivial solution. We, thus, need points of infinite order.

For $N=13$, the curve is $v^2=h^3+342h^2+28561h$, which has a point of infinite order with $h=26896/1521$, which gives
$x=80624763025$ and $y=-2127258432$, and
\begin{equation*}
x^2+13x\,y=65350793585^2 \hspace{2cm} y^2-13x\,y=47266810968^2
\end{equation*}

The BSD calculations gave large heights for several curves with $N \in [1,999]$. To find generators, we had to make extensive use
of the fact that the points of order $4$ give a 4-isogeny to another elliptic curve.
The formulae underlying this are given in an appendix, and, for this problem, give the curve
\begin{equation}
s^2=t^3+8(N^2-1)t^2+16(N^2+1)^2t
\end{equation}

The interesting fact is that this curve also has points of order $4$ when $t=4(N^2+1)$, which allows us to derive a further
4-isogeny to the curve
\begin{equation}
j^2=f^3-32(2N^2+1)f^2+256f
\end{equation}

It might be hoped that there was a cascade downwards of heights, along this sequence of curves, but this never seems to happen.
The unsolved values of $N$ are given in Table $5.10$, where the heights are either from the original curve or the first
4-isogenous curve (with the asterisk).

\begin{table}\label{dd40uns}
\begin{center}
\caption{Unsolved values and heights for \eqref{ecdd40}}
\begin{tabular}{llllllrrrrr}
$\,$&$\,$&$\,$&$\,$&$\,$&$\,$&$\,$&$\,$&$\,$&$\,$&$\,$\\
N&Ht&$\,$&N&Ht&$\,$&N&Ht&$\,$&N&Ht\\
$\,$&$\,$&$\,$&$\,$&$\,$&$\,$&$\,$&$\,$&$\,$&$\,$&$\,$\\
263&106*&$\,$&367&175*&$\,$&373&145*&$\,$&389&166*\\
433&145*&$\,$&466&279*&$\,$&503&102*&$\,$&545&114*\\
547&149*&$\,$&559&99.5*&$\,$&587&105*&$\,$&617&130*\\
634&463&$\,$&673&236&$\,$&674&129&$\,$&683&324\\
685&92.6&$\,$&689&149&$\,$&691&93&$\,$&717&110\\
769&112*&$\,$&773&206*&$\,$&778&124*&$\,$&779&135*\\
781&224&$\,$&802&229&$\,$&815&165&$\,$&823&296\\
839&128*&$\,$&859&107*&$\,$&869&131*&$\,$&872&219*\\
878&142*&$\,$&883&484*&$\,$&907&150*&$\,$&921&161*\\
937&752*&$\,$&941&185*&$\,$&942&187*&$\,$&947&92.2*\\
949&117&$\,$&953&125&$\,$&955&217&$\,$&962&216\\
981&109&$\,$&989&97.2&$\,$&997&160&$\,$&998&208
\end{tabular}
\end{center}
\end{table}

\vspace{2cm}

Another problem one might analyse is
\begin{equation*}
x^2+(N-1)y^2=Nz^2 \hspace{1cm} (N-1)x^2+y^2=Nw^2
\end{equation*}
where we look for solutions with $|x| \ne |y|$.

Setting $s=z/x$ and $t=y/x$ transforms the first quadric to $Ns^2=(N-1)t^2+1$ which has a rational solution $s=t=1$. The line $s=1+m(t-1)$
meets the curve again where
\begin{equation*}
t=\frac{N(m-1)^2-1}{N(m^2-1)+1}
\end{equation*}

The second quadric gives $N(t^2+N-1)=\Box$, which means that we have to consider the quartic
\begin{equation*}
y^2=x^4-4x^3-2(N^2-5N+2)x^2+4N(1-N)x+N^2(1-N)^2
\end{equation*}
where $x=mN$.

This can be transformed to an equivalent elliptic curve by Mordell's method, giving
\begin{equation}
v^2=u^3+(N^2-2N+2)u^2+(N-1)^2u=u(u+1)(u+(N-1)^2)
\end{equation}
with
\begin{equation}
m=\frac{v+u+(N-1)^2}{N(u+N-1)}
\end{equation}

The above elliptic curve is, however, just the curve for Leech's problem as discussed previously. All the generators
have been computed, so we just use different transformations. From the Leech data, if $N=8$ the curve $v^2=u^3+50u^2+49u$ has a
point $u=1, v=10$. This gives $m=15/16$ which gives (ignoring signs) $x=1$ and $y=31$ with
\begin{equation*}
1^2+7 \times 31^2=8\times29^2 \hspace{2cm} 31^2+7\times1^2=8\times11^2
\end{equation*}

\vspace{2cm}

Finally, we consider the problem
\begin{equation}
a^2+Nb^2=c^2 \hspace{3cm} Na^2+b^2=d^2
\end{equation}
which was asked about on the math.stackexchange.com web-site by Tito Piezas in February 2016.

The first quadric is common to both the congruent number problem and the concordant number problem, and does seem a natural question.
This quadric can be parameterized by $a=N-k^2$ and $b=2k$, which, substituted into $Na^2+b^2=d^2$ gives the quartic
\begin{equation}
d^2=Nk^4+(4-2N^2)k^2+N^3
\end{equation}

The problem is that this quartic has no obvious rational solution for general $N$ and so cannot be transformed to a specific family
of elliptic curves.
In fact, for some values of $N$, the quartic has no rational solutions at all, such as for $N=2$ and $N=5$.

There is, however, a rational solution when $N=M^2$. The problem then is
\begin{equation*}
a^2+(Mb)^2=c^2 \hspace{3cm} (Ma)^2+b^2=d^2
\end{equation*}
with the quartic
\begin{equation*}
d^2=M^2k^4+(4-2M^4)k^2+M^6
\end{equation*}

We can transform this to the elliptic curve
\begin{equation}
v^2=u(u+1)(u+M^4)
\end{equation}
with
\begin{equation}
k=\frac{v}{M(u+1)}
\end{equation}
and $a=M^2-k^2$ and $b=2k$. This gives (usually) rational solutions which can be scaled up to integer ones.

The elliptic curve has $3$ finite torsion points of order $2$ - namely $(0,0)$, $(-1,0)$ and $(-M^4,0)$. There are also $4$ points of order $4$
\begin{equation*}
(\,M^2 \, , \, \pm M^2(M^2+1) \, ) \hspace{2cm} (\,-M^2 \, , \, \pm M^2(M^2-1) \, )
\end{equation*}
Since the first pair of these points has a square u-coordinate, it is possible that there could be points of order $8$. Investigating this
shows that this only occurs if $M^2+1=\Box$, which is impossible for integer $M \ne 0$. Thus the torsion subgroup is isomorphic to
$\mathbb{Z}/2\mathbb{Z} \oplus \mathbb{Z}/4\mathbb{Z}$. None of these torsion points lead to a non-trivial solution of the original problem
so we need to find curves with rank at least $1$.

As an example, $M=17$ gives a curve of rank $1$ with a point of infinite order $u=760384/12321$ and $v=24598269800/1367631$. This
gives $k=2323880/137751$ which leads to
\begin{equation*}
a=834\,5442\,7889 \hspace{2cm} b=6402\,3358\,7760
\end{equation*}

\newpage

\section{Cubic Representations of Integers}

\subsection{Two Cubes Problem}
Apart from the congruent number problem, this is (probably) the most studied problem leading to an elliptic curve. The problem's importance
is seen from the fact that it takes up several pages in Chapter XXI of Dickson's History \cite{dick1}.

{\bf Given an integer $N$, can we find two rational numbers $a,b$ such that $N=a^3+b^3$? }

Writing $a=X/Z$ and $b=Y/Z$, this is equivalent to
the problem of finding three integers $X,Y,Z$ such that
\begin{equation}
N=\frac{X^3+Y^3}{Z^3}
\end{equation}
Note that, if $(a,b)$ is a solution for $N$ then $(-a,-b)$ is a solution for $-N$, so we assume $N > 0$. It is also clear that if $N=N_1N_2^3$
with $N_1, N_2 \in \mathbb{Z}$ then
\begin{equation*}
\left( \frac{a}{N_2} \right)^3 + \left( \frac{b}{N_2} \right)^3 = N_1
\end{equation*}
so we can assume, without loss of generality, that $N$ is cubefree.

Let $c=a-b$, so that
\begin{equation*}
N=a^3+(c-a)^3=3ca^2-3c^2a+c^3
\end{equation*}
so $a$ must satisfy the quadratic equation
\begin{equation}
3ca^2-3c^2a+(c^3-N)=0
\end{equation}

For $a$ to be rational, the discriminant must be a rational square, so there must exist $d \in \mathbb{Q}$ such that
\begin{equation}
d^2=(3c^2)^2-4(3c)(c^3-N)=-3c^4+12cN
\end{equation}

Since $N$ is cubefree, $c \neq 0$, so let $e=1/c$ and thus
\begin{equation}
d^2e^4=12Ne^3-3
\end{equation}
and, multiplying both sides by $(12N)^2$, we have
\begin{equation}
(12Nde^2)^2=(12Ne)^3-432N^2
\end{equation}

Define $g=12Nde^2$ and $h=12Ne$ so we finally have the elliptic curve
\begin{equation}\label{ectwocub}
g^2=h^3-432N^2
\end{equation}

We have
\begin{equation*}
c=\frac{12N}{h} \hspace{2cm} d = \frac{12Ng}{h^2}
\end{equation*}
and it is standard algebra to show
\begin{equation}
a=\frac{36N+g}{6h} \hspace{2cm} b=\frac{36N-g}{6h}
\end{equation}

The first rank-1 curve is at $N=6$ and it is an easy search to find the point $(28,80)$. This gives $a=37/21$ and $b=17/21$.

The torsion of the curve \eqref{ectwocub} is fairly simple. There will be points of order $2$ if $2N^2$ is a cube. Since $N$ is assumed
cubefree, this only happens for $N=2$ when the torsion subgroup is isomorphic to $\mathbb{Z} / 2 \mathbb{Z}$ with the torsion point
at $(12,0)$. Points of inflexion occur when $N^2$ is also a cube, which is only true when $N=1$ and again $x=12$. For all other values
of $N$ there are no torsion points.

Using the BSD conjecture, for $N \in [1,999]$ there are $399$ rank $0$ curves, $490$ with rank $1$ and $110$ with rank at least $2$. The average
height for the rank-1 curves is a lowly $6.39$ with the maximum estimated height being $52.73$ for $N=690$. This latter value, however, is
an overestimate by a factor of $9$. I actually have solutions for all $N$ in $[1,9999]$. At the time of writing, I am expanding this
calculation to eventually cover all $N$ up to $99999$ - I have finished up to $N=19999$.

As with the congruent number problem, so much has been written on the two-cubes problem that I cannot describe it all. I will
describe what are known as Sylvester's Conjectures, which attempt to characterize when a solution exists for $N$ prime.

A descent shows that the rank is $0$ when $N \equiv 2,5 \, (\bmod \, \, 9)$ so no solution exists. The rank is at most $1$ when $N \equiv 4,\,7,\,8 \, (\bmod \, \, 9)$
and the BSD and Parity conjecture both suggest exactly equal to $1$. Noam Elkies has an unpublished proof for the values $4$ and $7$.
For $N \equiv 1 \, (\bmod \, \,9)$, the analysis suggests the rank is even, either $0$ or $2$.

\subsection{Knight's Problem}
In \cite{bgn}, Andrew Bremner, Richard Guy and Richard Nowakowski discussed the problem of (if possible) finding integers $X,Y,Z$ such that
\begin{equation}
N=(X+Y+Z)\left( \frac{1}{X} + \frac{1}{Y} + \frac{1}{Z} \right )
\end{equation}
which they stated was originally proposed by Benjamin Knight.

We have
\begin{equation*}
N=\frac{X^2Y+Y^2Z+Z^2X+XY^2+YZ^2+ZX^2+3XYZ}{X\,Y\,Z}
\end{equation*}

Define $x=X/Z$ and $y=Y/Z$, and we have the quadratic
\begin{equation}
(y+1)x^2+(y^2+(3-N)y+1)x+y(y+1)=0
\end{equation}

If this is to have a rational solution $x$, then the discriminant must be a rational square, so there must exist $d \in \mathbb{Q}$
with
\begin{equation*}
d^2=y^4+2(1-N)y^3+(N^2-6N+3)y^2+2(1-N)y+1
\end{equation*}

Standard methods give that this quartic is birationally equivalent to the elliptic curve
\begin{equation}\label{eckn}
g^2=h^3+(N^2-6N-3)h^2+16Nh
\end{equation}
with
\begin{equation}
y=\frac{g+(N-1)h}{2(h-4N)}
\end{equation}

The elliptic curve has discriminant
\begin{equation}
\Delta=4096N^2(N-9)(N-1)^3
\end{equation}
and so is singular when $N=0,1,9$, which we now exclude from the discussion.

From a point $(h,g)$ on the elliptic curve we can find $y=Y/Z$. It might be thought that we need to slot this $y$-value into
the quadratic for $x$ and solve. It is possible to show, however, that
\begin{equation}
x=\frac{-g+(N-1)h}{2(h-4N)}
\end{equation}
and we summarize as
\begin{equation}
\frac{X,Y}{Z}=\frac{\pm g +(N-1)h}{2(h-4N)}
\end{equation}

The elliptic curve \eqref{eckn} has an obvious point of order $2$ at $(0,0)$. There are $3$ points of order $2$ when $(N-9)(N-1)$ is an
non-zero integer square. This only seems to happen when $N=10$, though I have not tried to prove this!

Points of order $3$ are points of inflexion, which occur where
\begin{equation*}
3h^4+4(N^2-6N-3)h^3+96Nh^2-256N^2=0
\end{equation*}
which factors to
\begin{equation*}
(h-4)(3h^3+4N(N-6)h^2+16N^2h+64N^2)=0
\end{equation*}
so that $(4, \pm 4(N-1)\,)$ are points of order $3$.

Points of order $2$ and $3$ imply points of order $6$, which will occur when
\begin{equation*}
\frac{(h^2-16N)^2}{4(h^3+(N^2-6N-3)h^2+16Nh)}=4
\end{equation*}
which can be easily solved to give $(4N, \pm 4N(N-1)\,)$.

For a point of order $12$ to exist, we would require $N=M^2$ and $M^2+2M-3=\Box$. Numerical tests suggest this cannot happen, but, again,
I have not proven the result.

Thus, we can assume there are $6$ torsion points and the torsion subgroup is isomorphic to $\mathbb{Z} / 6 \mathbb{Z}$. None of the finite
torsion points give acceptable solutions, so we would need the rank of the curve to be greater than zero.

For the $997$ non-singular curves for $N \in [1,999]$, the BSD calculations suggest $439$ rank zero, $494$ rank one and $64$ rank
greater than one.

Bremner, Guy and Nowakowski stated that they had solutions for all but about 5\% of the values in this range. I was able, by about $2005$, to
find all the solutions, with the help of others, especially Mark Watkins of the Magma group in Sydney. I have, in fact, solutions for
all $N$ in $[-999,999]$.

Since the curve \eqref{eckn} has a point of order $2$, there is a standard $2$-isogeny to
\begin{equation*}
v^2=u^3-2(N^2-6N-3)u^2+(N-9)(N-1)^3u
\end{equation*}

Using all the torsion points, we can find a $6$-isogeny to
\begin{equation*}
s^2=t^3-2(N^2+18N-27)t^2+(N-1)(N-9)^3t
\end{equation*}
with
\begin{equation*}
t=\frac{(h^2+4(N-3)h+16N)^2(h^2+(N^2-6N-3)h+16N)}{h(h-4)^2(h-4N)^2}
\end{equation*}

As an example, $N=48$ gives a curve with rank $1$ and the height estimate of the generator is $5.4$. We quickly find $h=1587/191^2, \, g= 42509244/191^3$.
These give $X=-402523, \, Y=200445, \, Z=18972030$.

Parametric solutions also are reasonably easy to find. If $N=-(k-2)(k+1), \, h=4(k+1)$ we have $g=4(k+1)(k^2-k+1)$. These give $X=-k(k-1)$, $Y=1$ and
$Z=k-1$.

\subsection{Bremner-Guy Problems}
In \cite{bg}, Andrew Bremner and Richard Guy considered the problems: if $N$ is a non-zero integer find non-zero integers $X,Y,Z$ such that
\begin{equation}\label{bga}
N=\frac{(X+Y+Z)^3}{X\,Y\,Z}
\end{equation}
or
\begin{equation}\label{bgb1}
N=\frac{X}{Y}+\frac{Y}{Z}+\frac{Z}{X}\
\end{equation}

As pointed out by Bremner and Guy, the latter problem is intimately linked to
\begin{equation}\label{bgb2}
N=\frac{a^3+b^3+c^3}{a\,b\,c}
\end{equation}
since, from a solution $(a,b,c)$, we can define $X=a^2b$, $Y=b^2c$ and $Z=c^2a$, which gives a solution to \eqref{bgb1}.

For problem \eqref{bga}, we have the cubic
\begin{equation*}
X^3+3(Y+Z)X^2+(3Y^2+(6-N)Y\,Z+3Z^2)X+(Y+Z)^3=0
\end{equation*}
and we cannot do very much.

Let $S=X+Y+Z$ and substitute $Z=S-X-Y$ into this cubic giving a quadratic in $Y$
\begin{equation}
N\,X\,Y^2+N\,X(X-S)Y+S^3=0
\end{equation}
which, if we are to have a rational solution, must have a discriminant which is a rational square, so
there must be a rational solution of
\begin{equation*}
d^2=N^2X^4-2N^2X^3S+N^2X^2S^2-4NXS^3
\end{equation*}

Defining $y=Nd/S^2$ and $x=NX/S$ gives the quartic
\begin{equation}
y^2=x^4-2Nx^3+N^2x^2-4N^2x
\end{equation}
which is birationally equivalent to the elliptic curve
\begin{equation}
g^2=h^3+N^2h^2+8N^3h+16N^4=h^3+(Nh+4N^2)^2
\end{equation}
with
\begin{equation}
\frac{X}{S}=\frac{g+N(h+4N)}{2Nh}
\end{equation}

It might be thought that, from a point $(h,g)$ on the elliptic curve, we find $X/S$ and substitute into the quadratic
to get $Y$ and hence $Z$. It can be shown, however, that one of the two roots of the quadratic in $Y$ satisfies
\begin{equation}
\frac{Y}{S}=\frac{-g+N(h+4N)}{2Nh}
\end{equation}

The elliptic curve has discriminant
\begin{equation*}
\Delta=4096N^8(N-27)
\end{equation*}
so is singular when $N=27$, which has the solution $X=Y=Z=1$.

Over the years, I have gradually reduced the number of unsolved values, until in Spring $2013$ I finished all values in $[1,999]$ that
the BSD conjecture predicts will have a solution.
As an example, for $N=41$ we find the point on the curve given by $(-378056/37^2, 14752972/37^3)$. This gives $X=27270901,\, Y=30959144, \, Z=85147693$.

The only parametric solution found comes from $N=-k^2,\, h=4k^2,\, g=8k^3$, which gives $X=-1,\, Y=1,\, Z=k$.

\vspace{1cm}

For \eqref{bgb1}, we have the quadratic in $X$,
\begin{equation*}
Z X^2+Y(Y-N Z)X+Y Z^2=0
\end{equation*}
which must have a rational square discriminant, so
\begin{equation}
D^2=-4YZ^3+NY^2Z^2-2NY^3Z+Y^4
\end{equation}
which gives
\begin{equation*}
\frac{16D^2}{Y^4}=-64\frac{Z^3}{Y^3}+16N\frac{Z^2}{Y^2}-32N\frac{Z}{Y}+16
\end{equation*}

Defining $G=4D/Y^2$ and $H=-4Z/Y$ leads to
\begin{equation}\label{ecbgb1}
G^2=H^3+N^2H^2+8NH+16=H^3+(NH+4)^2
\end{equation}
with
\begin{equation}
\frac{Z}{Y}=\frac{-H}{4} \hspace{2cm} \frac{X}{Y}=\frac{(NH+4) \pm G}{2H}
\end{equation}

The elliptic curve has discriminant
\begin{equation*}
\Delta=4096(N-3)(N^2+3N+9)
\end{equation*}
so is singular when $N=3$, which has the obvious solution $X=Y=Z=1$.

As a simple example, curve for $N=6$ has rank $1$ with point of infinite order $(-3,13)$. This gives
$Z/Y=3/4$ and $X/Y=1/6$ giving $X=2$, $Y=12$ and $Z=9$ as a solution.

The elliptic curves \eqref{ecbgb1} have points of order $3$ at $(0,\pm 4)$. These seem to be the only ones, apart from $N=5$ where there is also
a point of order $2$. Thus the torsion subgroup seems to be isomorphic to $\mathbb{Z}/3\mathbb{Z}$ except when
$N=5$ when it is isomorphic to $\mathbb{Z}/6\mathbb{Z}$.

The presence of the order $3$ torsion points allows us to find a 3-isogenous curve. Using the formulae from the appendix this curve
is
\begin{equation}\label{ecbgb13}
V^2=U^3-27(N\,U-4(N^3-27))^2
\end{equation}

There are a few parametric solutions. We found $N=-(k+1)^2, \, h=4(k+1)$ gives $g=4(k+2)(k^2+k+1)$. This gives $X=1, \, Y=(k+1), \, Z=-(k+1)^2$.

\vspace{1cm}

If we look at problem \eqref{bgb2}, we have
\begin{equation}
a^3+b^3+c^3-Nabc=0
\end{equation}

As a cubic, this fairly unusable. Define $s=a+b+c$, which gives a quadratic
\begin{equation}
(b(N-3)+3s)a^2+(b-s)(b(N-3)+3s)a+3b^2s-3bs^2+s^3=0
\end{equation}

To simplify, define $e=a/s$ and $f=b/s$, giving
\begin{equation}
(f(N-3)+3)e^2+(f-1)(f(N-3)+3)e+3f^2-3f+1=0
\end{equation}

To get a rational solution for $e$, this quadratic must have a discriminant which is a rational square. So there exists $D \in \mathbb{Q}$
such that
\begin{equation}
D^2=(N-3)^2f^4+2N(3-N)f^3+(N^2-6N-18)f^2+(2N+12)f-3
\end{equation}

Using Mordell's method as described in the appendices, this is equivalent to
\begin{equation*}
q^2=p^3-27N(N^3+216)p+54(N^6-540N^3-5832)
\end{equation*}
which does not have a point of order $2$. Points of inflexion give points of order $3$, and this curve has two, one at $p=-9N^2$
and the other at $p=3(N+6)^2$. Making the first point of inflexion $r=0$ gives
\begin{equation*}
q^2=r^3-27(Nr-4(N^3-27))^2
\end{equation*}
which is just the same equation as \eqref{ecbgb13}. The points on this curve all have estimated height one-third of \eqref{ecbgb1}.
There is, thus, no point in considering \eqref{ecbgb1}.

If we take the second point of inflexion and transfer it to the zero point we get
\begin{equation}\label{ecbgb2}
G^2=H^3+((N+6)H+4(N^2+3N+9))^2
\end{equation}
with
\begin{equation}
\frac{a,b}{s}=\frac{\pm G+NH+4(N^2+3N+9)}{2H(N-3)}
\end{equation}

The elliptic curve has discriminant
\begin{equation*}
\Delta=4096(N-3)^3(N^2+3N+9)^3
\end{equation*}
so is singular again when $N=3$. It also seems to have torsion subgroup isomorphic to $\mathbb{Z}/3\mathbb{Z}$ except when $N=5$
when it is isomorphic to $\mathbb{Z}/6\mathbb{Z}$.

In general, therefore none of the torsion points give a finite solution to problem, so we look for points of infinite order.
This problem was taken up by Dave Rusin, then at Northern Illinois University, now at the University of Texas. With the help
of several people, including myself, he eventually produced solutions for \eqref{bgb2} and hence for \eqref{bgb1}. Readers interested
in these solutions should go to the web-site

http://www.math.niu.edu/~rusin/research-math/abcn/

\subsection{Other Cubic Problems}
It is, of course, possible to think of other problems similar to those in the previous two sections. Possibly the most obvious is
\begin{equation}\label{cuberat}
N=\frac{(a+b+c)^3}{a^3+b^3+c^3}
\end{equation}
which will strike a chord with anyone who has taught algebra, in that many students seem convinced that $(a+b+c)^3=a^3+b^3+c^3$ always!

Thus we have
\begin{equation*}
(N-1)a^3-3(b+c)a^2-3(b+c)^2a+(b+c)((b^2+c^2)(N-1)-bc(N+2))=0
\end{equation*}

Equivalent cubics in $b$ and $c$ can be derived, but these are all fairly intractable. We can reduce to quadratics, however, by
defining $s=a+b+c$ and substituting $c=s-a-b$ into the above giving
\begin{equation*}
3N(s-b)a^2-3N(b-s)^2a+s(3b^2N-3bsN+s^2(N-1))=0
\end{equation*}

For simplicity, define $w=a/s$ and $z=b/s$, giving the quadratic in $z$
\begin{equation}\label{quada}
3N(w-1)z^2+3(w-1)^2Nz-N(3w^2-3w+1)+1=0
\end{equation}

For this to have rational solutions, the discriminant must be a rational square, giving the quartic
\begin{equation*}
D^2=9N^2w^4-18N^2w^2+12N(N-1)w-3N(N-4)
\end{equation*}
and defining $D=y/3N$ and $b=x/3N$ gives
\begin{equation} \label{quart1}
y^2=x^4-18N^2x^2+36N^2(N-1)x-27N^3(N-4)
\end{equation}

This quartic is birationally equivalent to an elliptic curve. Using the method described in Mordell \cite{mord},
we find one form of this curve to be
\begin{equation*}
q^2=p^3-34992N^3p-314928N^4(N^2-6N-3)
\end{equation*}

This curve has a point of inflexion at $p=108N^2$, so we transform this point to zero, and simplifying gives
\begin{equation}
v^2=u^3+36N^2u^2+(432N^4-432N^3)u+1296N^6-2592N^5+1296N^4
\end{equation}
which we can write in the simple form
\begin{equation} \label{ec5}
v^2=u^3+36N^2(\,u+6N(N-1)\,)^2
\end{equation}

This curve has discriminant given by
\begin{equation*}
\Delta=2^{12}\,3^9\,N^8\,(1-N)^3\,(N-9)
\end{equation*}
and so is singular for $N=1,9$ corresponding to $(a,b,c)=(1,0,0)$ and $(a,b,c)=(1,1,1)$.

The points $u=0, v= \pm 36N^2(N-1)$ are points of inflexion and so are of order $3$ in the group of rational points.
In general, these are the only finite torsion points so the torsion subgroup
is (usually) isomorphic to $\mathbb{Z} /3\mathbb{Z}$.

The transformations from quartic to elliptic give
\begin{equation}\label{ecrat}
w=\frac{v+36N^2(1-N)}{6Nu}
\end{equation}

It might be thought that we use this formula to give $w$, and substitute into \eqref{quada} to solve for $z$. It is possible, however,
to use the transformation formulae to show that one of the two possible $z$ values is equal to
\begin{equation*}
z=\frac{-v+36N^2(1-N)}{6Nu}
\end{equation*}
and we summarize these results in the form
\begin{equation}\label{trans5}
\frac{a,b}{s}=\frac{36N^2(1-N) \pm v}{6Nu}
\end{equation}

Clearly there is a problem when $u=0$. Thus, to get a solution to the problem, we would need points of infinite order
and so the rank of the elliptic curve would need to be greater than zero.

For example, $N=11$ gives
\begin{equation*}
v^2=u^3+4356u^2+5749920u+1897473600
\end{equation*}
which has only $2$ finite torsion points at $(0,\pm 43560)$, but rank $1$ with point $(-140,34280)$. This gives $b/s=232/231$,
and substituting $b=232$ and $s=231$ into \eqref{quada} gives
\begin{equation*}
33a^2+33a-125032446=0=33(a-1946)(a+1947)
\end{equation*}
The two roots correspond to values for $a$ and $c$, giving $a=1946$, $b=232$, and $c=-1947$.

For $N \in [1,999]$, Table $6.1$ gives the values and heights for the unsolved values. The heights are for the original
curves \eqref{ecrat}, except for $N=772, \, 997$ where the height is for the 3-isogenous curve
\begin{equation*}
g^2=h^3-3(18Nh+324N^2(N-9))^2
\end{equation*}

\begin{table}\label{cratuns}
\begin{center}
\caption{Unsolved values and heights for \eqref{ecrat}}
\begin{tabular}{llllllrrrrr}
$\,$&$\,$&$\,$&$\,$&$\,$&$\,$&$\,$&$\,$&$\,$&$\,$&$\,$\\
N&Ht&$\,$&N&Ht&$\,$&N&Ht&$\,$&N&Ht\\
$\,$&$\,$&$\,$&$\,$&$\,$&$\,$&$\,$&$\,$&$\,$&$\,$&$\,$\\
179&114  &$\,$&302&313  &$\,$&314&178  &$\,$&367&250  \\
382&108  &$\,$&397&251  &$\,$&402&180  &$\,$&421&108  \\
447&218  &$\,$&464&138  &$\,$&478&183  &$\,$&488&213  \\
491&169  &$\,$&509&351  &$\,$&515&141  &$\,$&522&269  \\
527&277  &$\,$&537&236  &$\,$&542&395  &$\,$&549&169  \\
556&138  &$\,$&562&111  &$\,$&566&283  &$\,$&569&317  \\
577&171  &$\,$&587&269  &$\,$&592&128  &$\,$&604&156  \\
659&235  &$\,$&667&112  &$\,$&670&108  &$\,$&682&433  \\
683&167  &$\,$&691&346  &$\,$&719&536  &$\,$&727&220  \\
737&168  &$\,$&739&266  &$\,$&752&98  &$\,$&753&139  \\
758&777  &$\,$&761&346  &$\,$&772&106  &$\,$&773&179  \\
787&944  &$\,$&789&110  &$\,$&797&478  &$\,$&802&604  \\
807&239  &$\,$&808&230  &$\,$&817&242  &$\,$&827&550  \\
839&134  &$\,$&842&368  &$\,$&844&239  &$\,$&849&232  \\
863&1450  &$\,$&866&170  &$\,$&881&187  &$\,$&890&139  \\
899&167  &$\,$&907&1850  &$\,$&911&122  &$\,$&917&189  \\
927&131  &$\,$&929&399  &$\,$&933&248  &$\,$&934&200  \\
937&185  &$\,$&942&266  &$\,$&943&118  &$\,$&954&110  \\
957&109  &$\,$&977&360  &$\,$&982&939  &$\,$&984&101  \\
985&175  &$\,$&986&105  &$\,$&997&304  &$\,$&998&1120
\end{tabular}
\end{center}
\end{table}
\vspace{1cm}

Very similar to the last problem is
\begin{equation}\label{cubic321}
N=\frac{a^3+b^3+c^3}{(a+b+c)(a^2+b^2+c^2)}
\end{equation}
with $a,b,c,N \in \mathbb{Z}$.

We can use a very similar approach to the previous problem, and we find the related elliptic curve to be
\begin{equation}\label{ec321}
v^2=u^3+4(\,3(N-1)u+2(8N^3-24N^2+27N-9)\,)^2
\end{equation}
with
\begin{equation*}
\frac{a,b}{s}=\frac{\pm v + 2Nu + 4(8N^3-24N^2+27N-9)}{6u}
\end{equation*}

The curve has discriminant
\begin{equation*}
\Delta=2^{12} \, 3^3 \,(1-3N) \, (8N^3-24N^2+27N-9)^3
\end{equation*}
so is non-singular for integral $N$.

The torsion subgroup $\mathbb{Z}/3\mathbb{Z}$ with torsion points are $(\,0 , \pm  \, 4(8N^3-24N^2+27N-9)\,)$.

The interesting thing about this problem is that the estimated heights of the rank-$1$ curves are huge compared to any of the other problems in this survey.
For $N \in [1,999]$ the rank-one curves have an average height of $686$ with the largest value being $15640$. I do not know if these are true
values for the heights or are influenced by the unknown value for the size of the Tate-Safarevic group. Given these extreme heights, I have not pursued
this problem with any great enthusiasm.

\vspace{1cm}

The final problem in this section is
\begin{equation} \label{fr3}
N=\frac{a}{b+c}+\frac{b}{c+a}+\frac{c}{a+b}
\end{equation}
with $a,b,c,N \in \mathbb{Z}$.

The same approach as before leads to the elliptic curve
\begin{equation}\label{ecfr3}
v^2=u^3+(4N^2+12N-3)u^2+32(N+3)u
\end{equation}
and the reverse transformation is
\begin{equation*}
\frac{a,b}{s}=\frac{u \pm v - 8(N+3)}{2(N+3)(u-4)} \hspace{1.5cm}
\end{equation*}

The curve \eqref{ecfr3} has discriminant
\begin{equation*}
\Delta=16384(N+3)^2(2N-3)(2N+5)^3
\end{equation*}
so that $\Delta>0$ unless $N=-1,0,1$. Hence the elliptic curve has two components (usually).
The torsion subgroup is $\mathbb{Z}/6\mathbb{Z}$ where the torsion points are $(\,0,0\,)$ of order $2$, $(\,4 , \pm 4(2N+5) \, )$ of order $3$
and $(\,8(N+3) , \pm 8(N+3)(2N+5)\, )$ of order $6$, so all the torsion points are on the infinite component.

Numerical experimentation suggested two unusual behaviours for these elliptic curves, concerning the problem of finding totally
positive solutions. Firstly, we were unable to find a totally positive solution if $N$ was an odd number. Secondly, for some $N$,
the simplest totally positive solution was enormous in size.

I informed Andrew Bremner of these facts and he very quickly proved that they were true. This then led to the joint
paper \cite{brmac}, to which the interested reader is directed for all the specific information on this problem.

\newpage
\section{Equal Sums of Like Powers}
\subsection{Simple Example}
Consider the problem of finding a parametric solution to
\begin{equation}
A^3+B^3=C^3+D^3
\end{equation}

If we use the standard Euler substitution $A=p+q$, $B=r-s$, $C=p-q$ and $D=r+s$, the system reduces to
\begin{equation*}
q(3p^2+q^2)=s(3r^2+s^2)
\end{equation*}

Define $a=p/s$, $b=q/s$ and $c=r/s$ so
\begin{equation*}
3a^2b+b^3-3c^2-1=0
\end{equation*}

For $a \in \mathbb{Q}$, we must have $d \in \mathbb{Q}$ such that
\begin{equation*}
d^2=3b(3c^2+1-b^3)=-3b^4+3(3c^2+1)b
\end{equation*}
and, assuming $c$ is a free rational parameter, we find that $b=1$ gives $d=\pm 3c$. Thus the quartic is birationally equivalent
to an elliptic curve, which we find to be
\begin{equation}\label{ec33}
v^2=u^3-27(3c^2+1)^2
\end{equation}
with
\begin{equation}
b=1+\frac{1}{w} \hspace{1cm} w=\frac{2cv+3(1-c^2)u+18(3c^2+1)}{2(4c^2u-3(c^2+3)(3c^2+1))}
\end{equation}

Equating the denominator of the $w$ equation to zero, gives rational points
\begin{equation*}
u=\frac{3(c^2+3)(3c^2+1)}{4c^2} \hspace{1cm} v= \pm \frac{9(c^4-6c^2-3)(3c^2+1)}{8c^3}
\end{equation*}
on \eqref{ec33}.

Clearly we cannot use these points to get a sensible value of $w$. The elliptic curve is a Mordell curve which are known to have
trivial torsion, so we cannot add any torsion points. We could double the point, but it is better to see if the given points
are double a simpler point. Given the elliptic curve $y^2=x^3+D$ where $D \in \mathbb{Q}$, and a point $(P,Q)$ on the curve, then it
is straightforward algebra that $(P,Q)$ is double a point if the x-coordinate satisfies
\begin{equation*}
x^4-4Px^3-8Dx-4DP=0
\end{equation*}

Substituting the specific values of $D$ and $P$, we find a simpler pair of points on the curve
\begin{equation}
u=3(3c^2+1) \hspace{2cm} v= \pm 9c(3c^2+1)
\end{equation}
and the positive signed $v$ gives
\begin{equation*}
b=\frac{4c^2}{c^2+3} \hspace{2cm} a=\frac{c^4-6c^2-3}{2c(c^2+3)}
\end{equation*}

This allows us to set
\begin{equation*}
p=c^4-6c^2-3 \hspace{1cm} q=8c^3 \hspace{1cm} r=2c^2(c^2+3) \hspace{1cm} s=2c(c^2+3)
\end{equation*}
which give
\begin{equation}
A=c^4+8c^3-6c^2-3 \hspace{2cm} B=2c^4-2c^3+6c^2-6c
\end{equation}
\begin{equation}
C=c^4-8c^3-6c^2-3 \hspace{2cm} D=2c^4+2c^3+6c^2+6c
\end{equation}
where
\begin{equation*}
A^3+B^3=C^3+D^3=9(c^2-1)^3(c^6+33c^4+27c^2+3)
\end{equation*}

\subsection{$A^4+B^4=C^4+D^4$}
This was considered by the great Euler in \cite{euler1}. We have
\begin{equation*}
A^4-D^4=(A-D)(A+D)(A^2+D^2)=C^4-B^4=(C-B)(C+B)(C^2+B^2)
\end{equation*}
so Euler uses the substitutions $A=p+q$, $D=p-q$, $C=r+s$ and $B=r-s$ to give
\begin{equation}
pq(p^2+q^2)=rs(r^2+s^2)
\end{equation}

We follow the description of the method in chapter XXII of Dickson's History \cite{dick1}. Define $p=ax$, $q=by$, $r=kx$
and $s=y$, so that
\begin{equation*}
\frac{y^2}{x^2}=\frac{a^3\,b-k^3}{k-a\,b^3}
\end{equation*}

Setting $k=ab$, gives $y^2/x^2=a^2$. Thus if $x=1$, then $y=\pm a$ and $C = \pm A$ and $D = \pm B$. Euler then sets $k=a\,b\,(1+z)$
giving
\begin{equation}
\frac{y^2}{x^2}=\frac{-a^2(b^2z^3+3b^2z^2+3b^2z+b^2-1)}{z-b^2+1}
\end{equation}
so, for a rational solution, we require $D \in \mathbb{Q}$ such that
\begin{equation*}
D^2=-(b^2z^3+3b^2z^2+3b^2z+b^2-1)(z-b^2+1)
\end{equation*}
which can be written as the quartic
\begin{equation}
D^2=-b^2z^4+b^2(b^2-4)z^3+3b^2(b^2-2)z^2+(3b^4-4b^2+1)z+(b^2-1)^2
\end{equation}

When $z=0$, we have $D= \pm (b^2-1)$, so the quartic is birationally equivalent to the elliptic curve
\begin{equation}\label{eceuler4}
v^2=u^3-3b^4u-b^2(b^8+1)
\end{equation}
with
\begin{equation}\label{euler4z}
z=\frac{(1-b^2)(b^4+10b^2+1-4u)}{b^6+3b^4-3b^2(u+2)+u+2v}
\end{equation}

The elliptic curve \eqref{eceuler4} has discriminant $-432b^4(b^8-1)^2$ so is singular when $b=0, \pm 1$. Numerical
experimentation suggests that the curves have no finite torsion points and also have rank at least $1$. Equating the numerator
of \eqref{euler4z} to zero gives the points
\begin{equation*}
u=\frac{b^4+10b^2+1}{4} \hspace{2cm} v = \pm \frac{b^6-17b^4-17b^2+1}{8}
\end{equation*}
which give $z=0$ or $z$ undefined.

Experiments with Simon's ellrank package suggest that these are not generators. It is straightforward algebraic
manipulation to show that these points are, in fact, double the points with
\begin{equation}
u=b^4+b^2+1 \hspace{2cm} v = \pm (b^6+b^4+b^2+1)
\end{equation}
which, often, are generators.

The negative v-value gives
\begin{equation*}
z=\frac{-3(b^2-1)^3}{4b^6+b^4+10b^2+1}
\end{equation*}
and, hence,
\begin{equation}
y=a(b^2+1)(b^2+4b-1)(b^2-4b-1)
\end{equation}
\begin{equation*}
x=2(4b^6+b^4+10b^2+1)
\end{equation*}

Clearing all obvious common factors gives the parametric solution
\begin{equation*}
A=(b-1)(b^6+9b^5-8b^4-6b^3-23b^2-3b-2)
\end{equation*}
\begin{equation*}
B=(b+1)(2b^6-3b^5+23b^4-6b^3+8b^2+9b-1)
\end{equation*}
\begin{equation*}
C=(b-1)(2b^6+3b^5+23b^4+6b^3+8b^2-9b-1)
\end{equation*}
\begin{equation*}
D=-(b+1)(b^6-9b^5-8b^4+6b^3-23b^2+3b-2)
\end{equation*}
and, for example, $b=5$ gives the solution $2338^4+3351^4=3494^4+1623^4$.

\subsection{An order $4$ multigrade}
Consider the two sets $S_1=\{A+t,B+t,t-C,t-D,t-E\}$ and $S_2=\{t-A,t-B,C+t,D+t,E+t\}$. Then standard algebra
shows
\begin{equation}\label{mg4}
(A+t)^n+(B+t)^n+(t-c)^n+(t-D)^n+(t-E)^n=
\end{equation}
\begin{equation*}
(t-A)^n+(t-B)^n+(C+t)^n+(D+t)^n+(E+t)^n
\end{equation*}
for $n=1,2,3,4$ when
\begin{equation}\label{mg4a}
A+B=C+D+E \hspace{2cm} A^3+B^3=C^3+D^3+E^3
\end{equation}
and we assume $A,B,C,D,E,t \in \mathbb{Q}$. The homogeneity of the equations \eqref{mg4} mean that we can scale any
rational numbers to integers.

\eqref{mg4} is an example of a multigrade equation in the sense of Gloden \cite{glod}. In fact, Gloden discusses this problem, though
not from the perspective of elliptic curves. I consider the problem of solving \eqref{mg4a} and derive a simple parametric solution.

We have $E=A+B-C-D$ and so,
\begin{equation}
(B-C-D)A^2+(B-C-D)^2A-(B-C)(B-D)(C+D)=0
\end{equation}
and for this to give rational roots for $A$ we must have the discriminant being a rational square. Thus the quartic
\begin{equation}
d^2=B^4-2(C+D)^2B^2+4CD(C+D)B+(C^2-D^2)^2
\end{equation}
must have a rational $(B,d)$ solution. Clearly, $B=0$ gives $d=C^2-D^2$ so the curve is birationally equivalent
to an elliptic curve.

Using Mordell's formula gives the elliptic curve
\begin{equation}\label{mg4ec}
v^2=u^3+4(C+D)^2(u+2CD)^2
\end{equation}
with
\begin{equation}\label{mg4tr}
B=\frac{v-4CD(C+D))}{2u}
\end{equation}

The curve \eqref{mg4ec} is in the standard form for a curve with $(0,0)$ as a point of order $3$. There might
be other torsion points for special values of $C,D$, but I have not checked.

Numerical experiments suggested the rank was always at least one with a generator having $u=-4CD$, which gives
$v= \pm 4CD(C-D) $. Slotting this into \eqref{mg4tr} gives $B=C$ or $B=D$. $B=D$ gives $A=0$ or $A=D$
which just lead to the sets $S_1$ and $S_2$ just being permutations of each other. The same occurs if we use the points
obtained by adding the torsion points of order $3$.

We get a non-trivial solution if we double $(-4CD, 4CD(C-D))$. This gives $u=4C^2D^2/(C-D)^2$ which has $v=\pm 4CD(C^4+D^4)/(C-D)^3$.
The positive-sign value gives
\begin{equation*}
B=\frac{C^3-CD^2+D^3}{C(C-D)}
\end{equation*}
which gives $A=(C^3+D^3)/(CD)$ or $A=(C^3-C^2D+D^3)/(D(D-C))$.
Selecting the first gives $E=(C^3-C^2D+D^3)/(D(C-D))$.

For example $C=3, D=2$ gives $A=35/6$, $B=23/3$, and $E=17/2$. Multiplying by $6$ and using $t=70$ gives
\begin{equation*}
19^n+52^n+58^n+105^n+116^n=24^n+35^n+82^n+88^n+121^n \hspace{1cm} (n=1,2,3,4)
\end{equation*}

We can derive a much simpler parametric solution as follows. Since we can assume rational values and then scale up, let $D=1$. The elliptic curve
\eqref{mg4ec} becomes
\begin{equation*}
v^2=u^3+(2(C+1)u+4C(C+1))^2
\end{equation*}

The value $u=-6(C+1)$ would give a rational point when $16C^2-6C-18=W^2$ has a rational solution. It is easy to find that $C=-1$ gives
$W=\pm 2$. Thus the line $W=2+k(C+1)$ should intersect the curve at one further point and it is straightforward algebra to find
this second point is given by
\begin{equation*}
c=\frac{k^2+4k+22}{16-k^2}
\end{equation*}
which gives
\begin{equation*}
u=\frac{12(2k+19)}{k^2-16} \hspace{2cm} v=\frac{8(k^2+19k+16)(2k+19)}{(k^2-16)^2}
\end{equation*}

Using the various formulae above and clearing denominators gives $A=2(k^2-2k-35)$, $B=3(5k-2)$, $C=-3(k^2+4k+22)$, $D=3(k^2-16)$
and $E=(k+2)(2k+19)$.

\newpage

\section{Single quartic made square}
Most of the previous problems reduce to the basic "quartic=rational-square" form. In this section we consider directly problems
where we usually have a function of degree $4$ to start with.

\subsection{Simple Quartics}
In this section, we consider the following three problems: Find integer solutions of

\begin{equation}
1. \hspace{1cm} y^2=x^4+Ny^4
\end{equation}

\begin{equation}
2. \hspace{1cm} y^2=x^4+Nx^2y^2+y^4
\end{equation}

\begin{equation}
3. \hspace{1cm} y^2=x^4+Mx^2y^2+Ny^4
\end{equation}
where $M,N \in \mathbb{Z}$.

The first of these with $N=1$ was proved impossible by Fermat using infinite descent and is, thus, of historical importance. The first
two problems are discussed in Chapter XXII of Dickson's History, which the interested reader is highly recommended to read.
Problem $3$ is mentioned by R.D. Carmichael \cite{carm} as being worthy of investigation.

The first two problems are specializations of the third, so we will just consider it. We find the quartic to be birationally equivalent to
\begin{equation}
v^2=u^3-2Mu^2+(M^2-4N)u
\end{equation}
with
\begin{equation}
\frac{x}{y}=\frac{v}{2u}
\end{equation}

The two parameters $M$ and $N$ generate a vast amount of data, even when we restrict both to lie in $[-999,999]$. We summarize the relevant
information. Of the nearly $4$ million curves, the number of rank $0$ and rank $1$ curves are almost exactly the same. There are only $13$
curves with rank greater than one, which is surprisingly small given the distributions seen in other problems. The largest rank $1$ height
observed is $590/1180$ for $M=-985, N=-193$. Practically all the large height curves happen when both $M$ and $N$ are negative. The
largest height for positive $M,N$ is $101/202$ for $M=695, N=227$.

The first two problems have been completely solved for $N \in [-99999,99999]$.

\subsection{Variants of $(x^2-y^2)(z^2-w^2)=N\,x\,y\,z\,w$}
I first considered this problem after reading a math.stackexchange.com post in early $2016$. The question was specifically about $N=2$
and wanted solutions with $x,y,z,w \in \mathbb{Z}^+$. There are, in fact, no non-trivial rational solutions, but I started thinking about the general
problem of the title.

We have
\begin{equation*}
\left( \frac{x}{y} - \frac{y}{x} \right) \left( \frac{z}{w}-\frac{w}{z} \right)=N
\end{equation*}
or
\begin{equation*}
\left( g - \frac{1}{g} \right) \left( h-\frac{1}{h} \right)=N
\end{equation*}
which is equivalent to the quadratic equation
\begin{equation*}
g^2+\frac{hN}{1-h^2}g-1=0
\end{equation*}

For rational $g$, the discriminant must be a rational square, so there must exist $D \in \mathbb{Q}$ such that
\begin{equation}
D^2=4h^4+(N^2-8)h^2+4
\end{equation}
so we must exclude $N=4$ since the quartic reduces to $D^2=4(h^2+1)^2$.

This quartic has a rational solution $h=0, D=2$, so is birationally equivalent to the elliptic curve
\begin{equation}
V^2=U(U+16)(U+N^2)
\end{equation}
with
\begin{equation}
h=\frac{V}{4(U+N^2)}
\end{equation}

The curves have $3$ points of order $2$ at $(0,0)$, $(-16,0)$ and $(-N^2,0)$ and $4$ points of order $4$ at
\begin{equation*}
(\, 4N \, , \, \pm 4N(N+4) \, ) \hspace{2cm} (\, -4N \, , \, \pm 4N(N-4) \, )
\end{equation*}
all of which give $h=0$, $|h|=1$ or $h$ undefined, and, therefore, no non-trivial solution. We thus need curves
with strictly positive rank.

As an example, $N=47$ is predicted to give rank one. We quickly find $U = -699548/21609$ with $V = \pm 3411924956/3176523$,
which give $h=671/5439$ and finally $g=320/1937$. Thus $(x,y,z,w)=(320,1937,671,5439)$.

For $1 \le N \le 999$, the BSD conjecture predicts $405$ curves have rank $0$, with $505$ having rank $1$, and $88$ rank greater than one while one value
gives a singular curve. The average height for the rank one curves is $9.3$ with the maximum height $107$ at $N=937$. I have results for all $N \in [1,999]$
predicted to have solutions.
As with most problems, some values of $N$, even of moderate height, can cause lots of problems.

\vspace{0.5cm}
It is an obvious variant to consider
\begin{equation}
(x^2-y^2)(z^2+w^2)=Nxyzw
\end{equation}
and, using exactly the same definitions and approach we are led to the elliptic curve
\begin{equation}
V^2=U(U+16)(U+N^2+16)
\end{equation}
with
\begin{equation}
h=\frac{V}{4(U+N^2+16)} \hspace{1cm} g^2-\frac{Nh}{h^2+1}g-1=0
\end{equation}

As an example $N=34$ has a point $(U,V)=(-198068/441, 109821088/9261)$, giving $h=689/168$ and $g=485/61$. Thus,
$(x,y,z,w)=(485,61,689,168)$. For $N \in [1,999]$, the BSD Conjecture predicts $401$ rank $0$ curves, $502$ rank $1$ curves
and $96$ with rank greater than one. The average rank $1$ height is $17.1$ with the maximum height $152$. These are larger
than the previous problem and I still have some unsolved values of $N$ to discover. These are given in the Table $8.1$.
In each case, the 2-isogenous height is double the value given.

\begin{table}[H]
\begin{center}
\caption{Unsolved values of $N$ and heights}
\begin{tabular}{llcrr}
$\,$&$\,$&$\,$&$\,$ &$\,$\\
N&Ht.&$\,$&N&Ht.\\
$\,$&$\,$&$\,$&$\,$ &$\,$ \\
563&128.2&$\,$&587&124.7\\
758&84.81&$\,$&811&102.7\\
823&127.0&$\,$&863&97.96\\
869&72.99&$\,$&911&97.46\\
919&119.6&$\,$&934&92.23\\
983&152.2&$\,$&$\,$&$\,$
\end{tabular}
\end{center}
\end{table}

\vspace{0.5cm}
Clearly, the third variant of this problem is
\begin{equation}
(x^2+y^2)(z^2+w^2)=Nxyzw
\end{equation}
which gives, using the same ideas as before, the requirement that
\begin{equation*}
-4h^4+(N^2-8)h^2-4=\Box
\end{equation*}
but this quartic is qualitatively different than before, as might be expected from the minus signs.

In fact, practically all values of $N$ give a quartic which is not everywhere locally soluble, so has no rational solutions. There
are a few which do give a solution but I have been unable to determine a pattern amongst them.

\vspace{1cm}
A problem in this vein is
\begin{equation}
(x^2+Ay^2)zw=N(z^2+Bw^2)xy
\end{equation}
where $x,y,z,w,N \in \mathbb{Z}$ and $A,B = \pm 1$. The investigation of this problem is still progressing as this survey
is being written.

\newpage

\section{Other Problems}
Practically all the problems in this section are based on the Pythagorean equation
\begin{equation}
x^2+y^2=z^2
\end{equation}
where $x,y,z \in \mathbb{Z}$. It is a standard result that there exists $\alpha, p, q \in \mathbb{Z}$ with $\gcd(p,q)=1$ and
\begin{equation}
x=\alpha\, 2pq \hspace{1.5cm} y=\alpha(p^2-q^2) \hspace{1.5cm} z=\alpha(p^2+q^2)
\end{equation}

The $\alpha$ term can be a pain, but not if we use
\begin{equation}\label{pyrat}
\frac{x}{y}=\frac{2pq}{p^2-q^2}=\frac{2e}{1-e^2}
\end{equation}
where $e=q/p$ is rational.

Note that
\begin{equation*}
\frac{y}{x}=\frac{p^2-q^2}{2pq}=\frac{(p+q)(p-q)}{2pq}=\frac{2rs}{r^2-s^2}=\frac{2f}{1-f^2}
\end{equation*}
if we set $r=p+q$, $s=p-q$ and $f=s/r$.

\subsection{Rational Cuboid problem}
This is possibly the most famous unsolved problem in simple Diophantine Analysis. Anyone interested in this problem
should read the paper by the late John Leech \cite{leech1}, which is full of fascinating details.
Also well worth reading is the paper by Andrew Bremner \cite{brcub} and Ronald van Luijk's thesis.

Find a cuboid with integer sides $L,B,H$ such that the face diagonals and the interior space diagonals are integers.

Thus we
look for $I,J,M,N \in \mathbb{Z}$ such that
\begin{equation}\label{ratcub}
L^2+B^2=I^2 \hspace{2cm} L^2+H^2=J^2 \hspace{0.5cm}
\end{equation}
\begin{equation*}
B^2+H^2=M^2 \hspace{2cm} L^2+B^2+H^2=N^2
\end{equation*}

Since all the relations are homogenous, we can loosen the restrictions on the variables from integer to rational, without loss of
generality. There are $3$ subproblems

\begin{enumerate}
\item The space diagonal can be irrational.
\item One of the face diagonals can be irrational.
\item One of the sides can be irrational.
\end{enumerate}

Consider the first subproblem. Then
\begin{equation*}
\frac{L}{B}=\frac{2e}{1-e^2} \hspace{2cm} \frac{L}{H}=\frac{2f}{1-f^2} \hspace{2cm} \frac{B}{H}=\frac{2g}{1-g^2}
\end{equation*}
for some $e,f,g, \in \mathbb{Q}$.

Since
\begin{equation*}
\frac{L}{B} \frac{H}{L} \frac{B}{H} =1
\end{equation*}
we have the quadratic equation
\begin{equation}
e^2+\frac{2g(f^2-1)}{f(g^2-1)}e-1=0
\end{equation}

Such equations are very common in these types of problems. The first thing that is clear is that the product of the two e-roots is $-1$,
so there is always exactly one positive root.

For the quadratic to have rational roots, the discriminant of the quadratic must be a rational square, so
\begin{equation*}
D^2=g^2f^4+(g^4-4g^2+1)f^2+g^2
\end{equation*}
for some rational $D$.

Define $D=Y/g$ and $f=X/g$ to give the quartic
\begin{equation}
Y^2=X^4+(g^4-4g^2+1)X^2+g^4
\end{equation}
which can be transformed, using Mordell's method \cite{mord}, to the elliptic curve
\begin{equation}
v^2=u(u-(g^2-1)^2)(u-(g^4-6g^2+1))
\end{equation}

We can get an elliptic curve with integer coefficients by setting $g=m/n$, $u=H/n^4$ and $v=G/n^6$ giving
\begin{equation}\label{ecrat1}
G^2=H(H-(m^2-n^2)^2)(H-(m^4-6m^2n^2+n^4))
\end{equation}
with the transformation
\begin{equation}
f=\frac{G}{2mnH}
\end{equation}

We thus find points $(H,G)$ on \eqref{ecrat1}, get $f$, and then $e$ from the quadratic and find $L,B,H$.

The elliptic curve has $3$ simple points of order $2$ when $H=0$, $H=(m^2-n^2)^2$ and $H=m^4-6m^2n^2+n^4$. Numerical
experiments with small $(m,n)$ pairs also suggest points of order $4$. Since doubling a point on such an elliptic curve always give
a point with a square $H$ coordinate, we must have
\begin{equation}
\frac{(H^2-(m^2-n^2)^2(m^4-6m^2n^2+n^4))^2}{4H(H-(m^2-n^2)^2)(H-(m^4-6m^2n^2+n^4))}=T
\end{equation}
where $T$ denotes the $H$ coordinate of a torsion point that is an integer square. $T=0$ does not give a solution, but $T=(m^2-n^2)^2$
gives the following points of order $4$,
\begin{equation*}
((m^2-n^2)(m^2+2mn-n^2), \pm 2mn(m^2-n^2)(m^2+2mn-n^2))
\end{equation*}
and
\begin{equation*}
((m^2-n^2)(m^2-2mn-n^2), \pm 2mn(m^2-n^2)(m^2-2mn-n^2))
\end{equation*}

Thus, the torsion subgroup is isomorphic to $\mathbb{Z}/2\mathbb{Z} \oplus \mathbb{Z}/4\mathbb{Z}$ or
$\mathbb{Z}/2\mathbb{Z} \oplus \mathbb{Z}/8\mathbb{Z}$. I conjecture the latter is impossible.

From any of the finite torsion points, we have $f$ is undefined or $|f|=1$ which is unacceptable. Thus we need
points of infinite order.

Experiments with Denis Simon's {\bf ellrank} package find the first rank $1$ curve when $m+n=7$.
For $m=5$, $n=2$, we have the point $(36,170)$ not a
torsion point. So $g=5/2$, and $f=3/8$ giving $e=18/7$ from the quadratic. Thus
\begin{equation*}
\frac{L}{B}=\frac{-252}{275} \hspace{2cm} \frac{L}{H}=\frac{48}{55} \hspace{2cm} \frac{B}{H}=\frac{-20}{21}
\end{equation*}

Clearing denominators, and remembering we can effectively ignore the signs gives $L=1008$, $B=1100$ and $H=1155$ with face
diagonals $1492$, $1533$ and $1595$.

For each $(m,n)$ pair, we used Pari to find $R$ independent points $P_1, \ldots, P_R$ of infinite order on \eqref{ecrat1}.
If $R \ge 1$, we then generated points
\begin{equation*}
(H,G)= k_1P_1+\ldots + k_rP_R + T
\end{equation*}
where $-L \le k_i \le +L$ for some preset limit $L$, and $T$ is a possible torsion point. From $(H,G)$ we get $f$
and finally $e$ from the quadratic, and hence determine $L,B,H$. The following Table list the cuboid sides with maximum
size $999$. During the relatively short search, a curve of rank at least $3$ was found.

\begin{table}[H]
\begin{center}
\caption{Cuboids with irrational space diagonal}
\begin{tabular}{rrr}
$\,$&$\,$&$\,$\\
L&B&H\\
$\,$&$\,$&$\,$\\
44&117&240\\
240&252&275\\
140&483&693\\
85&132&720\\
150&231&792
\end{tabular}
\end{center}
\end{table}

\bigskip

Now, suppose one of the face diagonals is irrational. Thus
\begin{equation}
L^2+B^2=I^2 \hspace{1.5cm} L^2+H^2=J^2 \hspace{1.5cm}L^2+B^2+H^2=N^2
\end{equation}
where $A$ is also an integer.

Thus
\begin{equation*}
\frac{L}{B}=\frac{2e}{1-e^2} \hspace{1cm} \frac{L}{H}=\frac{2f}{1-f^2} \hspace{1cm} \frac{I}{H}=\frac{2g}{1-g^2}
\end{equation*}
and we also have
\begin{equation*}
\frac{L}{I}=\frac{2e}{1+e^2}
\end{equation*}

Thus
\begin{equation}
\frac{2f}{1-f^2}=\frac{2g}{1-g^2} \frac{2e}{1+e^2}
\end{equation}
giving the quadratic in $e$
\begin{equation}
e^2+\frac{2g(f^2-1)}{f(1-g^2)}e+1=0
\end{equation}

If $e$ is to be rational, the discriminant must be a rational square so rational $D$ must exist with
\begin{equation}
D^2=g^2f^4-(g^4+1)f^2+g^2
\end{equation}
and, defining $D=Y/g$ and $f=X/g$, we derive the quartic
\begin{equation}
Y^2=X^4-(1+g^4)X^2+g^4
\end{equation}

Using Mordell's method, this quartic can be transformed to the equivalent elliptic curve
\begin{equation}\label{ecrat2}
G^2=H(H+(m^2+n^2)^2)(H+(m^2-n^2)^2)
\end{equation}
where $g=m/n$ with $\gcd(m,n)=1$ and
\begin{equation*}
f=\frac{G}{2mnH}
\end{equation*}

The elliptic curve has $3$ finite points of order $2$ at $H=0$, $H=-(m^2+n^2)^2$ and $H=-(m^2-n^2)^2$. There are also $4$ points
of order $4$ at
\begin{equation*}
(m^4-n^4, \pm 2m^2(m^4-n^4)) \hspace{2cm} (n^4-m^4, \pm 2n^2(m^4-n^4))
\end{equation*}

As before, we generated the elliptic curve, found (if any) points of infinite order and found $e,f$ and hence solutions. The following Table gives
the smallest found so far, in a fairly short test.

\begin{table}[H]
\begin{center}
\caption{Cuboids with one irrational face diagonal}
\begin{tabular}{rrr}
$\,$&$\,$&$\,$\\
L&B&H\\
$\,$&$\,$&$\,$\\
104&153&672\\
117&520&756\\
448&495&840\\
264&495&952\\
264&448&975
\end{tabular}
\end{center}
\end{table}

\bigskip

Finally, consider the case of one side being a quadratic irrational. Thus we look for rational solutions of
\begin{equation}
L^2+B^2=I^2 \hspace{2cm} L^2+H=J^2
\end{equation}
\begin{equation*}
B^2+H=M^2 \hspace{2cm} L^2+B^2+H=N^2
\end{equation*}

We can assume, without loss of generality, that $L=2e$, $B=1-e^2$ and $I=1+e^2$ for some rational $e$. Also
\begin{equation*}
H=J^2-L^2=M^2-B^2=N^2-I^2
\end{equation*}
so that $B^2+J^2=L^2+M^2$, which we can express as a vector relation
\begin{equation}
\left| \left| \left( \begin{array}{c}B\\J \end{array} \right) \right| \right| ^2 =
\left| \left| \left( \begin{array}{c}L\\M \end{array} \right) \right| \right| ^2
\end{equation}

Thus one vector is just a rotation about the origin of the other vector, so that
\begin{equation}
\left( \begin{array}{c}B\\J \end{array} \right) = \left( \begin{array}{rr} \cos \theta&\sin \theta\\-\sin \theta&\cos \theta \end{array} \right)
 \left( \begin{array}{c}L\\M \end{array} \right)
\end{equation}

Thus
\begin{equation*}
B=L \, \cos \theta  +M \, \sin \theta \hspace{2cm} J = -L \, \sin \theta + M \, \cos \theta
\end{equation*}
which we can rearrange to give
\begin{equation}
J=\frac{B \, \cos \theta - L}{\sin \theta}
\end{equation}

It is clear that $\cos \theta$ and $\sin \theta$ must be rational so
\begin{equation*}
\cos \theta =\frac{1-f^2}{1+f^2} \hspace{2cm} \sin \theta = \frac{2f}{1+f^2}
\end{equation*}
for $f \in \mathbb{Q}$.

Thus $H$ is given by the form
\begin{equation*}
\frac{(f^2-1)\,(e^2(f+1)+2e(1-f)-f-1)\,(e^2(f-1)-2e(f+1)-f+1)}{4f^2}
\end{equation*}

Thus, to force $L^2+B^2+H$ to be a rational square we must have
\begin{equation*}
D^2=(f^2+1)^2e^4+4e^3(f+1)(1-f)(f^2+1)e^3+
\end{equation*}
\begin{equation*}
2(f^2+1)^2e^2+4(f^2-1)(f^2+1)e+(f^2+1)^2
\end{equation*}
with $D,e,f \in \mathbb{Q}$.

Clearly, $e=0$ gives $D=f^2+1$, so the quartic is birationally equivalent to an elliptic curve. We find
\begin{equation}\label{eccub3}
v^2=u^3+2(f^2+1)^2u^2+4f^2(f^2+1)^2u=u(u+2(f^2+1))(u+2f^2(f^2+1))
\end{equation}
with
\begin{equation}\label{trcub3}
e=\frac{f^2u-u+v}{(f^2+1)(4f^2+u)}
\end{equation}

The elliptic curve seems to only have the $3$ finite points of order $2$, $(0,0)$, $(-2(f^2+1),0)$ and $(-2f^2(f^2+1),0)$ which
give $e=0, \pm 1$ which do not give real-world cuboids. Making the denominator of \eqref{trcub3} zero, gives $u=-4f^2$ which has
$v = \pm 4f^2(f^2-1)$. Adding to $(0,0)$ gives the points $(-(f^2+1)^2, \pm (f^2-1)(f^2+1)^2)$.

The negative value of $v$ gives
\begin{equation}
e=\frac{2(f^2+1)}{f^2-1}
\end{equation}
and slotting into the various formulae, we get the parametric formulae (ignoring the sign in $B$)
\begin{equation*}
L=\frac{4(f^2+1)}{f^2-1} \hspace{2cm} B= \frac{(f^2+3)(3f^2+1)}{(f^2-1)^2}
\end{equation*}
\begin{equation*}
H=\frac{(f^4+8f^3-2f^2+8f+1)(f^4-8f^3-2f^2-8f+1)}{4f^2(f^2-1)^2}
\end{equation*}

If we wish a real-world cuboid we must have $|f|>1$ for positive $L$. $B$ is positive for all $f \ne \pm 1$. $H$, however, is
more restrictive. The denominator is always positive when $|f| \ne 1$, but the numerator has $4$ real roots at
$|f| \approx 8.35$ and $|f| \approx 0.1197$.

For example, $f=9$ gives $L=41/10, \, B=1281/400, \, H=23839/8100$.

This parametrization would give a perfect cuboid if the numerator of $H$ were a perfect square. Thus, we would need
$w=f^2$ to lie on the quartic
\begin{equation*}
t^2=w^4-68w^3-122w^2-68w+1
\end{equation*}
which is birationally equivalent to
\begin{equation*}
v^2=u^3-8u^2+20u \hspace{2cm} w=\frac{17u+2v}{u-80}
\end{equation*}

This elliptic curve just has the one finite torsion point at $(0,0)$ and has rank $1$ with generator $(2,4)$. Computing
multiples of the generator and the torsion point does not discover any square value of $w$. There might be one at larger
multiples, but I doubt it.

\vspace{2cm}

At the end of \cite{leech1}, Leech asks for solutions of the following related problems, where we assume $(X,Y,Z,A,B,C) \in \mathbb{Z}$
satisfy
\begin{equation*}\label{rceq}
X^2+Y^2=A^2 \hspace{1.5cm} X^2+Z^2=B^2 \hspace{1.5cm} Y^2+Z^2=C^2
\end{equation*}

{\bf Problem 1:} Find $T \in \mathbb{Z}$ such that
\begin{equation}\label{prob1}
T^2-X^2=E^2 \hspace{1.5cm} T^2-Y^2=F^2 \hspace{1.5cm} T^2-Z^2=G^2
\end{equation}
with $E, F, G$ all integers.

{\bf Problem 2:} Find $T \in \mathbb{Z}$ such that
\begin{equation}\label{prob2}
T^2-X^2-Y^2=E^2 \hspace{1.5cm} T^2-X^2-Z^2=F^2 \hspace{1.5cm} T^2-Y^2-Z^2=G^2
\end{equation}
with $E, F, G$ all integers.
Note that this can be written
\begin{equation*}
T^2-A^2=E^2 \hspace{1.5cm} T^2-B^2=F^2 \hspace{1.5cm} T^2-C^2=G^2
\end{equation*}

A solution to the perfect rational cuboid would give a solution to both problems. Conversely, if either problem has no solution,
then there cannot be a perfect rational cuboid.

In \cite{ajm5}, I show that the second problem always has a solution, but am unable to find a single solution to Problem $1$. Based on this,
I conjecture that a perfect rational cuboid does not exist, and have stopped looking for one.

\subsection{Tiling the unit square with triangles}
This problem is concerned with subdividing the unit square into $N$ triangles with rational sides. The basic reference is by Richard Guy
in \cite{guy2}, but this is quite hard to find. Follow-up papers by Bremner and Guy \cite{bg2} \cite{bg3} are easier to find.

For $N=2$, we only have one possibility
with two equal triangles with the diagonal of the square as one side. But $\sqrt{2}$ is probably the first number that was proven to be irrational.
Guy proves in \cite{guy2} that $N=3$ is impossible, so $N=4$ is the first possibility for a solution. There are $4$ different types of configuration
of triangles which Guy christens the $\Delta$, $\nu$, $\kappa$ and $\chi$ configurations, for obvious visual reasons.

For the $\Delta$ configuration, consider the unit square
\begin{center}
\begin{picture}(8,8)
\put(2,0){\line(1,0){8}}
\put(2,0){\line(0,1){8}}
\put(2,8){\line(1,0){8}}
\put(10,0){\line(0,1){8}}
\put(2,0){\line(2,1){8}}
\put(2,0){\line(1,2){4}}
\put(6,8){\line(1,-1){4}}
\put(1.5,0){O}
\put(10.25,4){P}
\put(6,8.25){Q}
\put(4,8.25){X}
\put(10.25,2){Y}
\end{picture}
\end{center}
with corners at $(0,0)$, $(1,0)$, $(1,1)$ and $(0,1)$. Let $P$ be $(1,Y)$ and $Q$ be $(X,1)$ with $X,Y \in \mathbb{Q}$.

Thus we need
\begin{equation}
1+Y^2=\Box \hspace{1cm} X^2+1=\Box \hspace{1cm} (1-X)^2+(1-Y)^2=\Box
\end{equation}
so that there must exist $e,f,g \in \mathbb{Q}$ with
\begin{equation}
Y=\frac{2e}{1-e^2} \hspace{1cm} X=\frac{2f}{1-f^2} \hspace{1cm} \frac{1-X}{1-Y}=\frac{2g}{1-g^2}
\end{equation}

This gives the quadratic in $e$
\begin{equation*}
e^2+\frac{2(f^2-1)(g^2-1)}{f^2(g^2+2g-1)+4fg-g^2-2g+1}e-1=0
\end{equation*}
and, for $e$ to be rational, the discriminant must be a rational square giving the quartic in $f$
\begin{equation*}
D^2=8(g^4+2g^3-2g+1)f^4+32g(g^2+2g-1)f^3-16(g^4+2g^3-4g^2-2g+1)f^2-
\end{equation*}
\begin{equation*}
32g(g^2+2g-1)f+8(g^4+2g^3-2g+1)
\end{equation*}

For $f=1$ we have $D=\pm 8g$, so define $z=1/(f-1)$ which gives the equivalent quartic
\begin{equation*}
w^2=64g^2z^4+64g(g^2+4g-1)z^3+32(g^4+5g^3+8g^2-5g+1)z^2+
\end{equation*}
\begin{equation*}
32(g^4+3g^3+2g^2-3g+1)z+8(g^4+2g^3-2g+1)
\end{equation*}

Finally, define $w=y/(8g)$ and $x=z/(8g)$ to give
\begin{equation*}
y^2=x^4+8(g^2+4g-1)x^3+32(g^4+5g^3+8g^2-5g+1)x^2+
\end{equation*}
\begin{equation*}
256g(g^4+3g^3+2g^2-3g+1)x+512g^2(g^4+2g^3-2g+1)
\end{equation*}

The quartic can be transformed, as described in Mordell \cite{mord}, to the equivalent elliptic curve
\begin{equation}\label{ecdelta1}
v^2=u(u-2g^2(g+1)^2)(u-2(g^2+1)(g^2+2g-1))
\end{equation}
with the reverse transformation
\begin{equation}\label{trdelta}
z=\frac{v-(g^2+4g-1)u+2g(g+1)^3(g^2+2g-1)}{4g(u-(g+1)^2(g^2+2g-1))}
\end{equation}
and $f=1+1/z$.

Suppose $g=m/n$ with $m,n \in \mathbb{Z}$ and $\gcd(m,n)=1$. Defining $v=G/n^6$ and $u=H/n^4$ gives the equivalent elliptic curve
with integer coefficients
\begin{equation}\label{ecdelta2}
G^2=H(H-2m^2(m+n)^2)(H-2(m^2+n^2)(m^2+2mn-n^2))
\end{equation}
which has $3$ points of order $2$, when $H=0, 2m^2(m+n)^2, 2(m^2+n^2)(m^2+2mn-n^2)$. These seem to be the only torsion points
for most $(m,n)$ pairs.

The denominator of \eqref{trdelta} is zero when $u=(g+1)^2(g^2+2g-1)$ which is equivalent to $H=(m+n)^2(m^2+2mn-n^2)$, which gives
$G=\pm (m+n)^2(m-n)^2(m^2+2mn-n^2)$. Adding to $(0,0)$ gives the slightly simpler points $(4m^2(m^2+n^2), \pm 4m^2(m-n)^2(m^2+n^2))$.
These points are (normally) of infinite order, though they might not be generators.

We could use these points to generate parametric solutions, but this would only be sensible if the ranks of \eqref{ecdelta2} were
always $1$, and the points of infinite order always gave a generator. But simple numerical experiments show that the ranks can be larger.
We have found examples of rank $4$ in experiments. Thus, a computational approach was taken.

For each $(m,n)$ pair, we found $R$ independent points $P_1, \ldots, P_R$ and generated points
\begin{equation*}
(H,G)= k_1P_1+\ldots + k_rP_R + T
\end{equation*}
where $-L \le k_i \le +L$ for some preset limit $L$, and $T$ is a possible torsion point. From $(H,G)$ get $(u,v)$, then $z$, then $f$
and $e$ from the quadratic, and finally $X$ and $Y$. If we wish a true splitting of the unit square we have to have $0 < X < 1$ and
$0 < Y < 1$.

I wrote a pari program for this computation, which finds solutions quickly. For $L=2$, and using Denis Simon's {\bf ellrank} software
to find the exact rank of the elliptic curves we get the following "small" solutions

\begin{table}[H]
\begin{center}
\caption{$\Delta$ configurations}
\begin{tabular}{rr}
$\,$&$\,$\\
X&Y\\
$\,$&$\,$\\
7/24&28/45\\
87/416&451/780\\
3737/4416&2387/2484\\
1652/3285&22895/24528
\end{tabular}
\end{center}
\end{table}

The following conjecture seems reasonable:

{\bf Conjecture:} $\forall \epsilon$ with $0 < \epsilon < 1/\sqrt{2}$ there exists a $\Delta$ configuration such that $|PQ| < \epsilon$ and one such that
$|PQ| > \sqrt{2} - \epsilon$.

\bigskip

For the $\nu$ configuration, consider the unit square
\begin{center}
\begin{picture}(8,8)
\put(2,0){\line(1,0){8}}
\put(2,0){\line(0,1){8}}
\put(2,8){\line(1,0){8}}
\put(10,0){\line(0,1){8}}
\put(2,0){\line(1,4){2}}
\put(4,8){\line(1,-2){4}}
\put(8,0){\line(1,4){2}}
\put(1.5,0){O}
\put(8.25,0.15){P}
\put(4,8.15){Q}
\put(4,0.15){X}
\put(3,7.5){Y}
\end{picture}
\end{center}
with corners at $(0,0)$, $(1,0)$, $(1,1)$ and $(0,1)$. Let $P$ be $(X,0)$ and $Q$ be $(Y,1)$ with $X,Y \in \mathbb{Q}$.

Thus
\begin{equation*}
1+Y^2=\Box \hspace{1.5cm} 1+(1-X)^2=\Box \hspace{1.5cm} (X-Y)^2+1=\Box
\end{equation*}
so there exists $e,f,g \in \mathbb{Q}$ such that
\begin{equation*}
Y=\frac{2e}{1-e^2} \hspace{1.5cm} 1-X=\frac{2f}{1-f^2} \hspace{1.5cm} X-Y=\frac{2g}{1-g^2}
\end{equation*}

Since $1-(X-Y)=(1-X)+Y$ we have the quadratic equation in $e$
\begin{equation}
e^2+2\frac{(f^2-1)(g^2-1)}{f^2(g^2+2g-1)+2f(g^2-1)-(g^2+2g-1)}e-1
\end{equation}

For this to have rational solutions, the discriminant must be a rational square, so there must exist $D \in \mathbb{Q}$ such that
\begin{equation*}
D^2=2(g^4+2g^3-2g+1)f^4+4(g^2-1)(g^2+2g-1)f^3-
\end{equation*}
\begin{equation*}
8g(g^2+g-1)f^2+4(1-g^2)(g^2+2g-1)f+2(g^4+2g^3-2g+1)
\end{equation*}

Substituting $f=1$, gives $D^2=4(g+1)^2(g-1)^2$, so we hake the substitution $h=1/(f-1)$, which gives
\begin{equation*}
E^2=(2(g^2-1)^2)^2h^4+16(g^4+g^3-2g^2-g+1)h^3+
\end{equation*}
\begin{equation*}
8(3g^4+5g^3-4g^2-5g+3)h^2+4(3g^4+6g^3-2g^2-6g+3)h+2(g^4+2g^3-2g+1)
\end{equation*}

Defining $E=y/(2(g^2-1))$ and $h=x/(2(g^2-1))$ gives
\begin{equation*}
y^2=x^4+8(g^2+g-1)x^3+8(3g^4+5g^3-4g^2-5g+3)x^2+
\end{equation*}
\begin{equation}
8(g^2-1)(3g^4+6g^3-2g^2-6g+3)x+8(g^2-1)^2(g^4+2g^3-2g+1)
\end{equation}

This quartic can be transformed, by Mordell's method, to the equivalent elliptic curve
\begin{equation}\label{ecnuconf}
v^2=u^3+(3g^4+4g^3-2g^2-4g+3)u^2+(g^2-1)^4u
\end{equation}
with
\begin{equation}
h=\frac{v-2(g^2+g-1)u-(g^2-1)^3}{2(g^2-1)(u+(g^2-1)^2)}
\end{equation}

The elliptic curve has a point of order $2$ at $(0,0)$, and $2$ points of order $4$ at $(-(g^2-1)^2, \pm (g^2-1)^2(g^2+2g-1) \, )$,
and, in general, these are the only ones.

Using the same computational procedure as for the $\Delta$ configuration, the following table of small results was easily found

\begin{table}[H]
\begin{center}
\caption{$\nu$ configurations}
\begin{tabular}{rr}
$\,$&$\,$\\
X&Y\\
$\,$&$\,$\\
7/24&7/12\\
12/35&79/420\\
33/56&73/168\\
71/2520&391/420
\end{tabular}
\end{center}
\end{table}

\newpage

\begin{center}
\begin{picture}(8,8)
\put(0,0){\line(1,0){11}}
\put(0,0){\line(0,1){11}}
\put(0,11){\line(1,0){11}}
\put(11,0){\line(0,1){11}}
\put(0,0){\line(1,5){2.2}}
\put(1,5){\line(5,3){10}}
\put(1,5){\line(2,-1){10}}
\put(0.5,4.75){P}
\put(2,11.3){Q}
\end{picture}
\end{center}

For the $\kappa$ configuration, consider the unit square
with corners at $(0,0)$, $(1,0)$, $(1,1)$ and $(0,1)$. Let $P$ be $(X,Y)$ with $X,Y \in \mathbb{Q}$.

The first point to note is that $Q$ is $(Y/X,1)$. If $O$ is the origin, then $OP^2=X^2+Y^2$ and
$OQ^2=(X^2+Y^2)/X^2$. Thus, if $OP$ is rational then so is $OQ$ and also $PQ$. Thus there are only $3$ independent conditions
\begin{equation*}
X^2+Y^2=\Box \hspace{1.5cm} (1-X)^2+Y^2=\Box \hspace{1.5cm} (1-X)^2+(1-Y)^2=\Box
\end{equation*}
which gives
\begin{equation*}
\frac{X}{Y}=\frac{2e}{1-e^2} \hspace{1.5cm} \frac{1-X}{Y}=\frac{2f}{1-f^2} \hspace{1.5cm} \frac{1-X}{1-Y}=\frac{2g}{1-g^2}
\end{equation*}

We have
\begin{equation*}
\frac{1}{Y}=\frac{2e}{1-e^2}+\frac{2f}{1-f^2}=\frac{g(1-f^2)}{f(1-g^2)}+1
\end{equation*}
which gives the quadratic equation
\begin{equation}
e^2+\frac{2g(f^2-1)}{(f^2g+f(g^2+2g-1)-g)}e-1=0
\end{equation}

For this to give rational roots, the discriminant must be a rational square, so $\exists D \in \mathbb{Q}$ such that
\begin{equation*}
D^2=8g^2(f^4+1)+8g(g^2+2g-1)(f^3-f)+4(g^4+4g^3-2g^2-4g+1)f^2
\end{equation*}

If $f=1$, then $D=2(g^2+2g-1)$, so make the transformation $z=1/(f-1)$, which gives the quartic
\begin{equation*}
w^2=4(g^2+2g-1)z^4+8(g^4+6g^3+6g^2-6g+1)z^3+
\end{equation*}
\begin{equation*}
4(g^4+10g^3+22g^2-10g+1)z^2+8g(g^2+6g-1)z+8g^2
\end{equation*}

Define $w=y/(2(g^2+2g-1))$ and $z=x/(2(g^2+2g-1))$ to give
\begin{equation}
y^2=x^4+4(g^2+4g-1)x^3+4(g^4+10g^3+22g^2-10g+1)x^2+
\end{equation}
\begin{equation*}
16g(g^2+2g-1)(g^2+6g-1)x+32g^2(g^2+2g-1)^2
\end{equation*}

Mordell's method transforms this to the equivalent elliptic curve
\begin{equation}\label{eckappa}
v^2=u^3+(g^4+4g^3+10g^2-4g+1)u^2+4g^2(g^2+2g-1)^2u
\end{equation}
with
\begin{equation}
z=\frac{v-(g^2+4g-1)u-4g^2(g^2+2g-1)}{2(g^2+2g-1)(u+4g^2)}
\end{equation}

The curve \eqref{eckappa} has a point of order $2$ at $(0,0)$. Three points of order $2$ will occur only when
$g^4+8g^3+18g^2-8g+1$ is a rational square. This quartic is birationally equivalent to $s^2=t^3+6t^2+t$, which has
rank $0$, and finite points $(0,0)$ and $(-1,\pm 2)$, which give $g=0$ or $g$ undefined. So there are never $3$ points
of order $2$.

There are points of order $4$ when $u=-2g(g^2=2g-1)$ which gives $v=\pm 2g(g^2+2g-1)(g^2+1)$. Investigations suggest
that these are the only torsion points.

The denominator of the $z$ transformation suggests $u=-4g^2$ should give a rational point, and we find $v=\pm 8g^3$.
If we add this point to $(0,0)$ we get $u=-(g^2+2g-1)^2$ which has $v=\pm 2g(g^2+2g-1)^2$. Taking the positive
value gives the following
\begin{equation*}
f=\frac{-(g^4-10g^2+1)}{g^4+8g^3+6g^2-8g+1} \hspace{1.5cm} e=\frac{g^4+4g^3-10g^2-4g+1}{g^4+4g^3+6g^2-4g+1}
\end{equation*}
from which we could derive parametric formulae for $X$ and $Y$.

Instead, we use the same computational procedure as for the $\Delta$ configuration, and the following table of small results
was easily found.

\begin{table}[H]
\begin{center}
\caption{$\kappa$ configurations}
\begin{tabular}{rr}
$\,$&$\,$\\
X&Y\\
$\,$&$\,$\\
7/13&45/52\\
11/39&16/65\\
85/148&42/185\\
76/175&297/700
\end{tabular}
\end{center}
\end{table}

\vspace{1cm}

\begin{center}
\begin{picture}(8,8)
\put(0,0){\line(1,0){12}}
\put(0,0){\line(0,1){8}}
\put(0,8){\line(1,0){12}}
\put(12,0){\line(0,1){8}}
\put(0,0){\line(2,1){8}}
\put(8,4){\line(1,-1){4}}
\put(0,8){\line(2,-1){8}}
\put(8,4){\line(1,1){4}}
\end{picture}
\end{center}

Finally, for the $\chi$ configuration, consider the rectangle shown.

We have
\begin{equation*}
a^2+c^2=\Box \hspace{1cm} a^2+d^2=\Box \hspace{1cm} b^2+c^2=\Box \hspace{1cm} b^2+d^2=\Box
\end{equation*}
with a solution to the $\chi$ configuration when $a,b,c,d >0$ and $a+b=c+d$.

From
\begin{equation*}
\frac{c}{a}=\frac{2e}{1-e^2} \hspace{1cm} \frac{d}{a}=\frac{2f}{1-f^2} \hspace{1cm} \frac{c}{b}=\frac{2g}{1-g^2} \hspace{1cm}
\frac{d}{b}=\frac{2h}{1-h^2} \hspace{1cm}
\end{equation*}
and
\begin{equation*}
\frac{c}{a} \, \, \, \frac{d}{b}=\frac{c}{b} \, \, \, \frac{d}{a}
\end{equation*}
we have the quadratic equation
\begin{equation}
h^2+\frac{e(1-f^2)(1-g^2)}{fg(1-e^2)}h -1 =0
\end{equation}
which must have a discriminant being a rational square to give rational $h$.

Thus, there exists $D \in \mathbb{Q}$, such that
\begin{equation*}
D^2=e^2(f^2-1)^2g^4+(4e^4f^2-2e^2(f^2+1)^2+4f^2)g^2+e^2(f^2-1)^2
\end{equation*}

Defining $D=y/(e(f^2-1))$ and $g=x/(e(f^2-1))$ gives the quartic
\begin{equation}
y^2=x^4+2(2e^4f^2-e^2(f^2+1)^2+2f^2)x^2+e^4(f^2-1)^4
\end{equation}

Using Mordell's method, this can be transformed to the equivalent elliptic curve
\begin{equation}
v^2=w(w+4e^2(f^2-1)^2)(w+4f^2(e^2-1)^2)
\end{equation}
with
\begin{equation}
g=\frac{v}{2e(f^2-1)(w+4f^2(e^2-1)^2)}
\end{equation}

There are, thus, two free parameters $e$ and $f$. For each possible pair, we try to find curves with rank greater than zero,
and possible generators. Then generate points on the curve, and values of $g$ and solve the quadratic for $h$. From these
find $(a,b,c,d)$ and test if they are all positive and if $a+b=c+d$. No solution has been found, and, in fact, Bremner and Guy
conjecture that a solution does not exist. I must admit that I agree with them.

An interesting question is how small can $|(a+b)/(c+d)-1|$ be?

It is easy to adapt the code for this.

\subsection{$3 \times 3$ Magic Square of Squares}
Consider the $3 \times 3$ magic square
\begin{equation}
\begin{array}{ccccc}
a+b&$\hspace{1cm}$&a-b-c&\hspace{1cm}&a+c\\\,&\,&\,\\a-b+c&\hspace{1cm}&a&\hspace{1cm}&a+b-c\\\,&\,&\,\\a-c&\hspace{1cm}&a+b+c&\hspace{1cm}&a-b
\end{array}
\end{equation}
where we assume $a,b,c \in \mathbb{Z}$.

Can we find such a square where the nine entries are distinct squares?

So far, only a single square with $7$ square entries has been found, by Andrew Bremner \cite{bremsq1} \cite{bremsq2}. It is
\begin{equation*}
\begin{array}{ccccc}
373^2&\,&289^2&\,&565^2\\\,&\,&\,\\360721&\,&425^2&\,&23^2\\\,&\,&\,\\205^2&\,&527^2&\,&222121
\end{array}
\end{equation*}

On the principle that one solution usually implies more than one, I have spent a lot of computer time looking for more solutions.
They are all based on the elliptic curve approach described by Bremner, starting from making $6$ entries square.

To illustrate the approach, suppose we take the first and third rows to all be squares, so
\begin{equation*}
\begin{array}{lcr}
a+b=p^2 \hspace{1cm}&a+c=r^2&\hspace{1cm}a+b+c=u^2\\a-b=q^2\hspace{1cm}&a-c=s^2&\hspace{1cm}a-b-c=v^2
\end{array}
\end{equation*}
and we can, without loss of generality, consider all these variables to be rational.

We have
\begin{equation*}
2a=p^2+q^2=r^2+s^2=u^2+v^2
\end{equation*}
and
\begin{equation*}
2(b+c)=u^2-v^2=p^2-q^2+r^2-s^2
\end{equation*}

If we think in terms of vectors, the first relation says that the lengths of
\begin{equation*}
\left( \begin{array}{c}p\\q\end{array} \right) \hspace{2cm} \left( \begin{array}{c}r\\s\end{array} \right)
\hspace{2cm} \left( \begin{array}{c}u\\v\end{array} \right)
\end{equation*}
are the same.

Thus each vector must just be a rotation of another, so
\begin{equation*}
\left( \begin{array}{c}r\\s\end{array} \right) = \left( \begin{array}{rr}\cos \theta&\sin \theta\\-\sin \theta&\cos \theta \end{array} \right)
\left( \begin{array}{c}p\\q\end{array} \right)
\end{equation*}
and
\begin{equation*}
\left( \begin{array}{c}u\\v\end{array} \right) = \left( \begin{array}{rr}\cos \mu&\sin \mu\\-\sin \mu&\cos \mu \end{array} \right)
\left( \begin{array}{c}p\\q\end{array} \right)
\end{equation*}
for some angles $\theta$ and $\mu$.

These relations imply that the sines and cosines are rational, so
\begin{equation*}
\cos \theta=\frac{1-f^2}{1+f^2} \hspace{2cm} \sin \theta=\frac{2f}{1+f^2},
\end{equation*}
and
\begin{equation*}
\cos \mu=\frac{1-g^2}{1+g^2} \hspace{2cm} \sin \mu=\frac{2g}{1+g^2},
\end{equation*}
where $f,g \in \mathbb{Q}$. These give $r,s,u,v$ in terms of $p,q,f,g$. Substituting into the second relation gives a quadratic
\begin{equation*}
p^2+R(f,g)pq-q^2=0
\end{equation*}
where $R(f,g)$ is a fairly complicated rational function in $f,g$.

If this quadratic is to give rational values for $p,q$, the discriminant $R^2+4$ must be a rational square. Computing this, and clearing
as many square terms as possible, we eventually arrive at the quartic in $f$, assuming $g$ has been selected,
\begin{equation*}
d^2=(g^4+34g^2+1)f^4+32g(1-g^2)f^3+2(g^4-30g^2+1)f^2+
\end{equation*}
\begin{equation*}
32g(g^2-1)f+g^4+34g^2+1
\end{equation*}

This quartic has a rational solution when $f=1$ and $d=\pm 2(g^2+1)$. Thus the quartic is birationally equivalent to an elliptic curve.
After a reasonable amount of computer algebra, we find the relevant curve to be
\begin{equation}
w^2=z^3+2(m^4+18m^2n^2+n^4)z^2+Tz
\end{equation}
where $g=m/n$ with $m,n \in \mathbb{Z}$ and $\gcd(m,n)=1$ and $T=(m^2+4mn+n^2)^2(m^2-4mn+n^2)^2$.

The inverse transformation is
\begin{equation*}
f=\frac{-(m^2+n^2)w-(m^4+8m^3n+2m^2n^2-8mn^3+n^4)z-T }
{-(m^2+n^2)w+(m^4-8m^3n+2m^2n^2+8mn^3+n^4)z+T}
\end{equation*}

The elliptic curve has $3$ finite points of order $2$, at $z=0$, $z=-(m^2+4mn+n^2)^2$ and $z=-(m^2-4mn+n^2)^2$. There are also
points of order $4$ when $z=\pm (m^2+4mn+n^2)(m^2-4mn+n^2)$. The torsion subgroup is thus $\mathbb{Z}/2\mathbb{Z} \oplus \mathbb{Z}/4\mathbb{Z}$
or $\mathbb{Z}/2\mathbb{Z} \oplus \mathbb{Z}/8\mathbb{Z}$. Numerical testing suggests the former is always the case.

If we take the denominator of the f-transformation, solve for $w$, and substitute into the elliptic curve we find a point
$z=-T/(m^2+n^2)^2$ which gives a rational point on the curve. Adding to $(0,0)$ gives $z=-(m^2+n^2)^2$ which has $w=\pm 8mn(m^2-n^2)(m^2+n^2)$,
which suggests the curves have rank at least one.

As a numerical example, $m=5, n=2$, gives $z=-841, w= \pm 48720$. The positive value gives $f=-1$, but the negative value gives $f=18481/16481$.
This gives $p=11570426279$ and $q=224277439$ as possible values. From these, we can compute $a,b,c$ with
\begin{equation*}
a=66962532323709092281 \hspace{1cm} b=66912231954064693560
\end{equation*}
\begin{equation*}
 c=-66242061171706762440
\end{equation*}
with
\begin{equation*}
a+b=41^2*1231^2*229249^2 \hspace{2cm} a-b=23^2*3079^2*3167^2
\end{equation*}
\begin{equation*}
a+c=848805721^2 \hspace{2cm} a-c=41^2*103^2*2732993^2
\end{equation*}
\begin{equation*}
a+b+c=29^2*10427^2*27197^2 \hspace{2cm} a-b-c=11^2* 29^2* 109^2* 234161^2
\end{equation*}

Sadly, none of the other $3$ elements of the magic square are squares.

We can easily perform all these calculations in Pari, as well as considering the $15$ other distinct configurations of six
squares. Bremner's solution is found quickly in several configuration calculations, but no second solution has yet been discovered.

\section{Acknowledgement}
I would like to express my sincere thanks to Andrew Bremner for all his help, and his courtesy in answering lots of
questions.

\newpage

\appendix

\section{Basics of Elliptic Curves}
The standard elliptic curve over $\mathbb{Q}$ can be written
\begin{equation}\label{ecgenn}
y^2+a_1\,x\,y+a_3y=x^3+a_2x^2+a_4x+a_6
\end{equation}
with $a_i \in \mathbb{Z}$.

Firstly, we transform the equation of the curve to a simpler form.
From $(1.1)$ we have
\begin{equation}
4y^2+4a_1\,x\,y+4a_3y=4x^3+4a_2x^2+4a_4x+4a_6
\end{equation}
so
\begin{equation}
(2y+a_1x)^2+4a_3y=4x^3+(4a_2+a_1^2)x^2+4a_4x+4a_6
\end{equation}

Define $z=2y+a_1x$, so
\begin{equation}
z^2+2a_3(z-a_1x)=4x^3+(4a_2+a_1^2)x^2+4a_4x+4a_6
\end{equation}
giving
\begin{equation*}
z^2+2a_3z=4x^3+(4a_2+a_1^2)x^2+(4a_4+2a_1a_3)x+4a_6
\end{equation*}
so that
\begin{equation*}
(z+a_3)^2=4x^3+(4a_2+a_1^2)x^2+(4a_4+2a_1a_3)x+4a_6+a_3^2
\end{equation*}

Now, define $w=z+a_3$, giving
\begin{equation}
w^2=4x^3+(4a_2+a_1^2)x^2+(4a_4+2a_1a_3)x+4a_6+a_3^2
\end{equation}
and
\begin{equation}
16w^2=64x^3+(4a_2+a_1^2)16x^2+(4a_4+2a_1a_3)16x+64a_6+16a_3^2
\end{equation}

Finally, set $4w=v$ and $4x=u$ leading to
\begin{equation}
v^2=u^3+(4a_2+a_1^2)u^2+4(4a_4+2a_1a_3)u+16(4a_6+a_3^2)
\end{equation}
which we write as
\begin{equation}\label{ecgen}
E:\,v^2=u^3+au^2+bu+c=f(u)
\end{equation}
with $a,b,c \in \mathbb{Z}$. This is the form we use for all the elliptic curves in this survey. It has the great advantage of
symmetry about the u-axis, so, if $(p,q)$ is a point on the curve, so is $(p,-q)$.

We define addition of rational points in the usual chord-and-tangent way, see Silverman and Tate's
\cite{siltate} excellent introduction to elliptic curves.
If $p \ne r$ the line joining the points has equation
\begin{equation*}
v=\frac{s-q}{r-p} (u-p)+q
\end{equation*}
and so substituting this into the equation \eqref{ecgen} gives a cubic in $u$ for the points where the line and curve meet.

The constant of the cubic is
\begin{equation*}
c(p-r)^2-p^2s^2+qr(2ps-qr)
\end{equation*}
which is rational , and also is $(-p\,r\,t)$ where $t$ is the u-coordinate of the third intersection of the line and the curve.
Thus $t \in \mathbb{Q}$ and the corresponding v-coordinate (call it $w$) is also rational.

Thus combining two rational points in this way gives a third rational point. This gives a closed binary operation on the set
of rational points, which is often denoted by $\Gamma$. This is the first requirement for defining a group.

To get the associative property needed for a group, we actually need to take the result of the operation to be the point
$(t,-w)$ and not the actual point of intersection $(t,w)$.

We define the "addition" of rational points as
\begin{equation}
(p,q) \oplus (r,s) = (t,-w)
\end{equation}

With this binary operation, the set of rational points can be shown to be an Abelian group. The identity (denoted by $O$) is the
point at infinity, which is best approached by projectivising the curve, setting $u=U/W$ and
$v=V/W$ to give
\begin{equation}
V^2\,W=U^3+a\,U^2\,W+b\,U\,W^2+c\,W^3
\end{equation}
and the point at infinity is $(0,1,0)$. This gives the inverse of $(p,q)$ as $(p,-q)$.

$\Gamma$ is a finitely generated group, and is
isomorphic to $T_G \oplus \mathbb{Z}^r$, where $r$ is called the rank of the curve and $T_G$ is called the {\bf torsion subgroup} and is one of
\begin{itemize}
\item[{(i)}] $\hspace{1cm} \mathbb{Z}/n\mathbb{Z}, \hspace{0.5cm} n=1,2,\ldots,10,12$,
\item[{(ii)}] $\hspace{1cm} \mathbb{Z}/2\mathbb{Z} \oplus \mathbb{Z}/2n\mathbb{Z}, \hspace{0.5cm} n=1,2,3,4$.
\end{itemize}
as proven by Mazur \cite{mazur}, though first conjectured by Beppo Levi much earlier in the 20th century.

This means that there exist $r$ rational points $P_1,\ldots,P_r$ such that all rational points can be written
\begin{equation}
(u,v)=n_1P_1+\ldots+n_rP_r+T_i
\end{equation}
where $n_i \in \mathbb{Z}$ and $T_i \in T_G$. The hard part in all the problems is finding $P_1, \ldots , P_r$.

If $P$ is a point of order $2$ then $P+P=O$, so $P=-P$. The only way this can happen is if the v-coordinate of $P$ is zero,
so that the u-coordinate is a zero of $f(u)$. We want this root to be rational, but, since the coefficient of $v^3$ is $1$,
any rational root will actually be an integer root. Thus point of order $2$ in $\Gamma$ can only come from integer
roots of $f(u)=0$.

There is a theorem of Nagell and Lutz, which states that for curves of the form \eqref{ecgen}, all rational torsion points
have integer coordinates. There are two important points to note. This result does not always hold for curves of the form
\eqref{ecgenn}, and points of infinite order can have integer coordinates.

If $\Gamma$ has a point of order $2$, at $u=p$ say, if we define $w=u-p$ the elliptic curve reduces to the very simple but important form
\begin{equation}\label{ector2}
v^2=w^3+Aw^2+Bw
\end{equation}
with $A,B \in \mathbb{Z}$. Ten of the fifteen different torsion subgroup types have points of order $2$ so this form is a basic
structure.

Suppose $P=(p,q) \in \Gamma$ for a curve of this form. The point $2P=P+P$ comes from the intersection of the tangent at $P$ with the curve.
Some simple algebra shows that the w-coordinate of $2P$ is thus
\begin{equation}
\frac{(p^2-B)^2}{4q^2}
\end{equation}
so must be a rational square and, by Nagell-Lutz, must be an integer square. This is an incredibly useful result.

\newpage

\section{Transforming Quartics to Elliptic Curves}
We now show how to transform the quartic
\begin{equation}
y^2=ax^4+bx^3+cx^2+dx+e ,
\end{equation}
with $a,b,c,d,e \in \mathbb{Q}$, into an equivalent elliptic curve.

We first consider the case when $a=1$, which is very common. If $a=\alpha^2$ for some rational $\alpha$, we
substitute $y=Y/\alpha$ and $x=X/\alpha$, giving
\begin{equation*}
Y^2=X^4+\frac{b}{\alpha}X^3+cX^2+\alpha d X+\alpha^2 e
\end{equation*} Thus, suppose
\begin{equation}
y^2=x^4+bx^3+cx^2+dx+e
\end{equation}

We describe the method given by Mordell on page $77$ of \cite{mord}, with some minor modifications.

We first get rid of the cubic term by making the standard substitution $x=z-b/4$ giving
\begin{equation*}
y^2=z^4+fz^2+gz+h
\end{equation*}
where
\begin{equation*}
f=\frac{8c-3b^2}{8} \hspace{1cm} g=\frac{b^3-4bc+8d}{8} \hspace{1cm}
h=\frac{-(3b^4-16b^2c+64bd-256e)}{256}
\end{equation*}

We now get rid off the quartic term by defining $y=z^2+u+k$, where $u$ is a new variable and $k$ is a constant to be
determined. This gives the quadratic in $z$
\begin{equation*}
(f-2(k+u))z^2+gz+h-k^2-u(2k+u)=0
\end{equation*}

If $x$ were rational, then $z$ would be rational, so the discriminant of this quadratic would have to be a rational square.
The discriminant is a cubic in $u$. We do not get a term in $u^2$ if we make $k=f/6$, giving
\begin{equation*}
D^2=-8u^3+2\frac{f^2+12h}{3}u+\frac{2f^3-72fh+27g^2}{27}
\end{equation*}
and, if we substitute the formulae for $f,g,h$, and clear denominators we have
\begin{equation}
G^2=H^3+27(3bd-c^2-12e)H+27(27b^2e-9bcd+2c^3-72ce+27d^2)
\end{equation}
with
\begin{equation}\label{q4xtr}
x=\frac{2G-3bH+9(bc-6d)}{12H-9(3b^2-8c)}
\end{equation}
and
\begin{equation}
y= \pm \frac{18x^2+9bx+3c-H}{18}
\end{equation}

If we set the denominator of \eqref{q4xtr} to zero, we have $H=3(3b^2-8c)/4$ which gives
\begin{equation*}
G= \pm \frac{27(b^3-4bc+8d)}{8}
\end{equation*}
automatically giving a point on the curve. Often, it will just be a torsion point, but, occasionally, we get
the coordinates of a point of infinite order. Note also, that if $b=d=0$, then $G=0$ so we always have a point
of order $2$ from $y^2=x^4+cx^2+e$.

In this latter case, if we transfer the origin to this point of order $2$ and simplify we get
\begin{equation*}
S^2=T^3-2cT^2+(c^2-4e)T
\end{equation*}
with
\begin{equation*}
x=\frac{S}{2T} \hspace{2cm} y=\frac{2x^2+c-T}{2}
\end{equation*}

\bigskip

If $a \ne \alpha^2$, we need a rational point $(p,q)$ lying on the quartic.

Let $z=1/(x-p)$, so that $x=p+1/z$ giving
\begin{equation*}
y^2z^4=(ap^4+bp^3+cp^2+dp+e)z^4+(4ap^3+3bp^2+2cp+d)z^3+
\end{equation*}
\begin{equation*}
(6ap^2+3bp+c)z^2+(4ap+b)z+a
\end{equation*}
where the coefficient of $z^4$ is, in fact, $q^2$. Define, $y=w/z^2$, and then
$z=u/q^2, \, w=v/q^3$
giving
\begin{equation}
v^2=u^4+(4ap^3+3bp^2+2cp+d)u^3+q^2(6ap^2+3bp+c)u^2+
\end{equation}
\begin{equation*}
q^4(4ap+b)u+aq^6 \equiv u^4+fu^3+gu^2+hu+k
\end{equation*}

We now, essentially, complete the square. We can write
\begin{equation*}
y^2=u^4+fu^3+gu^2+hu+k=(u^2+mu+n)^2+(su+t)
\end{equation*}
if we set
\begin{equation*}
m=\frac{f}{2} \hspace{1cm} n=\frac{4g-f^2}{8} \hspace{1cm} s=\frac{f^3-4fg+8h}{8}
\end{equation*}
and
\begin{equation*}
 t=\frac{-(f^4-8f^2g+16(g^2-4k))}{64}
\end{equation*}

This gives
\begin{equation*}
(y+u^2+mu+n)(y-u^2-mu-n)=su+t
\end{equation*}
and if we define $y+u^2+mu+n=Z$ we have
\begin{equation}
2(u^2+mu+n)=Z-\frac{su+t}{Z}
\end{equation}

Multiply both sides by $Z^2$, giving
\begin{equation*}
2u^2Z^2+2muZ^2+suZ=Z^3-2nZ^2-tZ
\end{equation*}
which, on defining $W=uZ$, gives
\begin{equation*}
2W^2+2mWZ+sW=Z^3-2nZ^2-tZ
\end{equation*}

Define, $Z=X/2$ and $W=Y/4$ giving
\begin{equation}
Y^2+2mXY+2sY=X^3-4nX^2-4tX
\end{equation}
which is an elliptic curve.

If we define $Y=G-s-mX$, we transform to the form
\begin{equation}
G^2=X^3+(m^2-4n)X^2+(2ms-4t)X+s^2
\end{equation}

All of the above are easy to program in any symbolic package.

A variant of the above problem, that occurs often enough to be important, is when $b=d=0$, so the quartic is
\begin{equation*}
y^2=ax^4+cx^2+e
\end{equation*}

Applying the above formulae gives an elliptic curve where the cubic right-hand-side has a linear factor, so the curve
has a point of order $2$. Transforming to make this point the origin, and doing a variable scaling gives
\begin{equation}
F^2=H^3-2cH^2+(c^2-4ae)H=H(H^2-2cH+(c^2-4ae))
\end{equation}

\newpage

\section{Searching for rational points on $y^2=x^3+ax^2+bx$}
In this section, we describe a very simple method for searching for rational points on
\begin{equation}\label{ecdesc}
y^2=x^3+ax^2+bx
\end{equation}
where $a,b \in \mathbb{Z}$.

Suppose $(p,q) \in \Gamma$ where we set $p=N/D \, , \, q=M/E$, with $\gcd(N,D)=1$ and $\gcd(M,E)=1$, and $N, D, M, E \in \mathbb{Z}$.
Then
\begin{equation*}
M^2D^3=E^2(N^3+aN^2D+bND^2)
\end{equation*}

It is clear that $D=1$ iff $E=1$, so suppose both $D,E >1$.

Now $D$ cannot divide $N^3+aN^2D+bND^2$ since, if it did, that would give $D \mid N$ against $\gcd(N,D)=1$. Thus $D^3 \, \mid \, E^2$.
Similarly, we must have $E^2 \, \mid \, D^3$, so $D^3=E^2$ which implies there must exist an integer $v$ with $D=v^2$ and $E=v^3$.

Thus $x=N/v^2$ and we can always write $N$ as $d\,u^2$ where $d$ is squarefree. Putting $x=du^2/v^2$ and $y=M/v^3$ into \eqref{ecdesc}
gives
\begin{equation}
M^2=du^2(d^2u^4+adu^2v^2+bv^4)
\end{equation}
so that $d \mid M^2$ and $u^2 \mid M^2$. Thus $u \mid M$ and, since $d$ is squarefree by definition, $d \mid M$. Thus there must exist
an integer $w$ such that $y=duw/v^3$.

Thus
\begin{equation}
d^2w^2=d^3u^4+ad^2u^2v^2+bdv^4
\end{equation}
so that $d \mid bv^4$, and, since $\gcd(d,v)=1$ we must have $d \mid b$. So we can determine possible values for $d$ from the squarefree
factors of $b$. If $a^2<4b$, there will be no negative $d$ values.

Now,
\begin{equation*}
4dw^2=4d^2u^2+4adu^2v^2+4bv^4=(2du^2+av^2)^2-(a^2-4b)v^4
\end{equation*}
and we can easily write $a^2-4b=\alpha\beta^2$ where $\alpha$ is squarefree.

Define $F=2du^2+av^2$, $G=\beta v^2$ and $H=2w$, so that
\begin{equation}\label{desc1}
F^2=\alpha G^2+dH^2
\end{equation}
and
\begin{equation}\label{desc2}
\frac{u^2}{v^2}=\frac{\beta F-aG}{2dG}
\end{equation}

Thus we can look for solutions $(F,G,H)$ of \eqref{desc1} and substitute into the right-hand-side of \eqref{desc2} to see if we
get a rational square.

Quite often $\alpha=1$ and we can parameterize $F^2=G^2+dH^2$ as $F=p^2+dq^2$, $G=p^2-dq^2$ and $H=2pq$, which makes for a very
simple program.

If $\alpha \ne 1$, we can find an initial non-trivial solution $(F_0,G_0,H_0)$ of \eqref{desc1} and use standard quadric parameterisation
to give
\begin{equation*}
\frac{F}{G}=\frac{d(F_0-2mH_0)+m^2F_0}{G_0(d-m^2)}
\end{equation*}
where $m \in \mathbb{Q}$.

This gives
\begin{equation}
\frac{u^2}{v^2}=\frac{(aF_0+\beta G_0)m^2-2\beta dH_0m+d(\beta F_0-aG_0)}{2dG_0(d-m^2)}
\end{equation}

This simple method provides lots of points of small heights.

We can develop further by taking the numerator and denominator of this equation separately as
\begin{equation*}
-2dG_0m^2+2d^2G_0=ks^2
\end{equation*}
\begin{equation*}
(aF_0+\beta G_0)m^2-2\beta dH_0m+d(\beta F_0-aG_0)=kr^2
\end{equation*}
where $k$ is squarefree. We try to satisfy both simultaneously.

Values of $k$ must divide the resultant of the two quadratics on the left-hand-sides of these two relations, so
\begin{equation*}
k \mid 16\, d^4 \, G_0^3 \, (a^2-4b)
\end{equation*}
so there are only a computable finite set of possible values.

Suppose we find an initial solution $(m_0,s_0)$ to the first quadratic, then the line $s=s_0+t(m-m_0)$ meets it again where
\begin{equation*}
m=\frac{-(2dG_0m_0+kt(2s_0-m_0t))}{2dG_0+kt^2}
\end{equation*}

Substituting into the second quadratic, we must have
\begin{equation}\label{descq4}
k(c_0t^4+c_1t^3+c_2t^2+c_3t+c_4)= \Box
\end{equation}
with
\begin{equation*}
c_0=k^2(\beta(d(F_0-2H_0m_0)+G_0m_0^2)-a(dG_0-F_0m_0^2))
\end{equation*}
\begin{equation*}
c_1=-4k^2s_0(aF_0m_0+\beta(G_0m_0-dH_0))
\end{equation*}
\begin{equation*}
c_2=4k(\beta G_0(d^2F_0-dG_0m_0^2+ks_0^2)-a(d^2G_0^2+dF_0G_0m_0^2-F_0ks_0^2))
\end{equation*}
\begin{equation*}
c_3=8dG_0ks_0(aF_0m_0+\beta(dH_0+G_0m_0))
\end{equation*}
\begin{equation*}
c_4=4d^2G_0^2(\beta(d(F_0+2H_0m_0)+G_0m_0^2)-a(dG_0-F_0m_0^2))
\end{equation*}

Finding rational solutions to \eqref{descq4} is a non-trivial problem. There is not enough space to describe the many
tricks of the trade that should be employed to do this efficiently.

This entire method has been the workhorse for finding hundreds of solutions.

\newpage

\section{Isogenies of Elliptic Curves}
Isogenies are a crucial tool in finding points with large heights by allowing us to consider curves which {\bf might} have
rational points with much smaller heights. The classic reference is the paper of V{\'e}lu \cite{velu}.

We simplify the notation to considering curves of the form
\begin{equation}\label{eci1}
y^2=x^3+ax^2+bx+c
\end{equation}

In this case, points of order $2$ have coordinates $(p,0)$, and other torsion points occur in pairs $(p,\pm q)$. If
$T$ denotes the torsion subgroup, let $T^* \subset T$ be the subset consisting of all points of order $2$ together with {\bf ONE}
point from each of the pairs $(p, \pm q)$.

Define the two functions
\begin{equation}
g(x)=3x^2+2ax+b \hspace{2cm} h(y)=-2y
\end{equation}

For each point $P=(r,s) \in T^*$, define
\begin{equation}
t_P=\delta g(r) \hspace{1cm} u_P=h^2(s) \hspace{1cm} w_P=u_P+r\,t_P
\end{equation}
with $\delta=1$ if $P$ is of order $2$ and $\delta=2$ otherwise.

Let
\begin{equation}
t=\sum_{P \in T^*} t_P \hspace{2cm} w=\sum_{P \in T^*} w_P
\end{equation}
then the elliptic curve \eqref{eci1} is isogenous to
\begin{equation}\label{eci2}
v^2=u^3+au^2+(b-5t)u+(c-4at-7w)
\end{equation}

The relations between the variables are
\begin{equation}\label{eqisu}
u=x+\sum_{P \in T^*} \left( \frac{t_P}{x-r_P} +\frac{u_P}{(x-r_P)^2} \right)
\end{equation}
and
\begin{equation}\label{eqisv}
v=y-\sum_{P \in T^*} \left( \frac{2u_Py}{(x-r_P)^3} + \frac{t_P(y-s_P)}{(x-r_P)^2}-\frac{g(r_P)h(s_P)}{(x-r_P)^2} \right)
\end{equation}

Consider, first, the curve
\begin{equation}
y^2=x^3+ax^2+bx
\end{equation}
with just the torsion point $(0,0)$. Then $t_P=b=t$, $u_P=0$ and $w_P=0=w$. Thus equation \eqref{eci2} is
\begin{equation*}
v^2=u^3+au^2-4bu-4ab=(u+a)(u^2-4b)
\end{equation*}
and setting $u+a=h$ gives the standard 2-isogenous curve
\begin{equation}
v^2=h^3-2ah^2+(a^2-4b)h
\end{equation}

Equations \eqref{eqisu} and \eqref{eqisv} give
\begin{equation*}
h=\frac{x^2+ax+b}{x}=\frac{y^2}{x^2} \hspace{2cm} v=\frac{y(x^2-b)}{x^2}
\end{equation*}

Suppose we find a point on the 2-isogenous curve $P=(h_1,v_1)$. Solving the left-hand equation for $h=h_1$ will
usually give a non-rational value for $x$. If, however, we form $2P$ we find it has h-coordinate
\begin{equation*}
h=\frac{(h_1^2-(a^2-4b))^2}{4v_1^2}
\end{equation*}
Putting this into the left-hand equation we find $x=v_1^2/(4h_1^2)$.

\hspace{1cm}

Next, consider
\begin{equation}
y^2=x^3+(ax+b)^2
\end{equation}
which has $2$ points of order $3$ at $(0,\pm b)$.

Then, $T^*=\{(0,b)\}$ say and $t_P=4ab=t$, $u_P=4b^2$ and $w_P=4b^2=w$, giving
\begin{equation*}
v^2=u^3+a^2u^2-18abu-27b^2-16a^3b
\end{equation*}

For a curve of the form $y^2=f(u)$, a point of inflexion would occur when $2f\,f^{\prime\prime}-(f^{\prime})^2=0$.
This gives
\begin{equation*}
(u^3-36abu-4b(4a^3+27b))(3u+4a^2)=0
\end{equation*}
but there is not a point of inflexion at $u=-4a^2/3$ since $f$ is negative at this point. We can, however, get a
simpler form for the 3-isogenous curve by using this point.

Defining $v=g/27$ and $u=(h-12a^2)/9$ gives
\begin{equation}
g^2=h^3-27(ah+27b-4a^3)^2
\end{equation}
with
\begin{equation*}
h=\frac{3(3x^3+4a^2x^2+12abx+12b^2)}{x^2} \hspace{1.5cm} g=\frac{27y(x^3-4abx-8b^2)}{x^3}
\end{equation*}

If we find a point $P=(h_1,g_1)$ on the 3-isogenous curve, we first find $3P=(h_3,g_3)$ and solve the first
of these equations with $h=h_3$ to give the corresponding x-value.

\hspace{1cm}

As a final example, consider
\begin{equation}
y^2=x^3+(a^2+2b)x^2+b^2x
\end{equation}
where $a,b \in \mathbb{Z}, ab \ne 0$, which has $(0,0)$ as a point of order $2$, and $(-b,\pm ab)$ of order $4$.

Let $T^*=\{(0,0),(-b,ab)\}$, which give $t=b(b-4a^2)$ and $w=8a^2b^2$, so the 4-isogenous curve is
\begin{equation*}
v^2=u^3+(a^2+2b)u^2+4b(5a^2-b)u+4b(4a^4-7a^2b-2b^2)
\end{equation*}
which can be factored to
\begin{equation*}
v^2=(u+a^2-2b)(u^2+4bu+4b(4a^2+b))
\end{equation*}

Define $u=h+2b-a^2$ and we arrive at
\begin{equation}
v^2=h^3+2(4b-a^2)h^2+(a^2+4b)^2h
\end{equation}
with
\begin{equation}\label{4ish}
h=\frac{(x-b)^2(x^2+(a^2+2b)x+b^2)}{x(x+b)^2}
\end{equation}

\begin{equation*}
v=\frac{y(x-b)(x^4+4bx^3+2b(2a^2+3b)x^2+4b^3x+b^4)}{x^2(x+b)^3}
\end{equation*}

If we find a point $P=(h_1,g_1)$ on the 4-isogenous curve, we first find $4P=(h_4,g_4)$ and solve \eqref{4ish}
with $h=h_4$ to give the corresponding x-value.

\end{document}